\documentclass[preprint,times]{elsarticle}
\makeatletter
\def\ps@pprintTitle{%
  \let\@oddhead\@empty
  \let\@evenhead\@empty
  \let\@oddfoot\@empty
  \let\@evenfoot\@oddfoot
}
\makeatother
\usepackage{psfig}
\usepackage{amsmath}
\usepackage{geometry}                
\usepackage{graphicx}
\usepackage{amssymb}
\usepackage{epstopdf}
\usepackage{amsmath,amssymb,amsthm}
\usepackage {natbib}
\usepackage{bm}
\usepackage{bbm}
\usepackage{mathtools}
\usepackage{gensymb}
\usepackage{algorithm}
\usepackage{algorithmic}
\usepackage{setspace}
\usepackage{xcolor}

\usepackage{stmaryrd}

\usepackage{epsfig}
\usepackage{url}

\theoremstyle{example} \theoremstyle{assumption}
\theoremstyle{remark} \theoremstyle{lemma}
\theoremstyle{definition} \theoremstyle{corol}
\theoremstyle{proposition} \theoremstyle{condition}

\newtheorem{problem}{\n{Problem}}[section]

\newtheorem{lemma}{\n{Lemma}}[section]
\newtheorem{definition}{\n{Definition}}[section]

\newtheorem{proposition}{\n{Proposition}}[section]

\font\n=cmcsc10

\DeclareMathOperator*{\esssup}{ess\,sup}
\DeclareMathOperator\supp{supp}
\definecolor{tealblue}{rgb}{0.21, 0.46, 0.53}
   
\newcommand{\norm}[1]{\left\lVert #1 \right\rVert}
\newcommand{\jump}[1]{\Bigl\llbracket #1 \Bigr\rrbracket}
\newcommand{\abs}[1]{\left\lvert #1 \right\rvert}

\RequirePackage{fix-cm} 

\DeclareMathSizes{10}{10}{6}{4}

\begin{document}

\begin{frontmatter}

\title{{A} posteriori error estimator for elliptic interface problems in the fictitious formulation}

%

\author[1]{Najwa Alshehri}
\ead{najwa.alshehri@kaust.edu.sa}
\address[1]{King Abdullah University for Science and Technology, K.S.A.}

\author[1,2]{Daniele Boffi}
\ead{daniele.boffi@kaust.edu.sa}
\address[2]{Dipartimento di Matematica ``F. Casorati'', Universit\`a di
Pavia, Italy.}

\author[3]{Lucia Gastaldi}
\ead{lucia.gastaldi@unibs.it}
\address[3]{DICATAM, Universit\`a di Brescia, Italy.}

\begin{abstract}
A posteriori error estimator is derived for an elliptic interface problem in the fictitious domain formulation with distributed Lagrange multiplier considering a discontinuous Lagrange multiplier finite element space. A posteriori error estimation plays a pivotal role in assessing the accuracy and reliability of computational solutions across various domains of science and engineering.  This study delves into the theoretical underpinnings and computational considerations of a residual-based estimator.\\
Theoretically, the estimator is studied for cases with constant coefficients which jump across an interface as well as generalized scenarios with smooth coefficients that jump across an interface. Theoretical findings demonstrate the reliability and efficiency of the proposed estimators under all considered cases.\\
Numerical experiments are conducted to validate the theoretical results, incorporating various immersed geometries and instances of high coefficients jumps at the interface. Leveraging an adaptive algorithm, the estimator identifies regions with singularities and applies refinement accordingly. Results substantiate the theoretical findings, highlighting the reliability and efficiency of the estimators. Furthermore, numerical solutions exhibit optimal convergence properties, demonstrating resilience against geometric singularities or coefficients jumps.
\end{abstract}

\begin{keyword}
interface elliptic problems\sep finite elements method\sep unfitted approach\sep
Lagrange multiplier\sep a posteriori analysis\sep immersed boundary method.
\end{keyword}

\end{frontmatter}
\section{Introduction}
Elliptic interface problems are a class of partial differential equations (PDEs) that arise in a wide range of applications, such as fluid dynamics, heat transfer, and structural mechanics. These problems are characterized by discontinuous coefficients or interfaces in the solution domain, which can pose significant challenges in terms of numerical solution and error estimation. The a priori estimates are significantly influenced by these singularities, making it challenging to achieve accurate numerical approximations. Additionally, the complexity of interface geometries further make the numerical approximation more challenging, underscoring the need for improved accuracy in computational schemes.

A powerful technique for assessing the accuracy of numerical solutions to partial differential equations is a posteriori error estimation. This method involves computing an error indicator that quantifies the discrepancy between the numerical and exact solutions, subsequently driving adaptive mesh refinement to enhance solution accuracy.

Numerous a posteriori error estimators have been proposed in the literature, spanning various methodologies such as residual-based approaches~\cite{babuvvska1978error}, \cite{demkowicz1984adaptive}, and \cite{demkowicz1985adaptive}, hierarchical basis error techniques~\cite{bank1985some}, solution of auxiliary local problems~\cite{lonsing2004posteriori}, and recovery-based estimators~\cite{zienkiewicz1992superconvergent1}, and \cite{zienkiewicz1992superconvergent2}. Verf\"{u}rth provides a comprehensive overview of a posteriori error estimators in \cite{verfurth1996review}, and in~\cite{verfurth2013posteriori}.

Over the years, a great deal of research has focused on a posteriori analysis for interface problems. In 2000, Bernardi and Verf\"{u}rth established a residual-based estimator for conforming finite elements for second-order elliptic problems with discontinuous or anisotropic coefficients in \cite{bernardi2000adaptive}. The estimator was proven to be robust under the assumption that no more than three subdomains share a common point on the distribution of the coefficients. 

In 2002, Petzoldt derived an a posteriori error estimator for linear finite elements in \cite{petzoldt2002posteriori}. This estimator provided global upper and local lower bounds that were independent of the jump of coefficients at the interface, assuming quasi-monotonic distribution of coefficients. In the same year, Chen and Dai, in \cite{chen2002efficiency}, studied the influence of the factors that affect the performance of the adaptive methods such as the a posteriori error indicators, the refinement strategies and the choice of parameters.

Since then, many other robust estimators have been established for interface problems \cite{luce2004local}, \cite{berrone2006robust}, \cite{vohralik2011guaranteed}, and \cite{zhao2015robust}. For example, Cai, Zhang and their collaborators introduced a residual-based estimator for non-conforming finite elements in \cite{cai2017residual}, and for discontinuous Galerkin elements in \cite{cai2011discontinuous}. In \cite{cai2010recovery}, a recovery-based estimator was presented for mixed nonconforming linear finite elements, and in \cite{cai2009recovery}, a recovery-based estimator was proposed for conforming linear finite elements.

Among those approaches, the residual-based error estimator is widely used for a posteriori error estimation in the context of elliptic interface problems. This approach is based on the residual equation, which is obtained by subtracting the discrete equation satisfied by the finite element solution from the original elliptic equation. The residual equation represents the error in the discrete solution, and can be used to derive a reliable and efficient estimate of the error. 

Methods for solving elliptic interface problems can be categorized to fitted boundary approach, where the mesh is forced to align with the interface, or non fitted boundary approach, where the interface is allowed to cut the interior of the mesh cells. The latter approach overcome the challenges of having degenerated mesh with time, which is easier and computationally cheaper. The focus of this paper lies within the realm of Fictitious Domain formulations with Distributed Lagrange Multipliers (FD-DLM), representing an unfitted boundary approach to elliptic interface problems. This approach, emerged from the idea of utilizing finite elements immersed boundary methods, see \cite{boffi2003finite}. FD-DLM offers promising methods for addressing interface challenges without necessitating mesh realignment. We refer the reader to articles~\cite{boffi2021existence}, and \cite{boffi2022existence} that studied in depth this formulation.

We acknowledge that there has been significant studies and developments in the a priori error estimation for elliptic interface problems within the FD-DLM formulation over the past years~ \cite{auricchio2015fictitious}, \cite{boffi2014mixed}, and \cite{Najwa2022elliptic}, to name some. Due to the singularity of the solution, the available schemes in literature were able to achieve the optimal convergence rate of the error of the solution with respect to its regularity. 

We aim at provide reliable and efficient tools for assessing the accuracy of numerical solutions and achieve optimal rate of convergence regardless of the singularities in the solution or the complexity of the geometry by studying the a posteriori error estimation of the elliptic interface problems in the FD-DLM formulation. To our knowledge, no estimator has been studied for such a problem in this formulation.

In this paper, we propose a residual-based a posteriori error estimator for an alliptic interface problem in the FD-DLM formulation in order to improve control over the error, the accuracy of our model, and the computational cost. In the presence of interfaces, the residual equation can be particularly challenging to handle, due to the discontinuities in the solution and its derivatives across the interfaces. Addressing this challenge is the primary focus of our work.

This paper is organized as follows: Section~\ref{se:problem_setting} introduces the problem at both continuous and discrete levels, along with pertinent a priori error estimation results studied in \cite{Najwa2022elliptic}. Section~\ref{se:aposteriori} presents our proposed estimators and provides proofs for the Global Upper Bound (GUB) and the local lower bound (LLB).  The GUB proof utilize  the Cl\' ement interpolation operator and the LLB make use of the local bubble functions that are polynomials defined properly in order to localize the contribution of the error. Within the proof, we also discuss challenges posed by interface presence and review techniques to mitigate them as needed. In Section~\ref{se:Adaptivity}, we elucidate our approach to leveraging the a posteriori error estimator in adaptive schemes. Finally, Section~\ref{se:numerical_results} presents numerical tests considering various immersed geometries and factors, validating the reliability and efficacy of our proposed method.
\section{Problem setting}\label{se:problem_setting}
In this section, we establish the groundwork for our analysis by outlining the fundamental framework, highlighting the key aspects of the problem under consideration. We begin by introducing the notation employed throughout our investigation, providing a concise review of the symbols and conventions utilized. Subsequently, we delve into the core of our study by presenting the model problem, which serves as the basis for our analysis. This model problem encapsulates the essential characteristics of the phenomena under examination. Following the presentation of the model problem, we transition to the continuous Fictitious Domain Approach with a Distributed Lagrange Multiplier (FD-DLM) problem. Finally, we proceed to the discrete FD-DLM problem, where we discretize the continuous problem to obtain a unique numerical solution following the work in \cite{Najwa2022elliptic}.
\subsection{Notation}
\label{se:notation}
In this paper, we utilize various notations, which we shall now review. Let $\omega \subset \mathbb{R}^d$ for $d=2,3$. For a non-negative integer $s$, we denote the standard Sobolev space by $ H^s(\omega) $, which is the space of $ s $-times differentiable and square-integrable functions $\phi$ whose derivatives up to order \(s\) are square-integrable, equipped with the norm $\norm{\phi}_{s,\omega}$ defined as follows:
\begin{align*}
\norm{\phi}_{s,\omega} &= \left[ \sum_{\abs{\alpha} \leq s} \int_{\omega} \abs{\nabla^\alpha \phi}^2 \right]^{\frac{1}{2}},\\
\abs{\phi}_{s,\omega} &= \left[ \sum_{\abs{\alpha} = s} \int_{\omega} \abs{\nabla^\alpha \phi}^2 \right]^{\frac{1}{2}}.
\end{align*}

Let \( H^0(\omega) \coloneqq L^2(\omega) \) be the space of square-integrable
functions on \(\omega\), equipped with the norm \(\norm{\phi}_{0,\omega}\).
The standard scalar product associated with this space is represented by
\(\left(\cdot, \cdot \right)_{\omega}\). Additionally, \( H^1(\omega) \)
denotes the Hilbert space of functions in \(\omega\) that are
square-integrable and possess square-integrable first derivatives, equipped
with its norm \(\norm{\phi}_{1,\omega}\). A particular subset of this space is
\( H^1_0(\omega) \), which denotes the subspace of  \( H^1(\omega) \)
containing functions that vanish at the boundary, endowed with the norm
\(\abs{\phi}_{1,\omega}\) which is a norm thanks to the Poincar\'{e}
inequality. Moreover, \([H^1(\omega)]^* \) denotes the dual space of \( H^1(\omega) \), associated with its dual norm $\norm{\mu}_{[H^1(\omega)]^*}$. This norm is
defined, for \(\mu \in [H^1(\omega)]^*\), as follows:
\[
\norm{\mu}_{[H^1(\omega)]^*} = \sup_{\phi \in H^1(\omega)} \dfrac{\left\langle \mu, \phi \right\rangle }{\norm{\phi}_{1,\omega}}.
\]

where $\left\langle \mu, \phi \right\rangle$ is the duality pairing between $H^1(\omega)$ and its dual space $[H^1(\omega)]^*$.
Similarly, we define the dual space $[H^1_0(\omega)]^*$ and we denote it by $H^{-1}(\omega)$. Clearly we have
\begin{equation}\label{eq:norm-1}
\norm{\mu}_{H^{-1}(\omega)}\leq \norm{\mu}_{[H^1(\omega)]^*}.
\end{equation} 
Furthermore, for $\omega_1 \subset \omega$, we have
\begin{equation}\label{eq:norm1*}
\norm{\mu}_{[H^1(\omega_1)]^*}\leq \norm{\mu}_{[H^1(\omega)]^*}.
\end{equation}
\subsection{Model Problem}
\label{se:Model_problem}
Let $\Omega$ be a domain in $\mathbb{R}^d \; (d=2,3)$, with a Lipschitz continuous boundary $\partial \Omega$. The domain $\Omega$ is separated by a Lipschitz interface $\Gamma$ into two subdomains $\Omega_1$ and $\Omega_2$, where $\Omega_2$ is immersed so that the boundary of $ \Omega $ coincide with the boundary of $ \Omega_1$ as the case in Figure~\ref{fig:IBM_FDDLM1}. Assume that $\beta_1 \in W^{1,\infty} (\Omega_1)$ and $\beta_2 \in W^{1,\infty} (\Omega_2)$ are bounded by positive constants from below as follows
\begin{align*}
\beta_1 &> \underline{\beta_1}>0, &\beta_2 &> \underline{\beta_2}>0.
\end{align*}
\begin{figure}[htp]
\begin{minipage}[c]{0.45\linewidth}
\centering
		\includegraphics[width=0.4\linewidth]{IBM_FDDLM1.pdf}
		\caption{Immersed boundary domains.}
		\label{fig:IBM_FDDLM1}
\end{minipage}
\begin{minipage}[c]{.45\linewidth}
\centering
		\includegraphics[width=0.4\linewidth]{IBM_FDDLM2.pdf} 
		\caption{Fictitiously extended domain with immersed interface.}
		\label{fig:IBM_FDDLM2}
\end{minipage}
\end{figure}
We consider the following elliptic interface problem.
\begin{problem}\label{pbm:0}
Given $(f_1,f_2):\Omega_1\times \Omega_2\rightarrow \mathbb{R}$, find $(u_1,u_2):\Omega_1\times \Omega_2 \rightarrow \mathbb{R}$ such that:
\begin{align*}
		- \nabla \cdot (\beta_1 \nabla u_1)&=f_1		
					&&		\textit{ in } \Omega_1\\
		- \nabla \cdot (\beta_2 \nabla u_2)&=f_2		
					&&		\textit{ in } \Omega_2\\
	   u_1-u_2&=0		
	  				&&		\textit{ on } \Gamma  \\
	  \beta_1 \nabla (u_1) \cdot \textbf{n}_1+\beta_2 \nabla (u_2) \cdot \textbf{n}_2&=0		
					&&		\textit{ on } \Gamma \\ 
	 u_1&=0		
	 				&&		\textit{ on } \partial\Omega_1, 
\end{align*}
where $\textbf{n}_1$ and $\textbf{n}_2$ are the exterior unit vectors normal to the interface $\Gamma$.
\end{problem}

Problem~\ref{pbm:0} is governed by two Poisson equations with two different coefficients $\beta_1$ and $\beta_2$ that reflect the diffusivity properties of the materials $\Omega_1$ and $\Omega_2$, respectively. This problem is subject to a homogeneous boundary condition on $\partial \Omega_1$ and two transmission conditions on $\Gamma$. The first transmission condition requires the solutions to coincide at the interface $\Gamma$ and the other allows for a jump of the normal derivatives across the interface $\Gamma$ that is inversely proportional to the ratio between $\beta_1$ and $\beta_2$ which enforces the continuity of the conormal derivatives there.

The a priori analysis of this problem has been explored previously in references \cite{babuvska1970finite}, \cite{bruce1974poisson}, and \cite{lemrabet1978interface}, among others. The regularity of the solutions $u_1$ and $u_2$ can be affected by the presence of re-entrant corners. Therefore, $u_1 \in H^{2-\epsilon}(\Omega_1)$ and $u_2 \in H^{2-\epsilon}(\Omega_2)$ , where $0 \leq \epsilon<\frac{1}{2}$. Here we consider the Fictitious Domain Approach with Distributed Lagrange Multiplier (FD-DLM), which belongs to the category of unfitted immersed boundary methods, see\cite{auricchio2015fictitious},\cite{boffi2014mixed}, and \cite{Najwa2022elliptic}. In the following section, we explain the formulation of the continuous problem.
%
\subsection{Continuous FD-DLM problem} \label{se:cont_pbm}

To reformulate Problem~\ref{pbm:0} using the FD-DLM approach, we start by extending $\Omega_1$ fictitiously to fill also the region occupied by $\Omega_2$ and denoting this extension by $\Omega$, as in Figure~\ref{fig:IBM_FDDLM2}. Moreover, we extend $\beta_1, f_1$ smoothly to $\Omega$ and we denote these extensions by $ \beta, f$ such that $\beta|_{\Omega_1} = \beta_1$, and  $f|_{\Omega_1} = f_1$. 
In addition, we extend $u_1$ to $\Omega$ and denote the extension by $u$ such that $u|_{\Omega_1} = u_1$, and we require that the solution $u$ restricted to the common domain $\Omega_2$ coincides with $u_2$; i.e. $u|_{\Omega_2} = u_2$. Therefore, a Lagrange multiplier is added to enforce such constraint. 

Let $\mathcal{H}\coloneqq V \times V_2 \times \Lambda$ denote the space of the solution triplet $(u,u_2,\lambda)$, where $V= H^1_0(\Omega)$ and $V_2=H^1(\Omega_2)$. Regarding the Lagrange multiplier space, there are two potential selections for $\Lambda$, as extensively discussed in \cite{boffi2022existence}. One option is to define $\Lambda$ as the dual space of $V_2$, while the other is to set $\Lambda= V_2$. Each choice allows for different families of discretizations. In view of a piecewise discontinuous discretization of the Lagrange multiplier, we focus on the first choice, whose stability and approximation properties have been studied in \cite{Najwa2022elliptic}.

As shown in \cite{auricchio2015fictitious}, the fictitious domain approach leads to the following equivalent formulation of Problem~\ref{pbm:0}.
\begin{problem}
\label{pbm:1}
Given $f \in L^2(\Omega)$, $f_2 \in L^2(\Omega_2)$, $\beta \in W^{1,\infty} (\Omega)$ and $\beta_2 \in W^{1,\infty} (\Omega_2)$ with $f|_{\Omega_1}=f_1$ and 
$\beta|_{\Omega_1}=\beta_1$, find
$(u,u_2,\lambda)\in \mathcal{H}$ such that
\begin{align*}
(\beta \nabla u,\nabla v)_{\Omega} + \left\langle \lambda,  v|_{\Omega_2} \right\rangle
&=(f,v)_{\Omega} && \forall v \in V\\
((\beta_2-\beta) \nabla u_2,\nabla v_{2})_{\Omega_2} -\left\langle \lambda, v_2 \right\rangle
&=(f_2-f,v_2)_{\Omega_2} && \forall v_2 \in V_2\\
\left\langle \mu , u|_{\Omega_2} -u_2 \right\rangle
&=0 && \forall \mu \in \Lambda.
\end{align*}
\end{problem}
One can obtain the following representation of
$\lambda$ from the above problem, which will be useful to estimate the error and study the regularity of the Lagrange multiplier,
( see~\cite{auricchio2015fictitious}). 
\begin{equation}
\label{eq:dual_lambda}
\left\langle \lambda, v_2 \right\rangle = -
\int_{\Omega_2}\left(\frac{\beta}{\beta_2} f_2-f \right)\; v_2\;dx
+\int_{\Gamma} (\beta_2 -\beta) \nabla u_2 \cdot \textbf{n}_2\; v_2\; d \gamma.
\end{equation}

The regularity of the fictitiously extended solution $u$ is affected by the jump of the coefficients across the interface $\Gamma$. Hence, $u \in H^{\frac{3}{2}-\epsilon}(\Omega)$, (see~\cite{nicaise}). From characterization of $\lambda$ in \eqref{eq:dual_lambda}, we can write $\lambda =\lambda_1+\lambda_2$, where clearly $\lambda_1 \in L^2(\Omega_2)$ and $\lambda_2 \in H^{\epsilon-1}(\Omega_2)$ with $0 \leq \epsilon < \frac{1}{2}$. We conclude $ \norm{\lambda_2}_{\epsilon-1 ,\Omega_2}\leq C \norm{u_2}_{2-\epsilon,\Omega_2},
$ which is bounded for $u_2\in H^{2-\epsilon}(\Omega_2)$, see \cite{Najwa2022elliptic} for the proof.

Problem~\ref{pbm:1} is characterized as a saddle point problem, whose stability was demonstrated in \cite{auricchio2015fictitious}. Furthermore, the existence of a unique solution was affirmed  by showing that the problem satisfies the following ellipticity on the kernel and inf-sup conditions.

There exist constants $C_1, C_2 >0$ such that:
\begin{align*}
(\beta \nabla u,\nabla u)_{\Omega}+((\beta_2-\beta) \nabla u_2,\nabla u_{2})_{\Omega_2} &\geq C_1 \left\lbrace \norm{u}_{V}^2+ \norm{u_2}_{V_2}^2\right\rbrace && \forall (u,u_2)  \textit{ such that }\left\langle \mu , u|_{\Omega_2} -u_2\right\rangle=0 \;\forall \mu \in \Lambda,\\
\sup_{(v,v_2)\in V \times V_2}\dfrac{\left\langle\lambda , v|_{\Omega_2} -v_2\right\rangle}{ \left\lbrace \norm{v}_{V}^2+ \norm{v_2}_{V_2}^2\right\rbrace^{\frac{1}{2}}}&\geq C_2 \norm{\lambda }_{\Lambda}&&\forall \lambda \in \Lambda.
\end{align*}
Equivalently to the two previous conditions, we write the following inf-sup condition, see \cite{xu2003some}.

There exists a constant $C >0$ such that:
\begin{align}
\label{eq:weak_continuous_L_infsup}
 \sup_{(v,v_2,\mu)\in \mathcal{H}} \dfrac{(\beta \nabla u,\nabla u)_{\Omega} +\left\langle \lambda,  v|_{\Omega_2} \right\rangle+((\beta_2-\beta) \nabla u_2,\nabla u_{2})_{\Omega_2} -\left\langle \lambda, v_2 \right\rangle +\left\langle \mu , u|_{\Omega_2} -u_2 \right\rangle }
 {\left\lbrace  \abs{u}_{1,\Omega}^2+\norm{u_2}_{1,\Omega_2}^2+\norm{\lambda}_{[H^{1}(\Omega_2)]^*}^2\right\rbrace^{\frac{1}{2}} \left\lbrace  \abs{v}_{1,\Omega}^2+\norm{v_2}_{1,\Omega_2}^2+\norm{\mu}_{[H^{1}(\Omega_2)]^*}^2 \right\rbrace^{\frac{1}{2}}  } \geq C >0&& \forall (u,u_2,\lambda)\in \mathcal{H}.
\end{align}
\subsection{Discrete FD-DLM problem} \label{se:disc_pbm}
Let $\mathcal{T}_1$ and $\mathcal{T}_2$ be two partitions of $\Omega$ and $\Omega_2$ consisting of quadrilaterals in 2D and hexahedra in 3D.
Denote by $\mathcal{E}_1, \mathcal{E}_2 $ the collections of edges in both meshes.
Moreover, let $h_{K_1},h_{K_2}, h_{E_1}, h_{E_2}$ be the elements ($K_1 \in \mathcal{T}_1$ and $K_2 \in \mathcal{T}_2$~) and edges ($E_1 \in \mathcal{E}_1$ and $E_2 \in \mathcal{E}_2$~) diameters. 

Assume that our partitions are shape regular; meaning there exist two constants, $C_{\mathcal{T}_1}$ and $C_{\mathcal{T}_2}$, associated with the partitions $\mathcal{T}_1$ and $\mathcal{T}_2$, respectively. These constants ensure that for all $K_1 \in \mathcal{T}_1$ and $K_2 \in \mathcal{T}_2$, the shape parameters bound the ratio between the diameter of each element and the diameter of the largest inscribed ball within that element. This bound remains constant and independent of the element and its size.

Let $\mathcal{H}_h \coloneqq V_h \times V_{2h} \times \Lambda_h  $ be the space of the discrete solution triplet $(u_h,u_{2h},\lambda_h)$. Notice that $\Lambda_h \subset L^2(\Omega_2)$. Therefore, the duality pairing can be evaluated in the discrete level as the $L^2$ scalar product in $\Omega_2$; i.e. $\left\langle \lambda_h, v_{2h} \right\rangle = (\lambda_h, v_{2h})_{\Omega_2}$. The discretization of Problem~\ref{pbm:1} reads 
\begin{problem}
\label{pbm:2}
Given $f \in L^2(\Omega)$ and $f_2 \in L^2(\Omega_2)$, find $(u_h,u_{2h},\lambda_h) \in \mathcal{H}_h$ such that
\begin{align*}
(\beta \nabla u_h,\nabla v_h)_{\Omega} + ( \lambda_h,  v_h|_{\Omega_2} )_{\Omega_2}
&=(f,v_h)_{\Omega} && \forall v_h \in V_h\\
((\beta_2-\beta) \nabla u_{2h},\nabla v_{2h})_{\Omega_2} -( \lambda_h, v_{2h} )_{\Omega_2}
&=(f_2-f,v_{2h})_{\Omega_2} && \forall v_{2h} \in V_{2h}\\
( \mu_h , u_h|_{\Omega_2} -u_{2h} )_{\Omega_2}
&=0 && \forall \mu_h \in \Lambda_h.
\end{align*}
\end{problem}
We consider two choices for the space $\mathcal{H}_h$, namely $Q_1-(Q_1+B)-P_0$ and $Q_2-Q_2-P_0$, that are defined as follows:
\begin{itemize}
\item  \textbf{Element 1 [$Q_1-(Q_1+B)-P_0$]:} In this selection of finite element spaces, we approximate the solutions in $V_h$ using continuous piecewise bilinear polynomials, while in $V_{2h}$, we utilize continuous piecewise bilinear polynomials augmented by a bubble  $B$, that is a biquadratic bubble in each element in 2D and triquadratic bubble in each element in 3D. Finally, we approximate the Lagrange multiplier using discontinuous piecewise constant polynomials. A depiction of the element and its associated degrees of freedom in two dimensions can be observed in Figure~\ref{fig:Element1}.
\begin{figure}[htp]
\centering
\includegraphics[scale=0.6]{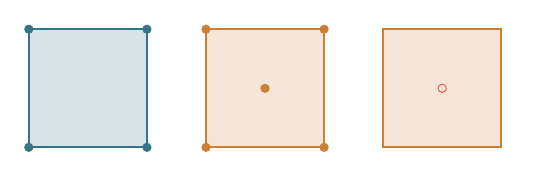} 
\caption{Element 1 [$Q_1-(Q_1+B)-P_0$] in 2D.}\label{fig:Element1}
\end{figure}
\item \textbf{Element 2 [$Q_2-Q_2-P_0$]:} In this approach, we approximate the solutions in $V_h$ and $V_{2h}$ using continuous piecewise biquadratic polynomials. Additionally, akin to the previous element, we approximate the Lagrange multiplier using discontinuous piecewise constant polynomials. Figure~\ref{fig:Element2} provides a two-dimensional representation of the element and its corresponding degrees of freedom. Notably, this element incorporates a higher number of degrees of freedom compared to the previous one, rendering it computationally more demanding.
\begin{figure}[htp]
\centering
\includegraphics[scale=0.6]{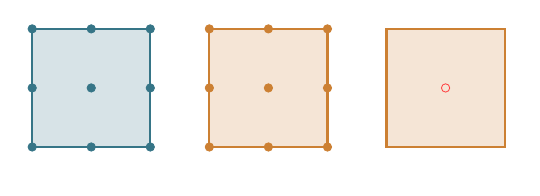} 
\caption{Element 2 [$Q_2-Q_2-P_0$ ] in 2D.}\label{fig:Element2}
\end{figure}
\end{itemize}
In \cite{Najwa2022elliptic}, it was demonstrated that the inclusion of element bubbles in the space $V_{2h}$ is essential for ensuring the stability of the finite element spaces when dealing with a discontinuous Lagrange multiplier space. Notably, the $Q_2$ finite element inherently incorporates a bubble, obviating the need for an additional one. Obviously, $V_h \subset V$, $V_{2h} \subset V_2$ and $\Lambda_h \subset \Lambda$. 

Furthermore, it was established in \cite{Najwa2022elliptic} that by employing these finite element spaces, Problem~\ref{pbm:2} attains stability, and a unique solution exists for the system, under the constraint that $\beta_2 > \beta > 0$ in the domain $\Omega_2$.
The same result was obtained in~\cite{auricchio2015fictitious} in the case of
continuous piecewise linear Lagrange multiplier. In~\cite{Najwa2022elliptic}
and~\cite{auricchio2015fictitious}, it was also observed numerically that the
constraint on $\beta_2$ and $\beta$ could be relaxed. This fact was proved
rigorously in~\cite{boffi2014mixed} in the case of continuous piecewise linear
Lagrange multiplier.

The stability and unique solvability of the Problem~\ref{pbm:2} were established by demonstrating the following two conditions.

There exist constants $C_3, C_4 >0$ such that:
\begin{align*}
(\beta \nabla u_h,\nabla u_h)_{\Omega} + ((\beta_2-\beta) \nabla u_{2h},\nabla u_{2h})_{\Omega_2} &\geq C_3 \left\lbrace \norm{u_h}_{V}^2+ \norm{u_{2h}}_{V_2}^2\right\rbrace &&\forall (u_h,u_{2h}) \; \textit {such that } (\mu_h,u_h |_{\Omega_2}-u_{2h})_{\Omega_2}=0 \; \forall \mu_h \in \Lambda_h, \\
 \sup_{(v_h,v_{2h})\in V \times V_2}\dfrac{( \lambda_h,v_h |_{\Omega_2} -v_{2h})_{\Omega_2}}
 { \left\lbrace \norm{v_h}_{V}^2+ \norm{v_{2h}}_{V_2}^2\right\rbrace^{\frac{1}{2}}}&\geq C_4 \norm{\lambda_h }_{\Lambda}&& \forall \lambda_h \in \Lambda_h, 
\end{align*}
where $C_3$ and $C_4$ are independent of the mesh sizes. 
%
\section{The a posteriori error estimator}\label{se:aposteriori}
Without loss of generality, we focus in this work on a two-dimensional scenario, with the extension to three dimensions being a straightforward continuation of this study. The objective of this paper is to propose error estimators that are proportional to the norms of the errors. In other words, we aim to introduce error estimators that both overestimate and underestimate the corresponding error norms, along with corrective terms tied to data perturbation. These estimators should not rely on knowledge of the exact solutions. This characteristic makes them particularly advantageous in real-world applications where the exact solutions are often unknown in practical applications of most partial differential equations. Furthermore, the a posteriori error estimators should approximate, within a constant factor, the behavior of the exact errors. Specifically, these estimators should exhibit the same convergence rate and be capable of identifying local regions where larger errors are expected. The latter feature is essential for leveraging the estimators in adaptive refinement strategies. The construction of such estimators typically unfolds in two stages.
\begin{enumerate}
\item The global error norms are constrained by both the estimators and the fluctuations in the data, a concept known by the Global Upper Bound (GUB) or reliability.
\item The local error norms are bounded from below by the local indicators and the local fluctuations in the data, known by the Local Lower Bound (LLB) or efficiency.
\end{enumerate}

\begin{figure}[htp]
\begin{center}
	\includegraphics[scale=0.3]{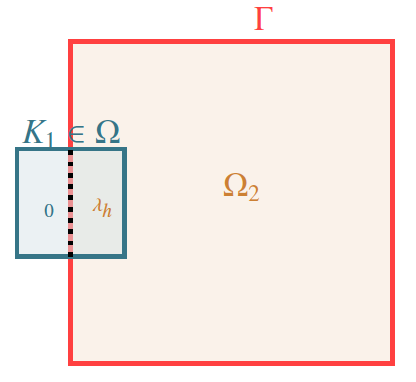}
	\caption{Extension of $\lambda_h$ out of $\Omega_2$.}
	\label{fig:extended_lambda}
\end{center}
\end{figure}

There exist various possibilities for a posteriori error estimators. We refer the reader to the book \cite{verfurth2013posteriori} by Verf{\"u}rth and R{\"u}diger for a review on a posteriori error techniques. In this study, we focus on residual-based estimators. For the design of the residual-based local error indicators, a crucial role is played by the evaluation of the residual element by element. In situations like the one depicted in Figure~\ref{fig:extended_lambda}, where the cell $K_1$ of the mesh $\mathcal{T}_1$ partially intersect the domain $\Omega_2$, we have to give a precise meaning to the quantity $\lambda-\lambda_h$ which is defined only in $\Omega_2$. This underscores the need for extensions of the Lagrange multipliers, $\lambda$ in the continuous problem and $\lambda_h$ in the discrete problem.

It is evident that $\lambda_h$ is a function. Therefore, its extension is straightforwardly defined in the following definition.\\

%
\begin{definition}
\label{def:extended_lambdah}
Recall that $\lambda_h \in L^2(\Omega_2)$ is a piecewise constant polynomial. Therefore, $\lambda_h\in L^2(\Omega)$ can be extended by zero in $\Omega$. We denote this extension by $\tilde{ \lambda}_h$ that is defined as follows
\begin{align*}
\tilde{ \lambda}_h:\begin{cases}
 0 &\Omega\setminus\Omega_2\\
\lambda_h  & \Omega_2,
\end{cases}
\end{align*}
such that $( \tilde{ \lambda}_h,v_h)_{\Omega}\coloneqq (\lambda_h,v_h|_{\Omega_2} )_{\Omega_2}$, as presented in Figure~\ref{fig:extended_lambda}. 
\end{definition}
Using the definition of $\tilde{\lambda}_h$, the first equation in Problem~\ref{pbm:2} can be written as follows
 \begin{align}
 \label{eq:extended_pbm2}
(\beta \nabla u_h,\nabla v_h)_{\Omega} + (\tilde{\lambda}_h, v_h)_{\Omega} &=(f,v_h)_{\Omega} && \forall v_h \in V_h.
 \end{align}
On the other hand, $\lambda$ is defined on $[H^1(\Omega_2)]^*$. Hence, the extension is not as trivial as for $\lambda_h$. The subsequent definition delineates this extension.\\

\begin{definition}
\label{def:extended_lambda}
Recall that $\lambda \in [H^1(\Omega_2)]^*$. We extend $\lambda$ to $\Omega$  and we denote this extension by $\tilde{ \lambda}$, which is  defined as follows
\[\tilde{ \lambda}:H^1_0(\Omega) \rightarrow \mathbb{R}.\]
The extension of $\lambda$ to $\Omega$, denoted by $\tilde{\lambda}$, is the element of $H^{-1}(\Omega)$ space such that:
\begin{align*}
{}_{H^{-1}(\Omega)}\left\langle \tilde{ \lambda},v\right\rangle{}_{H^1_0(\Omega)}\coloneqq {}_{[H^1(\Omega_2)]^*}\left\langle\lambda,v|_{\Omega_2} \right\rangle_{H^1(\Omega_2)}, && \forall v \in H^1_0(\Omega).
\end{align*}

\end{definition}

Moreover, we have the following bound for its $H^{-1}(\Omega)$- norm:
\begin{align}\label{eq:norm_lambda_extended}
&\norm{\tilde{\lambda}}_{H^{-1}(\Omega)}=\sup_{v \in H^{1}_0(\Omega)} \dfrac{\left\langle \tilde{\lambda},v\right\rangle }{\norm{v}_{1,\Omega}}=\sup_{v \in H^{1}_0(\Omega)} \dfrac{\left\langle \lambda,v|_{\Omega_2}\right\rangle }{\norm{v}_{1,\Omega}} \leq \sup_{v \in H^{1}_0(\Omega)} \dfrac{\left\langle \lambda,v|_{\Omega_2}\right\rangle }{\norm{v|_{\Omega_2}}_{1,\Omega_2}} \leq \sup_{w \in H^{1}(\Omega_2)} \dfrac{\left\langle \lambda,w\right\rangle }{\norm{w}_{1,\Omega_2}}= \norm{\lambda}_{[H^{1}(\Omega_2)]^*}.
\end{align}

Using the definition of $\tilde{\lambda}$, the first equation in Problem~\ref{pbm:1} can be written as follows
 \begin{align}
 \label{eq:extended_pbm}
(\beta \nabla u,\nabla v)_{\Omega} + \left\langle\tilde{\lambda}, v\right\rangle &=(f,v)_{\Omega} && \forall v \in V.
 \end{align}
%

On what follows, we denote by $C$ a general positive constant that depends solely on the shape parameters $C_{\mathcal{T}_1}$ and $C_{\mathcal{T}_2}$.
Below, we list different possibilities of collections of elements that are needed in this work.
\begin{figure}[htp]\centering
\begin{minipage}[c]{0.2225\linewidth}
\centering
		\includegraphics[width=0.8\linewidth]{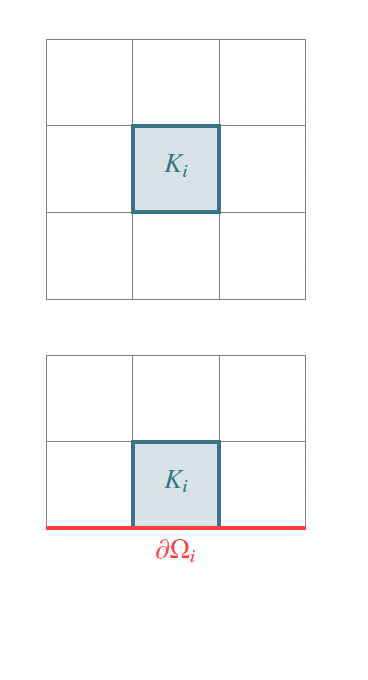}
		\caption{$\overline{\omega}_{K_i}$.}
		\label{fig:omegat}
\end{minipage}
\begin{minipage}[c]{.2225\linewidth}
\centering
		\includegraphics[width=0.8\linewidth]{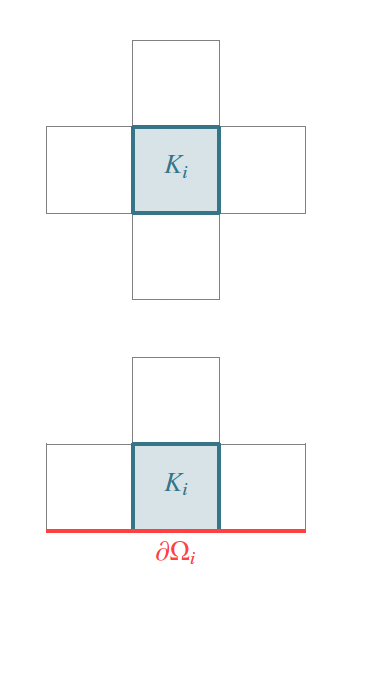} 
		\caption{$\omega_{K_i}$.}
		\label{fig:omegaT1}
\end{minipage}
\begin{minipage}[c]{.2225\linewidth}
\centering
		\includegraphics[width=0.8\linewidth]{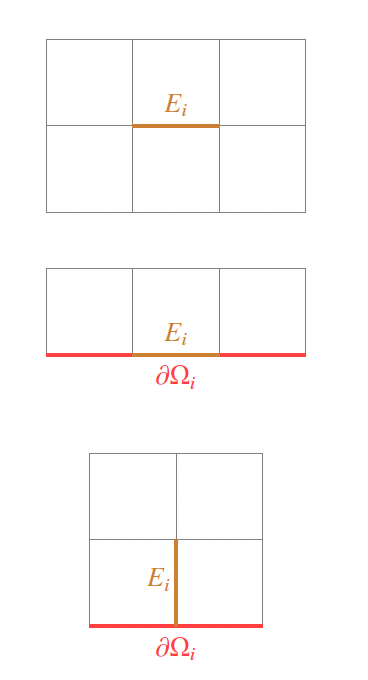} 
		\caption{$\overline{\omega}_{E_i}$.}
		\label{fig:omegae}
\end{minipage}
\begin{minipage}[c]{.2225\linewidth}
\centering
		\includegraphics[width=0.8\linewidth]{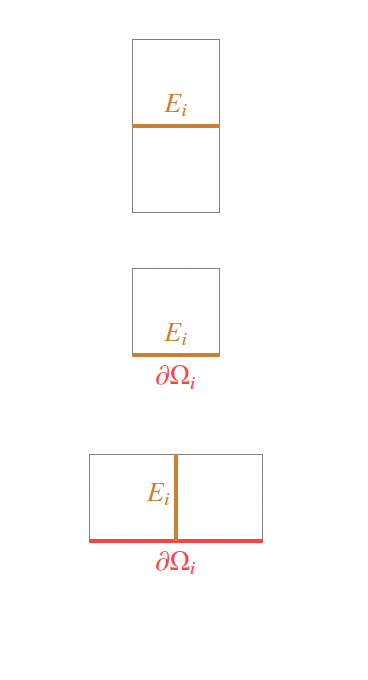} 
		\caption{${\omega_{E_i}}$.}
		\label{fig:omegaE1}
\end{minipage}
\end{figure}
\begin{itemize}
\item $\overline{\omega}_{K_1},\overline{\omega}_{K_2}$: The collection of all elements that share at least one vertex with the elements $K_1, K_2$ respectively,  as seen in Figure~\ref{fig:omegat}.
\item $\omega_{K_1},\omega_{K_2}$: The collection of all elements that share an edge with the elements $K_1$ or $K_2$ respectively,  as seen in Figure~\ref{fig:omegaT1}.
\item$\overline{\omega}_{E_1},\overline{\omega}_{E_2}$: The collection of all elements that share at least a vertex with the edges $E_1$ or $E_2$, respectively, as Figure~\ref{fig:omegae} shows.
\item$\omega_{E_1},\omega_{E_2}$: The collection of all elements that share the edges $E_1$ or $E_2$, respectively, as presented in Figure~\ref{fig:omegaE1}.
\end{itemize}

From now on, we write for short $f_3=f_2-f$ and $\beta_3 =\beta_2 -\beta$. Section~\ref{se:constant_beta} initiates our exploration by introducing estimators tailored for the scenario where the coefficients $\beta$ and $\beta_2$ are constant. We demonstrate the reliability and efficiency of these estimators, deliberately focusing on this case to streamline the core aspects of the proof while sidestepping the additional complexities stemming from approximating non-constant coefficients.
Subsequently, in Section~\ref{se:continuous_beta}, we extend our approach to encompass the more general case involving general smooth coefficients. While we refrain from presenting comprehensive proof in this section, we delineate the primary differences and elucidate strategies for addressing how to bound the new terms. Through this analysis, we establish that the estimators retain their reliability and efficiency even in the broader context of general smooth coefficients. 
\subsection{Constant diffusion coefficients}
\label{se:constant_beta} 
In this section, we assume that the diffusion coefficient $\beta_1$ and $\beta_2$ are constants in the subdomains $\Omega_1$ and $\Omega_2$, respectively. Hence, the extended coefficient $\beta$ is equal to $\beta_1$ over all $\Omega$.

Furthermore, we assume that $f\in L^2(\Omega)$ and $f_2\in L^2(\Omega_2)$ . Consequently, we must account for errors arising from data oscillation in our estimators. Let $\Pi_0$ be the $ L^2-$projection onto the space of piecewise constant polynomials in the cell $ K_i\in\mathcal{T}_i,  (i=1,2)$.
Subsequently, in Definition~\ref{def:residuals}, we define the residuals of elements and sides in $\Omega$, taking into account Equation~\eqref{eq:extended_pbm2}, as well as the residuals of elements and sides in $\Omega_2$, based on the second equation in Problem~\ref{pbm:2}.\\

\begin{definition}\label{def:residuals}
The residual of the element $K_1 \in \mathcal{T}_1$ is
\[
R_{K_1}(u_h)=\beta\; \Delta u_h - \tilde{ \lambda}_h+\Pi_0 f,
\]
where $\tilde{ \lambda}_h $   is the extension of $\lambda_h$  defined in Definition~\ref{def:extended_lambdah}. The residual of the edge $E_1 \in \mathcal{E}_1$ is
\[R_{E_1}(u_h)=\begin{cases} 
-\beta \;  \jump{\dfrac{\partial u_h}{\partial \textbf{n}_1}}_{E_1} &  E_1 \in \partial K_1 \setminus \partial \Omega\\
0 &  E_1 \in \partial \Omega,
\end{cases}
\]
where $\textbf{n}_1$ is the outward normal to $K_1$ and  $\jump{ \cdot } $ is the jump of the normal derivative  $\dfrac{\partial u_h}{\partial \textbf{n}_1} $  across the edge $E_1$. Similarly, the residual of the element $K_2 \in \mathcal{T}_2$ and the residual of the edge $E_2 \in \mathcal{E}_2$ are
\begin{align*}
	R_{K_2}(u_{2h})&=\beta_3 \; \Delta u_{2h} + \lambda_h+\Pi_0 f_3,\\
	R_{E_2}(u_{2h})&=\begin{cases} 
	- \beta_3\;\jump{ \dfrac{\partial u_{2h}}{\partial \textbf{n}_2}}_{E_2} & \forall E_2 \in \partial K_2 \setminus \Gamma \\
	- \beta_3\; \dfrac{\partial u_{2h}}{\partial \textbf{n}_2} &  \forall E_2 \in \partial K_2\cap \Gamma,
	\end{cases}
\end{align*}
where $\textbf{n}_2$ is the outward normal to $K_2$ and  $\jump{ \cdot } $ is the jump of the normal derivative $\dfrac{\partial u_{2h}}{\partial \textbf{n}_2} $  across the edge $E_2$. 
\end{definition}
We propose the following two a posteriori error estimators that are weighted combinations of the elements and edges residuals.
\begin{definition}\label{def:eta}
For any $K_1 \in \mathcal{T}_1$ and  $K_2 \in \mathcal{T}_2$, define the following local indicators
\begin{align*}
		\eta_{K_1}^2&= \left\lbrace h_{K_1}^2 \norm{R_{K_1}(u_h)}_{0,K_1}^2 
			+ \frac{1}{2}\sum_{E_1 \in\partial K_1 \setminus\partial \Omega} 
			h_{E_1} \norm{R_{E_1}(u_h)}_{0,E_1}^2 \right\rbrace, \\
		\eta_{K_2}^2&= \left\lbrace h_{K_2}^2 \norm{R_{K_2}(u_h)}_{0,K_2}^2
			+ \norm{u_h|_{\Omega_2} -u_{2h}}^2_{1,K_2} 
			+ \frac{1}{2} \sum_{E_2 \in\partial K_2\setminus \Gamma}h_{E_2} 
			\norm{R_{E_2}(u_{2h})}_{0,E_2}^2 
			+  \sum_{E_2 \in\partial K_2\cap\Gamma}h_{E_2} 
			\norm{R_{E_2}(u_{2h})}_{0,E_2}^2 \right\rbrace .
\end{align*}
To prevent redundant calculations of jumps from interior edges shared by neighboring elements, the terms associated with interior edges are scaled by $1/2$. This approach ensures efficiency by accounting for each jump only once from adjacent elements sharing the same edge. The aggregation of these local indicators across $\mathcal{T}_1$ and $\mathcal{T}_2$ yields the following global estimators
\begin{align*}
		\eta^2_1&=\sum_{K_1 \in \mathcal{T}_1 } \eta_{K_1}^2,
		&&\eta^2_2=\sum_{ K_2 \in \mathcal{T}_2 } \eta_{K_2}^2.
\end{align*}
\end{definition}
The estimators defined in Definition~\ref{def:eta} admit a global upper bound (reliability) and a local lower bound (efficiency) proven in the next section. In what follows, we write to simplify the notation
\begin{align*} osc_{K_1}= h_{K_1} \norm{f-\Pi_0 f}_{0,K_1}, && osc_{K_2}= h_{K_2} \norm{f_3-\Pi_0 f_3)}_{0,K_2},\end{align*}

where those terms are the correction terms that are higher order perturbations of the force. Also referred to as the oscillation terms. 
\subsubsection{Reliability and efficiency of the error estimator}
Our error estimation approach relies on a careful assessment of the residuals associated with the chosen stable elements, as detailed in Section~\ref{se:disc_pbm}. This section begins with the establishment of the Global Upper Bound (GUB) in Proposition~\ref{prop:GUB}, which provides a reliability bound. This bound ensures that the accuracy of the numerical solution remains within the desired tolerance. Moreover, in Proposition~\ref{prop:LLB}, we establish the Local Lower Bound (LLB), which offers an efficiency bound essential for managing the number of degrees of freedom utilized. Attaining the numerical solution with the desired accuracy while minimizing computational costs depends on both the reliability and efficiency bounds outlined in this section.

To show those bounds, we need some tools. The following lemma describes the Cl\' ement interpolation operator, in a general format, which is a classical mathematical tool usually used in the GUB construction along with the associated inequalities, (see~\cite{claement1975approximation}).

 Let $\overline{\omega}_K$ and  $\overline{\omega}_E$ be sets of elements as in Figures~\ref{fig:omegat} and \ref{fig:omegae} respectively, and let $Q_1(\overline{\omega}_K)$ be the finite element space of functions that, when restricted to each element, is a polynomial of order one in each variable. Moreover, $\mathcal{T}$ and $\mathcal{E}$ denote sets of elements and edges of the partition of a whole domain $\tilde{\Omega}$ and  $X_h =\{ w \in H^1(\Omega):w|_K \in Q_1(K)\}$.
\begin{lemma}
\label{lma:I_h_T1}
Let $\mathit{I}: H^1(\overline{\omega}_K) \rightarrow X_h$ be a Cl\' ement interpolation operator. Then, for any element $K$ in the partition $\mathcal{T}$, and any edge $E$ in the set of edges $\mathcal{E}$, there exists a constant $C>0$ that depends only the shape parameter of the partition, such that the following local error estimates are satisfied
\[\norm{v-\mathit{I} v}_{0,K} \leq C h_{K} \norm{\nabla v}_{0,\overline{\omega}_K} \hspace{1in} \forall v \in H^1_0(\overline{\omega}_K), \]
\[\norm{v-\mathit{I} v}_{0,E} \leq C h_{E}^{\frac{1}{2}} \norm{\nabla v}_{0,\overline{\omega}_E}\hspace{1in} \forall v \in H^1_0(\overline{\omega}_E). \]
The collection of these norms over all elements and edges leads to the following bounds
\begin{align*}
\Bigg[\sum_{K \in \mathcal{T}} \norm{v}_{r,\overline{\omega}_{K}}^2\Bigg]^{\frac{1}{2}} &\leq C \norm{v}_{r,\tilde{\Omega}}&& \forall v \in H^1_0(\tilde{\Omega}),  r=0,1,\\
\Bigg[\sum_{E \in \mathcal{E}} \norm{v}_{r,\overline{\omega}_{E}}^2\Bigg]^{\frac{1}{2}} &\leq C \norm{v}_{r,\tilde{\Omega}}&& \forall v \in H^1_0(\tilde{\Omega}), r=0,1.
\end{align*}
\end{lemma}

Let us denote the error of $u,u_2, \lambda,$ and $\tilde{\lambda}$ by $e,e_2,\varepsilon,$ and $\tilde{\varepsilon}$; that is $e \coloneqq u-u_h$, $e_2 \coloneqq u_2-u_{2h}$, $\varepsilon \coloneqq \lambda-\lambda_h$, $\tilde{\varepsilon} \coloneqq \tilde{\lambda}-\tilde{\lambda_h}$. In the following proposition, we establish the global upper bound of the proposed estimators.\\
\begin{proposition}{\textbf{Reliability:}\\}\label{prop:GUB} 
Let $(u,u_2,\lambda) \in \mathcal{H}$ and $(u_h,u_{2h},\lambda_h) \in \mathcal{H}_h$. Then, there exists a constant $\overline{C}>0$, which only depends on $\Omega,~\Omega_2$ and on the shape parameters $C_{\mathcal{T}_1}$ and $C_{\mathcal{T}_2}$ of the families $\mathcal{T}_1,~ \mathcal{T}_2$, such that
\begin{equation}
\label{eq:upper_a_posterori}
\abs{e}_{1,\Omega}+\norm{e_2}_{1,\Omega_2}+\norm{\varepsilon}_{\Lambda} \leq \overline{C} \left\lbrace 
 \sum_{K_1 \in\mathcal{T}_1} \left( \eta_{K_1}^2 + osc_{K_1}^2 \right) +
\sum_{K_2 \in \mathcal{T}_2} \left( \eta_{K_2}^2 + osc_{K_2}^2 \right) 
\right\rbrace ^{\frac{1}{2}},
\end{equation}
\end{proposition}
\begin{proof}
The objective of this proof is to constrain the left-hand side of inequality~\eqref{eq:upper_a_posterori} using terms containing solely the norms of the residuals and higher-order perturbations of the data. To achieve this, we commence with the Galerkin orthogonality. We consider Problems~\ref{pbm:1} and \ref{pbm:2} and we replace the first equations in both problems by their equivalent equations, \eqref{eq:extended_pbm2} and \eqref{eq:extended_pbm}, respectively, to get the following error equations
\begin{subequations}
\begin{align}
(\beta \nabla e,\nabla v_h)_{\Omega} +\left\langle \tilde{\varepsilon} , v_h \right\rangle &=0 && \forall v_h \in V_h \label{eq:orthognalityOfError1}\\
(\beta_3\nabla e_2,\nabla v_{2h})_{\Omega_2} 
			-\left\langle \varepsilon, v_{2h}\right\rangle & =0  && \forall v_{2h} \in V_{2h}\label{eq:orthognalityOfError2}\\
 \left\langle \mu_h,e|_{\Omega_2}-e_2\right\rangle&=0&&\forall \mu_h \in\Lambda_h.\label{eq:orthognalityOfError3}
\end{align}
\end{subequations}

For any $v \in V$ and $v_2 \in V_2$, let $\mathit{I}_h v, \mathit{I}_{2h} v_2 $ be their Cl\' ement like interpolation, defined in Lemma~\ref{lma:I_h_T1}. Then, using the error equation \eqref{eq:orthognalityOfError1}, the definition of $\tilde{\lambda}_h$ (Definition~\ref{def:extended_lambdah}), integration by parts, and Cauchy–Schwarz inequality,  we get
\begin{align*}
		(\beta \nabla e,\nabla v)_{\Omega}
			+\left\langle \tilde{\varepsilon}, v \right\rangle &=(\beta \nabla e,\nabla (v-\mathit{I}_h v))_{\Omega}
			+\left\langle  \tilde{\varepsilon}, (v-\mathit{I}_h v) \right\rangle \\
		&= (\beta \nabla u,\nabla\,(v-\mathit{I}_h v))_{\Omega}+\left\langle \tilde{\lambda},v-\mathit{I}_h 
			v\right\rangle 
			- (\beta \nabla u_h,\nabla\,(v-\mathit{I}_h v))_{\Omega}
			-( \tilde{\lambda_h},v-\mathit{I}_h 
			v)_{\Omega_2}\\
			&=(f,v-\mathit{I}_h v)_{\Omega}
			- (\beta \nabla u_h,\nabla\,(v-\mathit{I}_h v))_{\Omega}
			-( \tilde{\lambda}_h,v-\mathit{I}_h v)_{\Omega}\\
		&= \sum\limits_{\substack{K_1 \in \mathcal{T}_1}} \left\lbrace (f+\Pi_0 f 
			-\Pi_0 f,v-\mathit{I}_h v)_{K_1} + (\beta \Delta u_h ,v
			-\mathit{I}_h v)_{K_1}-( \tilde{\lambda}_h, v-\mathit{I}_h v)_{K_1}  \right\rbrace\\
			&\hspace{0.5in} -\frac{1}{2}\sum_{K_1 \in \mathcal{T}_1}\sum_{E_1 \in\partial 
			K_1\setminus\partial \Omega} (\beta  \jump{
			\frac{\partial u_h}{\partial \textbf{n}_1 }} _{E_1},v-\mathit{I}_h v)_{E_1}\\
		&=\sum_{K_1 \in \mathcal{T}_1} \left\lbrace ( R_{K_1}(u_h) ,
			v-\mathit{I}_h v)_{K_1}+(f-\Pi_0 f,
			v-\mathit{I}_h v)_{K_1}+\frac{1}{2}\sum_{E_1 \in\partial 
			K_1\setminus\partial \Omega}
			( R_{E_1}(u_h) ,v-\mathit{I}_h v)_{E_1} \right\rbrace.\\
		&\leq C \sum_{K_1 \in  \mathcal{T}_1} \left\lbrace \Bigg[\norm{ R_{K_1}(u_h) }_{0,K_1} +\norm{ f-\Pi_0 f}_{0,K_1} \Bigg]\norm{ v-\mathit{I}_h v}_{0,K_1}  + \frac{1}{2} \sum_{E_1 \in\partial K_1\setminus\partial \Omega} \norm{ R_{E_1}(u_h) }_{0,E_1} \norm{ v-\mathit{I}_h v}_{0,E_1} \right\rbrace \\
		&\leq C \sum_{K_1 \in \mathcal{T}_1} \left\lbrace\Bigg[ h_{K_1} 
			\norm{ R_{K_1}(u_h) }_{0,K_1} 
			+ osc_{K_1} \Bigg] \norm{\nabla v }_{0,\overline{\omega}_{K_1}} 
			+\frac{1}{2} \sum_{E_1 \in\partial K_1\setminus\partial \Omega} 
			h_{E_1} ^\frac{1}{2}\norm{ R_{E_1}(u_h) }_{0,E_1}
			\norm{\nabla v }_{0,\overline{\omega}_{E_1} } \right\rbrace. 
\end{align*}
This implies
\begin{equation}
\label{eq:fisrtupper}
(\beta \nabla e,\nabla v)_{\Omega}+\left\langle \varepsilon, v|_{\Omega_2}\right\rangle \leq C  \abs{v}_{1,\Omega} \left\lbrace \sum_{K_1 \in \mathcal{T}_1 } \left(\eta_{K_1}^2+osc_{K_1}^2 \right) \right\rbrace ^\frac{1}{2}.
\end{equation}
We deal similarly with the second equation in Problems~\ref{pbm:1} and \ref{pbm:2}.
\begin{align*}
		(\beta_3\nabla e_2,\nabla v_{2})_{\Omega_2} 
			-\left\langle \varepsilon, v_{2}\right\rangle&= (\beta_3\nabla e_2,\nabla (v_{2}-\mathit{I}_{2h} v_2)_{\Omega_2} 
			-\left\langle \varepsilon, v_{2}-\mathit{I}_{2h} v_2\right\rangle\\
		&= (\beta_3 \nabla u_2,\nabla (v_{2}-\mathit{I}_{2h} v_2))_{\Omega_2}
			-\left\langle \lambda, v_{2} -\mathit{I}_{2h} v_2\right\rangle 
			- (\beta_3 \nabla u_{2h},\nabla (v_{2}-\mathit{I}_{2h} v_2))_{\Omega_2} 
			+( \lambda_h, v_{2} -\mathit{I}_{2h} v_2)_{\Omega_2} \\
		&=(f_3,v_{2}-\mathit{I}_{2h} v_2)_{\Omega_2}
			- (\beta_3 \nabla u_{2h},\nabla (v_{2}-\mathit{I}_{2h} v_2))_{\Omega_2} 
			+( \lambda_h, v_{2} -\mathit{I}_{2h} v_2)_{\Omega_2} \\
		&=\sum_{K_2 \in \mathcal{T}_2} \left\lbrace (f_3+\Pi_0 f_3 
			-\Pi_0 f_3,v_2-\mathit{I}_{2h} v_2)_{K_2} + (\beta_3 \Delta u_{2h} ,v_2
			-\mathit{I}_{2h} v_2)_{K_2} +( \lambda_h, v_2-\mathit{I}_{2h} v_2)_{K_2}\right\rbrace\\
			&\hspace{0.5in}- \sum_{K_2 \in \mathcal{T}_2} \left\lbrace\frac{1}{2}
			\sum_{E_2 \in\partial K_2\setminus \Gamma}(\beta_3\; \jump{
			\frac{\partial u_{2h}}{\partial \textbf{n}_2 }} _{E_2},v_2
			-\mathit{I}_{2h} v_2)_{E_2}+\sum_{E_2 \in\partial K_2\cap\Gamma}( \beta_3
			\frac{\partial u_{2h}}{\partial \textbf{n}_2 },v_2
			-\mathit{I}_{2h} v_2)_{E_2}\right\rbrace \\
		&=\sum_{K_2 \in \mathcal{T}_2} \left\lbrace ( R_{K_2}(u_{2h}) 
			,v_2-\mathit{I}_{2h} v_2)_{K_2}+(f_3-\Pi_0 f_3,v_2-\mathit{I}_{2h} v_2)_{K_2} \right\rbrace\\  
			&\hspace{0.5in}+ \sum_{K_2 \in \mathcal{T}_2} \left\lbrace\frac{1}{2}
			\sum_{E_2 \in\partial K_2\setminus\Gamma}
			( R_{E_2}(u_{2h}),v_2-\mathit{I}_{2h} v_2)_{E_2}
			+\sum_{E_2 \in\partial K_2\cap\Gamma}
			( R_{E_2}(u_{2h}),v_2-\mathit{I}_{2h} v_2)_{E_2}
			\right\rbrace\\
		&\leq C \sum_{K_2 \in \mathcal{T}_2} 
				\Bigg[\norm{ R_{K_2}(u_{2h}) }_{0,K_2}
				+\norm{ f_3-\Pi_0 f_3}_{0,K_2}\Bigg] 
				\norm{ v_2 -\mathit{I}_{2h} v_2}_{0,K_2}\\
				&\hspace{0.5in}+ \sum_{K_2 \in \mathcal{T}_2}
				\left\lbrace  \frac{1}{2} \sum_{E_2 \in\partial K_2\setminus\Gamma}  
				\norm{ R_{E_2}(u_{2h}) }_{0,E_2} 
				\norm{ v_2-\mathit{I}_{2h} v_2}_{0,E_2}
				+\sum_{E_2 \in\partial K_2\cap\Gamma}  
				\norm{ R_{E_2}(u_{2h}) }_{0,E_2} 
				\norm{ v_2-\mathit{I}_{2h} v_2}_{0,E_2}\right\rbrace  \\
	&\leq C \sum_{K_2 \in \mathcal{T}_2}  \Bigg[ h_{K_2} 
		\norm{ R_{K_2}(u_{2h}) }_{0,K_2} 
		+ osc_{K_2} \Bigg]\norm{\nabla v_2 }_{0,\overline{\omega}_{K_2}} \\
		&\hspace{0.5in}+C \sum_{K_2 \in \mathcal{T}_2}\left\lbrace\frac{1}{2} \sum_{E_2 \in\partial K_2\setminus\Gamma} 
		h_{E_2} ^\frac{1}{2}\norm{ R_{E_2}(u_{2h}) }_{0,E_2}
		\norm{\nabla v_2 }_{0,w_{E_2} }
		+\sum_{E_2 \in\partial K_2\cap\Gamma} 
		h_{E_2} ^\frac{1}{2}\norm{ R_{E_2}(u_{2h}) }_{0,E_2}
		\norm{\nabla v_2 }_{0,w_{E_2} } \right\rbrace.
\end{align*}
Finally, for any $\mu \in \Lambda$, \eqref{eq:orthognalityOfError3} implies
\begin{align*}
		 \left\langle \mu,e|_{\Omega_2}-e_2\right\rangle
		&=-\left\langle \mu,u_h|_{\Omega_2} -u_{2h}\right\rangle  
		\leq\norm{\mu}_{\Lambda} \norm{u_h|_{\Omega_2}-u_{2h}}_{1,\Omega_2}.
\end{align*}
Hence, from those two previous results, we deduce that
\begin{equation}
\label{eq:secondupper}
 (\beta_3 \nabla e_2,\nabla v_2)_{\Omega_2} -\left\langle \varepsilon, v_2 \right\rangle +\left\langle \mu ,e|_{\Omega_2}-e_2\right\rangle 
\leq C \left\lbrace \abs{v_2}_{1,\Omega_2}^2
		+\norm{\mu}_{\Lambda}^2\right\rbrace^{\frac{1}{2}} \left\lbrace \sum_{K_2 \in \mathcal{T}_2 } 
		\eta_{K_2}^2+osc_{K_2}^2 
		\right\rbrace ^\frac{1}{2}.
\end{equation}
By adding the results in \eqref{eq:fisrtupper} and \eqref{eq:secondupper}, we are led to
$$
	\begin{aligned}
	(\beta \nabla e,\nabla v)_{\Omega}
			+\left\langle \varepsilon, v|_{\Omega_2}\right\rangle +(\beta_3\nabla e_2,\nabla v_{2})_{\Omega_2} 
			&-\left\langle \varepsilon, v_{2}\right\rangle + \left\langle \mu,e|_{\Omega_2}-e_2\right\rangle\\
		\leq C &\left\lbrace 
		\sum_{K_1 \in\mathcal{T}_1} \left( \eta_{K_1}^2 
		+ osc_{K_1}^2 \right) 
		+\sum_{K_2 \in \mathcal{T}_2} \left( \eta_{K_2}^2 
		+ osc_{K_2}^2 \right)
		 \right\rbrace ^{\frac{1}{2}}
		 \cdot \left\lbrace \abs{v}_{1,\Omega}^2+\norm{\mu}_{\Lambda}^2 + \abs{v_2}_{1,\Omega_2}^2 \right\rbrace^{\frac{1}{2}}.
	\end{aligned}
$$
Leveraging the fulfilled inf-sup condition in \eqref{eq:weak_continuous_L_infsup}, we derive the upper bound for the error in \eqref{eq:upper_a_posterori}.
\end{proof}
That conclude the proof for the reliability of the proposed estimators. It remains to show the efficiency of the estimators. The conventional mathematical tools employed in constructing the local lower bound typically include edges and elements bubble functions. Lemmas~\ref{lma:psiT} and \ref{lma:psiE} detail these functions and their associated inequalities, see \cite{verfurth2013posteriori} Section 1.3.4 and Section 3.6 for the related proofs. Additionally, extension operators are usually used with the edge bubble function which is defined in Definition~\ref{def:lifting_operator}.

\begin{figure}[htp]
\begin{minipage}[c]{.45\linewidth}
\centering
\includegraphics[scale=0.5]{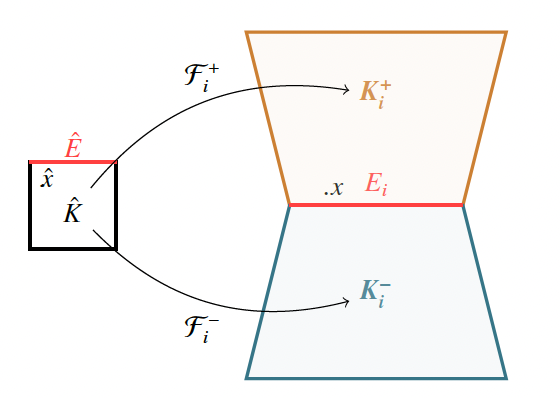} 
\caption{Mapping from the reference element $\hat{K}$.}\label{fig:lifting_oprator}
\end{minipage}
\begin{minipage}[c]{.45\linewidth}
\centering
		\includegraphics[width=0.68\linewidth]{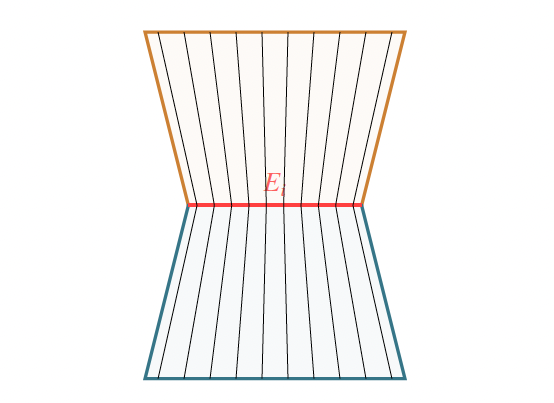} 
		\caption{Level lines.}
		\label{fig:extension_level_lines}
\end{minipage}
\end{figure}

Consider $\mathcal{F}^{-1}$ as a bilinear backward transformation that maps objects from the current element, $K_1 \in \mathcal{T}_1$ or $K_2 \in \mathcal{T}_2$, back to the reference element denoted by $\hat{K}$, a unit square in 2D or a unite cube in 3D. Consequently, $\hat{\varphi} = \varphi \circ \mathcal{F}$, where symbols with hats denote quantities evaluated in the reference domain $\hat{K}$.
\begin{definition}
\label{def:lifting_operator}
Let $\mathbf{P}$ be an extension operator that takes continuous function $\varphi$ defined on a 1D manifold and extends it to a 2D manifold. Consider an edge $ E_i \in \mathcal{E}_i\;(i=1,2)$ and two elements $K_i^+,K_i^- \in \mathcal{T}_i$ that share that edge. Let $\hat{K}$ be a reference unit square in $2D$ and let $\hat{E}$ be one of the edges in the reference element.

Let $\mathcal{F}_i^+$ be the map that transform $\hat{K}$ to $K_i^+$ and $\hat{E}$ to $E_i$, as in Figure~\ref{fig:lifting_oprator}. We define the extension operator of a function $\varphi$ at a point $x \in K_i^+$ as follows:
\[
\mathbf{P}(\varphi)(x)|_{K_i^+} \coloneqq \varphi \left( \mathcal{F}_i^+ \circ \Pi_{\hat{E}} \circ \left(\mathcal{F}_i^+\right)^{-1} (x)\right),
\]
where $\left(\mathcal{F}_i^+\right)^{-1} (x)=\hat{x}$ and  $\Pi_{\hat{E}}(\hat{x})$ is the $L^2-$projection of $\hat{x}$ to the edge $\hat{E}$.  
$\mathbf{P}(\varphi)$ is constant along the level lines presented in Figure~\ref{fig:extension_level_lines}.
\end{definition}

\begin{lemma}
\label{lma:psiT}
Let $\Psi_{K_i}\; (i=1,2)$ be the bubble function defined on elements $ K_1 \in \mathcal{T}_1$ and $ K_2 \in \mathcal{T}_2$ and constructed using their barycenter coordinates, such that
\begin{align*}
0 &\leq \Psi_{K_i }\leq 1, && \supp \Psi_{K_i} \subset K_i.
\end{align*}
Then, for all polynomials $\varphi$ of degree $k$ in each variable, there exists a constant $C$ that depends only on the degree of $\varphi$ and the shape parameter of the mesh such that this bubble function satisfies the following inequalities
\begin{align}
C \norm{\varphi}_{0,K_i } &\leq \norm{\Psi_{K_i}^{\frac{1}{2}} \varphi}_{0,K_i},\label{eq:psiT1}\\
\norm{\Psi_{K_i} \varphi}_{0,K_i} &\leq  C \norm{\varphi}_{0,K_i},\label{eq:psiT2} \\
\abs{\Psi_{K_i} \varphi}_{1,K_i}& \leq C h_{K_i}^{-1} \norm{\varphi}_{0,K_i}.\label{eq:psiT3}
\end{align}
\end{lemma}
\begin{lemma}
\label{lma:psiE}
Let $\Psi_{E_i}\; (i=1,2)$ be the bubble function defined on edges $E_i \in \mathcal{E}_i$ $(i=1,2)$ and constructed using their barycenter coordinates, such that
\begin{align*}
0 &\leq \Psi_{E_i} \leq 1, && \supp \Psi_{E_i} \in \omega_{E_i}.
\end{align*}

For all polynomials $\varphi$  of degree $k$ in each variable, let $\mathbf{P}\varphi$ be the extension of $ \varphi$, defined in Definition~\ref{def:lifting_operator}, such that for every point $x_i \in E_i$ the following holds
$$\mathbf{P}\varphi \coloneqq   \varphi(x_i).$$
Then, there exists a constant $C$ that depends only on the degree of $\varphi$ and the shape parameter of the mesh such that this bubble function satisfies the following inequalities
\begin{align}
 C \norm{\varphi}_{0,E_i} &\leq \norm{\Psi_{E_i}^{\frac{1}{2}} \varphi}_{0,E_i},\label{eq:psiE1}\\
 \norm{\Psi_{E_i} \mathbf{P}\varphi}_{0,\omega_{E_i}} &\leq C h_{E_i}^{\frac{1}{2}} \norm{\varphi}_{0,E_i},\label{eq:psiE2} \\
\abs{\Psi_{E_i} \mathbf{P}\varphi}_{1,\omega_{E_i}} &\leq C h_{E_i}^{-\frac{1}{2}} \norm{\varphi}_{0,E_i}. \label{eq:psiE3}
\end{align}
\end{lemma}
Recall that $\omega_{E_i}$ is the union of all sets that share the edge $E_i$, which contains two elements if $E_i$ is an interior edge and one element if $E_i$ is on the boundary as in Figure~\ref{fig:omegaE1}. In the forthcoming proposition, we derive the local lower bound of the proposed error estimators.\\
\begin{proposition}{\textbf{Efficiency:}\\}\label{prop:LLB}
There exist constants $\underline{C_1}>0$ and  $\underline{C_2}>0$, which only depend on $\Omega,~\Omega_2$ and the shape parameters $C_{\mathcal{T}_1}$ and $C_{\mathcal{T}_2}$ of the families $\mathcal{T}_1,~\mathcal{T}_2$ such that the local errors bound the local indicators as follows
\begin{equation}
\label{eq:lower_a_posterori_T1T2}
	\begin{aligned}
\eta_{K_1} &\leq \underline{C}_1 \left\lbrace  
						\beta \norm{\nabla e}_{0,\omega_{K_1}}
						+\norm{\tilde{\varepsilon}}_{[H^1(\omega_{K_1})]^*}
						+\sum_{K'_1 \in \omega_{E_1}} osc_{K'_1}
						\right\rbrace,\\
\eta_{K_2} &\leq \underline{C}_2  \left\lbrace
						\overline{\beta}_2 \norm{\nabla e_2}_{0
						,\omega_{K_2}}
						+\norm{\nabla e|_{\Omega_2}}_{0,\omega_{K_2}}  
						+\norm{\varepsilon}_{[H^{1}(\omega_{K_2})]^*}	
						+\sum_{K_2' \in \omega_{E_2}} osc_{K_2'}
						 \right\rbrace,
	\end{aligned}
\end{equation}
where $\overline{\beta}_2=\max{\left\lbrace \beta_3,1\right\rbrace}$.
\end{proposition}
The global efficiency is determined by summing up all terms across $\mathcal{T}_1$ and $\mathcal{T}_2$. It is important to note that the combined sums of the terms $\norm{\tilde{\varepsilon}}_{[H^1(\omega_{K_1})]^*}$ and $ \norm{\varepsilon}_{[H^{1}(\omega_{K_2})]^*}$ are bounded from above by $\norm{\tilde{\varepsilon}}_{[H^1(\Omega)]^*}$ and $ \norm{\varepsilon}_{[H^{1}(\Omega_2)]^*}$, respectively. This assertion extends straightforwardly from the proof of Lemma 1.1 in \cite{schatz1977interior}. Moreover, the global efficiency is evident in the results of the numerical experiments in Section~\ref{se:numerical_results}.
\begin{proof}
The objective of this proof is to locally bound the indicators $\eta_{K_1}$ and $\eta_{K_2}$, for all $K_1 \in \mathcal{T}_1$, and for all $K_2 \in \mathcal{T}_2$ by the norms of the errors and higher order oscillation terms.
We begin by establishing the first inequality in \eqref{eq:lower_a_posterori_T1T2}.

Fix an arbitrary element $K_1 \in \mathcal{T}_1$ and set $ v_{K_1}= R_{K_1}(u_h) \Psi_{K_1}$, where $\Psi_{K_1}$ denotes the bubble function described in Lemma~\ref{lma:psiT}. We obtain
	\begin{align*}
		(R_{K_1}(u_h),v_{K_1})_{K_1}&=(\Pi_0 f
											+\beta \Delta u_h,v_{K_1})_{K_1}
											-(\lambda_h,v_{K_1}|_{\Omega_2})_{K_1}
										=(f+ \beta \Delta u_h
											-(f- \Pi_0 f),v_{K_1})_{K_1}
											-(\lambda_h,v_{K_1}|_{\Omega_2})_{K_1}.
	\end{align*}
Integrating by parts and rearranging terms yields
\begin{equation}\label{eq:RK10}
(R_{K_1}(u_h),v_{K_1})_{K_1}= (\beta \nabla e ,\nabla v_{K_1})_{K_1} 
											+ {}_{[H^{1}(K_1)]^*}\left\langle \varepsilon,v_{K_1}|_{\Omega_2}\right\rangle_{H^1(K_1)} 
											-(f- \Pi_0 f ,v_{K_1})_{K_1}.
\end{equation}

There are three potential scenarios for the placement of the
element $K_1$ within the mesh $\mathcal{T}_1$, as depicted in
Figure~\ref{fig:possibilities_K1}. Nevertheless, all terms in \eqref{eq:RK10}
can be uniformly handled across these scenarios, except for the term involving the Lagrange multiplier.
Given that the Lagrange multiplier is defined on functions within $\Omega_2$,
the positioning of $K_1$ becomes crucial when evaluating the second term to
the right hand side of \eqref{eq:RK10}. In the leftmost scenario of
Figure~\ref{fig:possibilities_K1}, where the element lies entirely within
$\Omega_1=\Omega \setminus \Omega_2$, the Lagrange multiplier is not defined
in $K_1$. Therefore, $\tilde{\lambda}$ can be utilized since $\tilde{\lambda} =0 \in \Omega \setminus \Omega_2$. In the central scenario, where the element is entirely within $\Omega_2$, there are no complications in using $\lambda$. Likewise, $\tilde{\lambda}$ can be used since $\tilde{\lambda}=\lambda \in \Omega_2$. However, extra caution is warranted when dealing with the rightmost scenario in Figure~\ref{fig:possibilities_K1}. The necessity to consider the extended Lagrange multiplier $\tilde{\lambda}$ arises because $K_1$ partially intersects $\Omega_2$. We choose to write that term considering $\tilde{\varepsilon}$ over $\varepsilon$ and we consider the following equivalent Equation~\eqref{eq:RK1} instead of Equation~\eqref{eq:RK10} since the extended Lagrange multiplier is effective for all scenarios. Therefore, we opt to present the proof considering it, given that the other cases can be addressed more straightforwardly.\\
Using \eqref{eq:extended_pbm}, integrating by parts and rearranging terms yields
\begin{equation}\label{eq:RK1}
(R_{K_1}(u_h),v_{K_1})_{K_1}= (\beta \nabla e ,\nabla v_{K_1})_{K_1} 
											+ {}_{H^{-1}(K_1)}\left\langle \tilde{ \varepsilon},v_{K_1}\right\rangle_{H^1_0(K_1)} 
											-(f- \Pi_0 f ,v_{K_1})_{K_1}.
\end{equation}

\begin{figure}[htp]
\centering
 \includegraphics[scale=0.5]{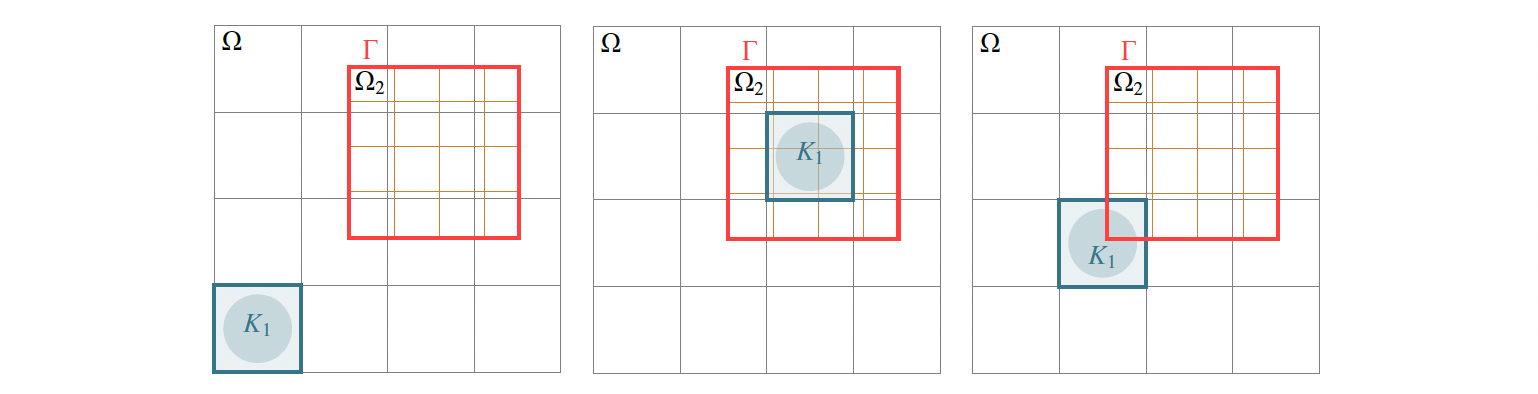} 
\caption{Three different possibilities for the position of the element $K_1$. }\label{fig:possibilities_K1}
\end{figure}

 From \eqref{eq:psiT1}, the left hand side of Equations~\eqref{eq:RK1} can be bounded as follows 
$$(R_{K_1}(u_h),v_{K_1})_{K_1}=(R_{K_1}(u_h)^2,\Psi_{K_1})_{K_1} \geq C \norm{R_{K_1}(u_h)}_{0,K_1}^2.$$
For the bound of the first term on the right hand side of Equations~\eqref{eq:RK1}  we use \eqref{eq:psiT3} to get
\begin{align*}
 (\beta \nabla e ,\nabla v_{K_1})_{K_1} 
			 &\leq \beta  \norm{ \nabla e  }_{0,K_1} 
			 	\norm{\nabla v_{K_1}}_{0,K_1} 
			=\beta  \norm{ \nabla e  }_{0,K_1} \abs{ R_{K_1}(u_h) 
				\Psi_{K_1}}_{1,K_1} 
			\leq C \; \beta \ h_{K_1}^{-1}  \norm{ e  }_{1,K_1}  
				\norm{R_{K_1}(u_h)}_{0,K_1}.
\end{align*}

Since $v_{K_1}$ is supported in $K_1$, which means it lives in $H^1_0(K_1)$, then restricting the error $\tilde{\varepsilon}$ to $K_1$ and using the same arguments in \eqref{eq:norm-1}, \eqref{eq:norm1*}, and \eqref{eq:norm_lambda_extended} we obtain
\begin{align*}
			{}_{H^{-1}(K_1)}\left\langle \tilde{\varepsilon}|_{K_1},v_{K_1} 
				\right\rangle_{H^1_0(K_1)}
			&\leq C \norm{\tilde{\varepsilon}}_{-1,K_1}
				 \norm{\nabla v_{K_1}}_{0,K_1}
			= C \norm{\tilde{\varepsilon}}_{-1,K_1} 
				\abs{ R_{K_1}(u_h) \Psi_{K_1}}_{1,K_1}
			\leq C h_{K_1}^{-1}\norm{\tilde{\varepsilon}}_{[H^1(\omega_{K_1})]^*}
				 \norm{R_{K_1}(u_h)}_{0,K_1}.
\end{align*}
For the last term of Equations~\eqref{eq:RK1}, we have the following bound
\begin{align*}
	(f- \Pi_0 f ,v_{K_1})_{K_1} 
				&\leq \norm{ f- \Pi_0 f}_{0,K_1} \norm{ R_{K_1}(u_h)}_{0,K_1}.
\end{align*}
Adding all results up,  multiplying by $h_{K_1}$, and dividing  by $\norm{R_{K_1}(u_h)}_{K_1}$ yields
\begin{equation}
\label{eq:R_T1LLB}
h_{K_1} \norm{R_{K_1}(u_h)}_{0,K_1} \leq C \left\lbrace \beta   \norm{ e  }_{1,K_1}+\norm{\tilde{\varepsilon}}_{[H^1(\omega_{K_1})]^*} +osc_{K_1}\right\rbrace.
\end{equation}
Now, fix an arbitrary interior edge $E_1 \in \mathcal{E}_1 $ and set $v_{E_1}= \Psi_{E_1} \mathbf{P}R_{E_1}(u_h)$, we have\\
\begin{align}\label{eq:RE10}
(R_{E_1}(u_h),v_{E_1})_{E_1}&= - (\beta \; \jump{ \frac{\partial u_h}{\partial \textbf{n}_1 }
									}_{E_1}, v_{E_1})_{E_1}\nonumber\\
									&= -{}_{H^{-1}(\omega_{E_1})}\left\langle \varepsilon ,
									 v_{E_1}|_{\Omega_2}\right\rangle_{H^1_0(\omega_{E_1})}
									 -\sum_{K_1 \in \omega_{E_1}}  
									\Bigg[ ( R_{K_1} (u_h),v_{E_1})_{K_1}
									-(f- \Pi_0 f,v_{E_1})_{K_1} 
									+\beta (\nabla e,\nabla  v_{E_1} )_{K_1}\Bigg]. 
\end{align}

Once again, similar to the consideration for $K_1$, the term
of utmost importance here is the Lagrange multiplier term which is the first
term to the right hand side of Equation~\eqref{eq:RE10}. The position of the
edge $E_1$ with respect to the immersed domain $\Omega_2$ is sensitive. For
the first term on the right-hand side of Equations~\eqref{eq:RE10}, various
scenarios exist for the position of the edge $E_1$, as depicted in
Figure~\ref{fig:possibilities_E1}. Clearly, the leftmost case represents a
scenario where the Lagrange multiplier term is not defined in $E_1$. Therefore, we can utilize $\tilde{\lambda}$ since it vanishes outside of $\Omega_2$ by definition. In the second leftmost case, where the edge lies within $\Omega_2$, $\lambda$ can be employed, or equivalently $\tilde{\lambda}$ since it coincides with $\lambda$ by definition in $\Omega_2$. The critical cases are the last two, where the support of the test function intersects $\Omega_2$ partially. For these cases, $\tilde{\lambda}$ needs to be utilized. Therefore, we choose to write the term considering $\tilde{\varepsilon}$ instead of $\varepsilon$. We present the proof in a general form, considering $\tilde{\lambda}$ and we write the following equivalent Equations~\eqref{eq:RE1} since the extended case is effective for all scenarios.
\begin{align}\label{eq:RE1}
(R_{E_1}(u_h),v_{E_1})_{E_1}&= - (\beta \; \jump{ \frac{\partial u_h}{\partial \textbf{n}_1 }
									}_{E_1}, v_{E_1})_{E_1}\nonumber\\
									&= -{}_{H^{-1}(\omega_{E_1})}\left\langle \tilde{\varepsilon} ,
									 v_{E_1}\right\rangle_{H^1_0(\omega_{E_1})}
									 -\sum_{K_1 \in \omega_{E_1}}  
									\Bigg[ ( R_{K_1} (u_h),v_{E_1})_{K_1}
									-(f- \Pi_0 f,v_{E_1})_{K_1} 
									+\beta (\nabla e,\nabla  v_{E_1} )_{K_1}\Bigg]. 
\end{align}

\begin{figure}[htp]
\centering
 \includegraphics[scale=0.5]{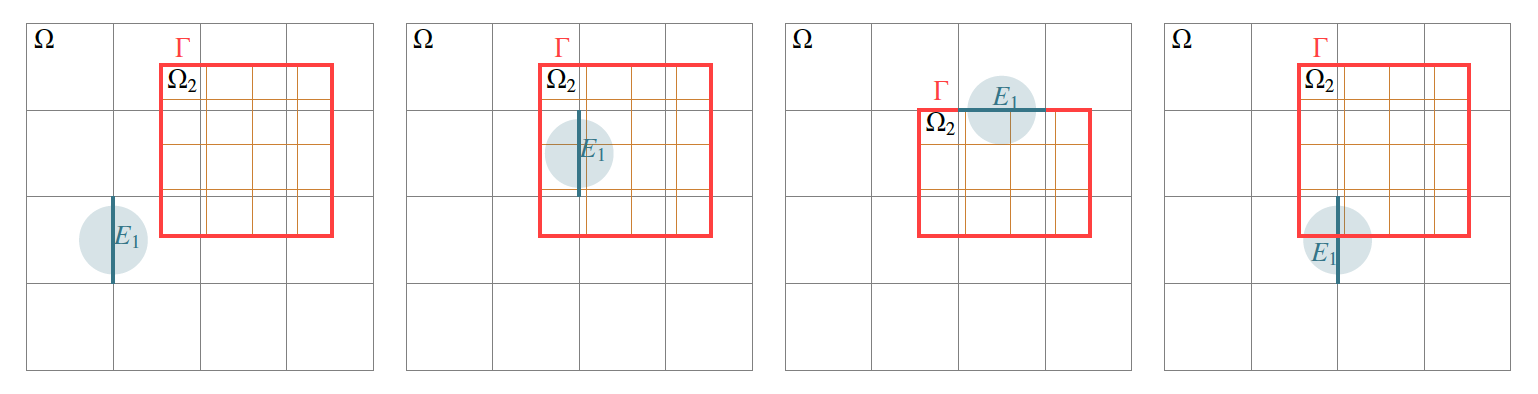} 
\caption{Four different possibilities for the position of the edge $E_1$. }\label{fig:possibilities_E1}
\end{figure}

We bound the left hand side of Equations~\eqref{eq:RE1} as follows
\begin{align*}
 (R_{E_1}(u_h),v_{E_1})_{E_1}&=(R_{E_1}(u_h),  \Psi_{E_1} R_{E_1}(u_h) )_{E_1}
									\geq C \norm{ R_{E_1}(u_h)}_{0,E_1}^2.
\end{align*}

With that in mind, we have the following bound
\begin{align*}									
{}_{H^{-1}(\omega_{E_1})}\left\langle \tilde{\varepsilon}|_{\omega_{E_1}},v_{E_1}\right\rangle_{H^1_0(\omega_{E_1})}
&\leq C \norm{\tilde{\varepsilon}}_{-1,\omega_{E_1}} \norm{\nabla v_{E_1}}_{0,\omega_{E_1}}
\\&= C \norm{\tilde{\varepsilon}}_{-1,\omega_{E_1}} \abs{  \Psi_{E_1}\mathbf{P}R_{E_1}(u_h)}_{1,\omega_{E_1}}
\leq C h_{K_1}^{-\frac{1}{2}}\norm{\tilde{\varepsilon}}_{-1,\omega_{E_1}} \norm{R_{E_1}(u_h)}_{0,E_1}.
\shortintertext{Using \eqref{eq:psiE3}, we bound the second term on the right hand side of Equations~\eqref{eq:RE1} as follows}
	\sum_{K_1 \in \omega_{E_1}}  (R_{K_1} (u_h),v_{E_1})_{K_1} 
								&\leq \sum_{K_1 \in \omega_{E_1}}  \norm{R_{K_1} (u_h) 
									}_{0,K_1} \norm {v_{E_1}}_{0,K_1} 
									=\sum_{K_1 \in \omega_{E_1}}  \norm{R_{K_1} (u_h) 
									}_{0,K_1} \norm { \Psi_{E_1}\mathbf{P}R_{E_1}(u_h)}_{0,K_1} \\
								&\leq C \sum_{K_1 
									\in \omega_{E_1}}   h_{E_1}^{\frac{1}{2}} 
									\norm{R_{K_1} (u_h) }_{0,K_1} \norm{ R_{E_1}(u_h)}_{0,E_1},\end{align*}
where $\norm{R_{K_1} (u_h) }_{K_1}$ is bounded already in~\eqref{eq:R_T1LLB}.\\
For the third term of Equations~\eqref{eq:RE1}, we have the following bound thanks to \eqref{eq:psiE2}
\begin{align*}	
	\sum_{K_1 \in \omega_{E_1}} (f- \Pi_0 f,v_{E_1})_{K_1}
								&\leq \sum_{K_1 \in \omega_{E_1}}  \norm{f- \Pi_0 f }_{0,K_1}
									 \norm{v_{E_1}}_{0,K_1}
								=\sum_{K_1 \in \omega_{E_1}}  \norm{f- \Pi_0 f }_{0,K_1}
									 \norm{ \Psi_{E_1}\mathbf{P}R_{E_1}(u_h)}_{0,K_1}\\
								& \leq C \sum_{K_1 \in \omega_{E_1}}
									 h_{E_1}^{\frac{1}{2}}  \norm{f- \Pi_0 f }_{0,K_1}
									 \norm{ R_{E_1}(u_h)}_{0,E_1}.
\end{align*}
Applying again \eqref{eq:psiE3}, we bound the last term of Equations~\eqref{eq:RE1}  as follows
\begin{align*}					
	 (\beta \nabla e,\nabla  v_{E_1} )_{0,\omega_{E_1}} 
								&\leq \sum_{K_1 \in \omega_{E_1}} \beta 
									\norm{\nabla e}_{0,K_1} 
									\norm{\nabla  v_{E_1}}_{0,K_1} 
								\\&=\sum_{K_1 \in \omega_{E_1}} \beta 
									\norm{\nabla e}_{0,K_1} 
									\abs{   \Psi_{E_1}\mathbf{P}R_{E_1}(u_h)}_{1,K_1}
								 \leq C \; \beta  \sum_{K_1 \in \omega_{E_1}}h_{E_1}^{-\frac{1}{2}}
								 	 \norm{\nabla e}_{0,K_1} \norm{ R_{E_1}(u_h)}_{0,E_1}.		
\end{align*}
Adding all results up, multiplying by $\frac{1}{\sqrt{2}} h_{E_1}^{\frac{1}{2}}/\norm{ R_{E_1}(u_h)}_{E_1}$, and using the bound in \eqref{eq:R_T1LLB} yields
\begin{equation}
\label{eq:R_E1LLB}
		\begin{aligned}
			\frac{1}{\sqrt{2}} h_{E_1}^{\frac{1}{2}} \norm{ R_{E_1}(u_h)}_{0,E_1} 
				&\leq C \Bigg[ \beta  \norm{\nabla e}_{0,\omega_{E_1}}+\norm{\tilde{\varepsilon}}_{-1,\omega_{E_1}}+\sum_{K_1 \in \omega_{E_1}} h_{E_1}  
					\norm{f- \Pi_0 f }_{0,K_1}\Bigg].
		\end{aligned}
\end{equation}

For a fixed $K_1$, we add the bound of $h_{K_1}^2 \norm{R_{K_1}(u_h)}_{0,K_1}^2$ considering inequality~\eqref{eq:R_T1LLB} and the bound for $\sum_{E_1 \in \partial K_1 }\frac{1}{2} h_{E_1} \norm{ R_{E_1}(u_h)}_{E_1}^2 $ considering inequality~\eqref{eq:R_E1LLB} and use the arguments in \eqref{eq:norm-1}, \eqref{eq:norm1*}, and \eqref{eq:norm_lambda_extended} to get the desired bound in Definition~\ref{def:eta} for local indicator, $\eta_{K_1}$, as follows:
\begin{equation}
\label{def:eta1}
\boxed{
			\eta_{K_1}^2 \leq \underline{C}_1 \left\lbrace  
						\beta  \norm{\nabla e}_{0,\omega_{K_1}}
						+\norm{\tilde{\varepsilon}}_{[H^1(\omega_{K_1})]^*}
						+\sum_{K_1' \in \overline{\omega}_{K_1}} osc_{K'_1}
						\right\rbrace \cdot \eta_{K_1}. 
			}
\end{equation}
This completes the first part of the proof of Proposition~\ref{prop:LLB} regarding the bound of $\eta_{K_1}$ in~\eqref{eq:lower_a_posterori_T1T2}.

Now, to prove the second inequality in~\eqref{eq:lower_a_posterori_T1T2} , one proceed following the structure of the above proof.

Fix an arbitrary $K_2 \in \mathcal{T}_2 $, and set $ v_{K_2}= R_{K_2}(u_{2h}) \Psi_{K_2}$, where $\Psi_{K_2}$ is the bubble function defined in Lemma~\ref{lma:psiT}. Then, we have
\begin{align}\label{eq:RK2}
		(R_{K_2}(u_{2h}),v_{K_2})_{K_2} & =(\Pi_0 f_3+ \beta_3 \Delta u_{2h} 
												+\lambda_h,v_{K_2})_{K_2}\nonumber\\
											& =(f_3+ \beta_3 \Delta u_{2h} 
												+\lambda_h-(f_3- \Pi_0 f_3),v_{K_2})_{K_2}\\
											& = (\beta_3 \nabla e_2 
												,\nabla v_{K_2})_{K_2} - {}_{H^{-1}(K_2)}\left\langle \varepsilon|_{K_2},v_{K_2} \right\rangle_{H^1_0(K_2)}
												-(f_3 - \Pi_0 f_3 ,v_{K_2})_{K_2}.\nonumber
\end{align}
$K_2$ presents two potential scenarios as depicted in the left of Figure~\ref{fig:possibilities_K2E2}. Due to the bubble function, $v_{K_2}$ vanishes at the boundary of $v_{K_2}$ in both scenarios; i.e. $v_{K_2} \in H^1_0(K_2)$. Consequently, we can restrict the Lagrange multiplier term to the element $v_{K_2}$. The proof proceed like before, therefore, we avoid repetition and we write the result
\begin{figure}[htp]
\centering
 \includegraphics[scale=0.5]{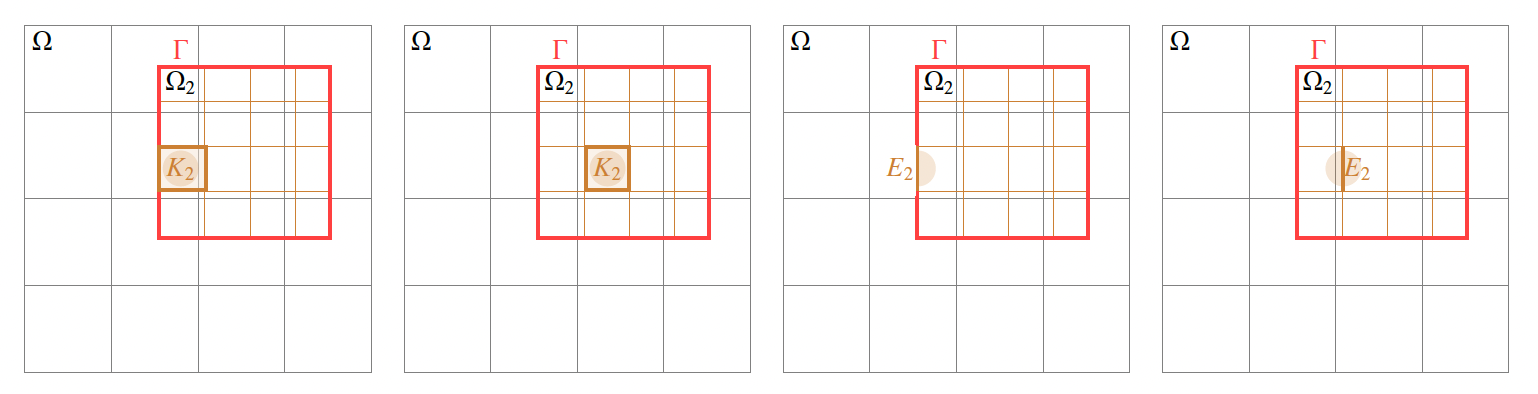} 
\caption{Two different possibilities for the position of the element $K_2$, and  two different possibilities for the position of the edge $E_2$.}\label{fig:possibilities_K2E2}
\end{figure}
%
\begin{equation}
\label{eq:R_T2LLB}
h_{K_2} \norm{R_{K_2}(u_{2h})}_{0,K_2} \leq C \left\lbrace  
\beta_3 \norm{ e_2  }_{1,K_2} +\norm{\tilde{\varepsilon}}_{[H^1(\omega_{K_2})]^*}+h_{K_2}\norm{ (f_3-\Pi_0 f_3)}_{0,K_2} 
 \right\rbrace.
\end{equation}
Regarding the edge $E_2$, we distinguish between two distinct positions of this edge, as illustrated in the rightmost of Figure~\ref{fig:possibilities_K2E2}. It may either be an interior edge or an edge lying on the interface $\Gamma$. For a fixed $E_2$, we set $v_{E_2}= \Psi_{E_2}\mathbf{P}R_{E_2}(u_{2h})$ in both cases.
Due to the bubble function, the support $\omega_{E_2}$ comprises two elements that share the edge if it is interior, and only one element if $E_2$ lies on the interface $\Gamma$. Despite this distinction, both cases are handled similarly for the most part. The primary discrepancy arises during the integration by parts. Specifically, the integral on the boundary of the support differs. In the case of an interior face, it involves the jump of the normal derivative in a fixed direction across the edge, factored by $1/2$. Conversely, if $E_2$ is on $\Gamma$, we only have the integral of the conormal derivative across that edge, which can be treated as a Neumann boundary.

Fix an arbitrary interior edge  $E_2 \in \mathcal{E}_2$ and set $v_{E_2}=  \Psi_{E_2}\mathbf{P}R_{E_2}(u_{2h})$, we obtain:
\begin{align}\label{eq:RE2}
(R_{E_2}(u_{2h}),v_{E_2})_{E_2}&= -   (\beta_3\; \jump{  \frac{\partial u_{2h}}
										{\partial \textbf{n}_2 } }_{E_2}, v_{E_2})_{E_2}\nonumber\\
&=  -{}_{H^{-1}(\omega_{E_2})}\left\langle \varepsilon,v_{E_2}\right\rangle_{H^1_0(\omega_{E_2})}
										-\sum_{K_2 \in \omega_{E_2} } \Bigg[
										( R_{K_2} (u_{2h}),v_{E_2})_{K_2}+(f_3- \Pi_0 f_3,v_{E_2})_{K_2} - (\beta_3 \nabla e_2,
										 \nabla  v_{E_2} )_{K_2}\Bigg].
\end{align}

Following the same structure of the proof, we jump into the result to avoid repetition
\begin{equation}
\label{eq:R_E2LLB2}
	\begin{aligned}
		 h_{E_2}^{\frac{1}{2}} \norm{ R_{E_2}(u_{2h})}_{0,E_2} 
								&\leq C \Bigg[
									 h_{E_2} \norm{f_3- \Pi_0 f_3 }_{0,K_2}
								+\beta_3 \norm{\nabla e_2}_{0,K_2}
								+\norm{\varepsilon}_{[H^{1}( \omega_{K _2})]^*}\Bigg].
	\end{aligned}
\end{equation}
Taking into consideration that $\left\langle\mu ,u|_{\Omega_2}-u_{2}\right\rangle=0$, for all $\mu \in \Lambda$, implies that $u|_{\Omega_2}=u_2$ in $\Omega_2$.
		 				
Hence, we deduce that $u_h|_{\Omega_2}-u_{2h}=-e|_{\Omega_2}+e_2$ and we have for a fixed element $K_2$ in $\mathcal{T}_2$
\begin{align}\label{eq:uhu2hLLB}
	\norm{u_h|_{\Omega_2} -u_{2h}}_{1,K_2} &\leq C \eta_{K_2}  \norm{e|_{\Omega_2} 
								 				- e_{2}}_{1,K_2}.
\end{align}
Adding the bounds in the inequalities~\eqref{eq:R_T2LLB},~\eqref{eq:R_E2LLB2}, and~\eqref{eq:uhu2hLLB}, we get the desired bound in Definition~\eqref{def:eta} for $\eta_{K_2}$ as follows
\begin{align*}
\eta_{K_2}^2 \leq & \underline{C}_2 \left\{\beta_3
								 \norm{\nabla e_2}_{0
								,\omega_{K_2}}
								+\norm{\nabla (e|_{\Omega_2}- e_{2}) }_{0
								,\omega_{K _2}}  
								+\norm{\varepsilon}_{[H^{1}( \omega_{K_2})]^*}
								+\sum_{K_2' \in \omega_{E_2}}
								 osc_{K_2'}\right\} \cdot \eta_{K_2} \\
					\leq &\underline{C}_2\left\{(\beta_3+1
								 \norm{\nabla e_2}_{0
								,\omega_{K_2}}
								+\norm{\nabla e|_{\Omega_2}}_{0
								,\omega_{K_2}} 
								+\norm{\varepsilon}_{[H^{1}( \omega_{K _2})]^*}
								 +\sum_{K_2' \in \omega_{E_2}}
								   osc_{K_2'} \right\} \cdot \eta_{K_2},
\end{align*}
which implies
\begin{equation}
\label{def:eta2}
\boxed{\eta_{K_2}^2 \leq \underline{C}_2  \left\{
							\max{(\beta_3,1)} \norm{\nabla e_2}_{0
							,\omega_{K_2}}
							+\norm{\nabla e|_{\Omega_2}}_{0,\omega_{K_2}}  
							+\norm{\varepsilon}_{[H^{1}( \omega_{K _2})]^*}					
							+\sum_{K_2' \in \omega_{E_2}}  osc_{K_2'}
							 \right\} \cdot \eta_{K_2}.
			}
\end{equation}
This completes the proof of the second inequality in \eqref{eq:lower_a_posterori_T1T2}.
\end{proof}
%
\subsection{Piecewise smooth diffusion coefficients}
\label{se:continuous_beta} 
In this section, we expand upon the investigation conducted in Section~\ref{se:constant_beta} to encompass a more general scenario where $\beta$ and $\beta_2$ are not assumed to be constants in $\Omega$ and $\Omega_2$, respectively. Before proceeding, we introduce the following norm.

For a non-negative integer $s$ and \( p = \infty \), we denote the standard Sobolev space by \( W^{s,\infty}(\omega) \), which is the space of \( s \)-times differentiable and essentially bounded functions \(\phi\) with all derivatives up to order $s$ must be essentially bounded, equipped with the norm \(\norm{\phi}_{s,\infty,\omega}\) defined as follows:
\[
\norm{\phi}_{s,\infty,\omega} = \max_{\abs{\alpha} \leq s} \left[ \esssup_{x \in \omega} \abs{\nabla^\alpha \phi(x)} \right].
\]

Recall that $\Pi_0 $ denotes the $L^2$-projection onto the space of piecewise constant polynomials in the cell $ K_i\in\mathcal{T}_i,  (i=1,2)$.
To simplify the notation, we write
\begin{align*}
osc_{K_1}&= h_{K_1} \norm{f-\Pi_0 f}_{0,K_1} + h_{K_1} \norm{\nabla \cdot((\beta- \Pi_0 \beta) \nabla u_h)}_{0,K_1} + \frac{1}{2} \sum_{E_1 \in \partial K_1} h_{E_1}^{\frac{1}{2}} \norm{\jump{ \dfrac{(\beta-\Pi_0 \beta) \partial u_h}{\partial\textbf{n}_1}}_{E_1}  }_{0,E_1},\\
osc_{K_2}&= h_{K_2} \norm{(f_3-\Pi_0 f_3)}_{0,K_2}+ h_{K_2} \norm{\nabla \cdot((\beta_3- \Pi_0 \beta_3) \nabla u_{2h})}_{0,K_2} 
\\&+ \frac{1}{2} \sum_{E_2 \in \partial K_2\setminus \Gamma} h_{E_2}^{\frac{1}{2}} \norm{\jump{\dfrac{(\beta_3 -\Pi_0 \beta_3) \partial u_{2h}}{\partial \textbf{n}_2}} _{E_2}  }_{0,E_2} 
+ \sum_{E_2 \in \partial K_2\cap \Gamma} h_{E_2}^{\frac{1}{2}} \norm{ \dfrac{(\beta_3 -\Pi_0 \beta_3)\partial u_{2h}}{\partial \textbf{n}_2}  }_{0,E_2}.
\end{align*}
Those terms are the correction terms that are a higher order perturbation of the data. Also referred to as the oscillation terms. One can see that $ h_{K_1} \norm{\nabla \cdot((\beta- \Pi_0 \beta) \nabla u_h)}_{0,K_1}$, $h_{E_1}^{\frac{1}{2}} \norm{\jump{ \dfrac{(\beta-\Pi_0 \beta) \partial u_h}{\partial\textbf{n}_1}}_{E_1}  }_{0,E_1}$, $h_{K_2} \norm{\nabla \cdot((\beta_3- \Pi_0 \beta_3) \nabla u_{2h})}_{0,K_2}$, and $h_{E_2}^{\frac{1}{2}} \norm{\jump{\dfrac{(\beta_3 -\Pi_0 \beta_3) \partial u_{2h}}{\partial \textbf{n}_2}} _{E_2}  }_{0,E_2}$ are bounded as follows
\begin{align*}
h_{K_1} \norm{\nabla \cdot((\beta- \Pi_0 \beta) \nabla u_h)}_{0,K_1}&\leq h_{K_1} \Bigg[\norm{\nabla (\beta- \Pi_0 \beta)}_{0,\infty,K_1} \norm{\nabla u_h}_{0,K_1}+ \norm{\beta- \Pi_0 \beta}_{0,\infty,K_1} \norm{\Delta u_h}_{0,K_1}\Bigg]\\
&\leq C h_{K_1} \norm{\beta- \Pi_0 \beta}_{1,\infty,K_1} \norm{\nabla u_h}_{0,K_1}.
\end{align*}
\begin{align*}
h_{K_2} \norm{\nabla \cdot((\beta_3- \Pi_0 \beta_3) \nabla u_{2h})}_{0,K_2}&\leq h_{K_2} \Bigg[\norm{\nabla (\beta_3- \Pi_0 \beta_3)}_{0,\infty,K_2} \norm{\nabla u_{2h}}_{0,K_2}+ \norm{\beta_3- \Pi_0 \beta_3}_{0,\infty,K_2} \norm{\Delta u_{2h}}_{0,K_2}\Bigg]\\
&\leq C h_{K_2} \norm{\beta_3- \Pi_0 \beta_3}_{1,\infty,K_2} \norm{\nabla u_{2h}}_{0,K_2}.
\end{align*}
\begin{align*}
h_{E_1}^{\frac{1}{2}} \norm{\jump{ \dfrac{(\beta-\Pi_0 \beta) \partial u_h}{\partial\textbf{n}_1}}_{E_1}  }_{0,E_1}& \leq h_{E_1}^{\frac{1}{2}}  \norm{\beta- \Pi_0 \beta}_{0,\infty,E_1} \norm{\jump{ \dfrac{ \partial u_h}{\partial\textbf{n}_1}}_{E_1}}_{0,E_1}
\leq C h_{E_1}^{\frac{1}{2}}  \norm{\beta- \Pi_0 \beta}_{0,\infty,\omega_{E_1}} \norm{ \nabla u_h}_{\frac{1}{2},\omega_{E_1}}\\
&\leq C h_{E_1}^{\frac{1}{2}}  \norm{\beta- \Pi_0 \beta}_{0,\infty,\omega_{E_1}} h_{E_1}^{-\frac{1}{2}} \norm{ \nabla u_h}_{0,\omega_{E_1}}
\leq C h_{E_1}  \norm{\beta- \Pi_0 \beta}_{1,\infty,\omega_{E_1}}  \norm{ \nabla u_h}_{0,\omega_{E_1}}.
\end{align*}
\begin{align*}
h_{E_2}^{\frac{1}{2}} \norm{\jump{ \dfrac{(\beta_3-\Pi_0 \beta_3) \partial u_{2h}}{\partial\textbf{n}_2}}_{E_2}  }_{0,E_2}& \leq h_{E_2}^{\frac{1}{2}}  \norm{\beta_3- \Pi_0 \beta_3}_{0,\infty,E_2} \norm{\jump{ \dfrac{ \partial u_{2h}}{\partial\textbf{n}_2}}_{E_2}}_{0,E_2}
\leq C h_{E_2}^{\frac{1}{2}}  \norm{\beta_3-\Pi_0 \beta_3}_{0,\infty,\omega_{E_2}} \norm{ \nabla u_{2h}}_{\frac{1}{2},\omega_{E_2}}\\
&\leq C h_{E_2}^{\frac{1}{2}}  \norm{\beta_3-\Pi_0 \beta_3}_{0,\infty,\omega_{E_2}} h_{E_2}^{-\frac{1}{2}} \norm{ \nabla u_{2h}}_{0,\omega_{E_2}}
\leq C h_{E_2} \norm{\beta_3-\Pi_0 \beta_3}_{1,\infty,\omega_{E_2}} \norm{ \nabla u_{2h}}_{0,\omega_{E_2}}. 
\end{align*}

 In this case, our definition for the residuals on both elements and edges in $\Omega$, and $\Omega_2$ are modified in the following definition.
 
\begin{definition}
The residual of the element $K_1 \in \mathcal{T}_1$ is
$$
R_{K_1}(u_h)=\nabla \cdot(\Pi_0 \beta\; \nabla u_h) - \tilde{ \lambda}_h+\Pi_0 f,
$$
where $\tilde{ \lambda}_h $   is the extension of $\lambda_h$  defined in Definition~\ref{def:extended_lambdah}. The residual of the edge $E_1 \in \mathcal{E}_1$ is
$$R_{E_1}(u_h)=\begin{cases} 
- \jump{ \Pi_0 \beta \; \dfrac{\partial u_h}{\partial \textbf{n}_1}}_{E_1} &  E_1 \in \partial K_1 \setminus \partial \Omega\\
0 &  E_1 \in \partial \Omega.
\end{cases}
$$
The residuals of the element $K_2 \in \mathcal{T}_2$ and the edge $E_2 \in \mathcal{E}_2$ are 
\begin{align*}
	R_{K_2}(u_{2h})&=\nabla \cdot(\Pi_0 \beta_3\nabla u_{2h}) + \lambda_h+\Pi_0 f_3,\\
	R_{E_2}(u_{2h})&=\begin{cases} 
	- \jump{\Pi_0 \beta_3\; \dfrac{\partial u_{2h}}{\partial \textbf{n}_2}}_{E_2} &  \forall E_2 \in \partial K_2 \setminus \Gamma \\
	- \Pi_0 \beta_3 \; \dfrac{\partial u_{2h}}{\partial \textbf{n}_2} & \forall E_2 \in \partial K_2\cap \Gamma.
	\end{cases}
\end{align*}
where  $\jump{ \cdot } $ is the jump of the normal derivatives  $\dfrac{\partial u_h}{\partial \textbf{n}_1} $ and $\dfrac{\partial u_{2h}}{\partial \textbf{n}_2} $  across the edges $E_1$ and $E_2$, respectively.
\end{definition}
The indicators we propose in this section are weighted combinations of the residuals and we define them as in Definition~\ref{def:eta}.  In order to show reliability, we proceed similarly to what we did for the proof of Proposition~\ref{prop:GUB}. We start by replacing $\beta$ and $\beta_3$ in the old proof by their approximations $\Pi_0 \beta$ and $\Pi_0 \beta_3$ respectively. This introduces extra terms that take into consideration the oscillation of the data $\beta$ and $\beta_3$. In this section, we skip the proof to avoid repetition. Instead, this section aims to point out the differences and to show how to deal with those new terms. The remainder of the proof remains unchanged. \\

Fix $v,v_2$ and let $\mathit{I}_h v, \mathit{I}_{2h} v_2 $ be their the Cl\' ement interpolation, defined in Lemma~\ref{lma:I_h_T1}, into the finite element spaces $V_h$ and $V_{2h}$, respectively.

Utilizing the error equations, the definition of $\tilde{\lambda}_h$ (Definition~\ref{def:extended_lambdah}), the integration by parts, and the Cauchy–Schwarz inequality,  we get
\begin{align*}
	(\beta \nabla e,\nabla v)_{\Omega}
			+\left\langle \varepsilon, v|_{\Omega_2}\right\rangle&=(\beta \nabla e,\nabla (v- \mathit{I}_h v))_{\Omega}
			+\left\langle \varepsilon, (v-\mathit{I}_h v)|_{\Omega_2}\right\rangle
		=(\beta \nabla e,\nabla (v-\mathit{I}_h v))_{\Omega}
			+\left\langle \tilde{\varepsilon}, v-\mathit{I}_h v\right\rangle \\
		&=\sum_{K_1 \in \mathcal{T}_1} \left\lbrace ( R_{K_1}(u_h) ,
			v-\mathit{I}_h v)_{K_1}+(f-\Pi_0 f,
			v-\mathit{I}_h v)_{K_1} 
		+	(\nabla \cdot((\beta- \Pi_0 \beta) \nabla u_h), v-\mathit{I}_h v)_{K_1} 		
			\right\rbrace\\
			&\hspace{0.5in}+\frac{1}{2}\sum_{K_1 \in \mathcal{T}_1}\sum_{E_1 \in\partial 
			K_1\setminus\partial \Omega} \left\lbrace
			( R_{E_1}(u_h) ,v-\mathit{I}_h v)_{E_1}
			+ h_{E_1}^{\frac{1}{2}} (\jump{ \dfrac{(\beta-\Pi_0 \beta) \partial u_h}{\partial \textbf{n}_1}},v-\mathit{I}_h v)_{E_1}\right\rbrace .
\end{align*}
where, for a fixed $K_1$, the new terms can be bounded as follows
\begin{align*}
(\nabla \cdot((\beta- \Pi_0 \beta) \nabla u_h)&, v-\mathit{I}_h v)_{K_1} + \frac{1}{2} \sum_{E_1 \in \partial K_1} h_{E_1}^{\frac{1}{2}} (\jump{ \dfrac{(\beta-\Pi_0 \beta) \partial u_h}{\partial \textbf{n}_1}},v-\mathit{I}_h v)_{E_1}   \\
 &\leq  \norm{\nabla \cdot((\beta- \Pi_0 \beta) \nabla u_h)}_{0,K_1} \norm{v-\mathit{I}_h v}_{0,K_1} + \frac{1}{2} \sum_{E_1 \in \partial K_1} \norm{\jump{ \dfrac{(\beta-\Pi_0 \beta) \partial u_h}{\partial \textbf{n}_1}}_{E_1}   }_{0,E_1}\norm{v-\mathit{I}_h v}_{0,E_1} \\
&\leq h_{K_1} \norm{\nabla \cdot((\beta- \Pi_0 \beta) \nabla u_h)}_{0,K_1} \norm{\nabla v}_{0,\overline{\omega}_{K_1}} + \frac{1}{2} \sum_{E_1 \in \partial K_1} h_{E_1}^{\frac{1}{2}} \norm{\jump{ \dfrac{(\beta-\Pi_0 \beta) \partial u_h}{\partial \textbf{n}_1}}_{E_1}  }_{0,E_1} \norm{\nabla v}_{0,\overline{\omega}_{E_1}}.
\end{align*}
From the second error equation we obtain
\begin{align*}
		(\beta_3\nabla e_2,\nabla v_{2})_{\Omega_2} 
			-\left\langle \varepsilon, v_{2}\right\rangle&= (\beta_3\nabla e_2,\nabla (v_{2}-\mathit{I}_{2h} v_2)_{\Omega_2} 
			-\left\langle \varepsilon, v_{2}-\mathit{I}_{2h} v_2\right\rangle
		= ((\beta_2-\beta)\nabla e_2,\nabla (v_{2}-\mathit{I}_{2h} v_2))_{\Omega_2} 
			-\left\langle \varepsilon, v_{2}-\mathit{I}_{2h} v_2\right\rangle\\
		&=\sum_{K_2 \in \mathcal{T}_2} \left\lbrace ( R_{K_2}(u_{2h}) 
			,v_2-\mathit{I}_{2h} v_2)_{K_2}+(f_3-\Pi_0 f_3,v_2-\mathit{I}_{2h} v_2)_{K_2}
			+ (\nabla \cdot((\beta_3- \Pi_0 \beta_3) \nabla u_{2h}),  v_{2}-\mathit{I}_{2h} v_2)_{K_2}\right\rbrace\\  
			&\hspace{0.5in}+ \frac{1}{2} \sum_{K_2 \in \mathcal{T}_2}
			\sum_{E_2 \in\partial K_2\setminus\Gamma} \left\lbrace
			( R_{E_2}(u_{2h}),v_2-\mathit{I}_{2h} v_2)_{E_2}
			+h_{E_2}^{\frac{1}{2}} (\jump{ \dfrac{(\beta_3-\Pi_0 \beta_3) \partial u_{2h}}{\partial \textbf{n}_2}}, v_{2}-\mathit{I}_{2h} v_2)_{E_2}\right\rbrace\\
			&\hspace{0.5in}+\sum_{K_2 \in \mathcal{T}_2}\sum_{E_2 \in\partial K_2\cap\Gamma}\left\lbrace
			( R_{E_2}(u_{2h}),v_2-\mathit{I}_{2h} v_2)_{E_2}+h_{E_2}^{\frac{1}{2}} (\jump{ \dfrac{(\beta_3-\Pi_0 \beta_3) \partial u_{2h}}{\partial \textbf{n}_2}}, v_{2}-\mathit{I}_{2h} v_2)_{E_2}\right\rbrace,
\end{align*}
where, for a fixed $K_2$, the new terms can be bounded as follows
\begin{align*}
(\nabla \cdot((\beta_3- \Pi_0 \beta_3) \nabla u_{2h})&,  v_{2}-\mathit{I}_{2h} v_2)_{K_2} + \frac{1}{2} \sum_{E_2 \in \partial K_2\setminus \Gamma} h_{E_2}^{\frac{1}{2}} (\jump{ \dfrac{(\beta_3-\Pi_0 \beta_3) \partial u_{2h}}{\partial \textbf{n}_2}}, v_{2}-\mathit{I}_{2h} v_2)_{E_2}   \\
&\hspace{2in}+\sum_{E_2 \in \partial K_2 \cap \Gamma} h_{E_2}^{\frac{1}{2}} ( \dfrac{(\beta_3-\Pi_0 \beta_3) \partial u_{2h}}{\partial \textbf{n}_2}, v_{2}-\mathit{I}_{2h} v_2)_{E_2}\\
 &\leq  \norm{\nabla \cdot((\beta_3- \Pi_0 \beta_3) \nabla u_{2h}}_{0,K_2} \norm{ v_{2}-\mathit{I}_{2h} v_2}_{0,K_2} + \frac{1}{2} \sum_{E_2 \in \partial K_2 \setminus \Gamma} \norm{\jump{ \dfrac{(\beta_3- \Pi_0 \beta_3) \partial u_{2h}}{\partial \textbf{n}_2}}_{E_2}   }_{0,E_2}\norm{ v_{2}-\mathit{I}_{2h} v_2}_{0,E_2} \\
& \hspace{2in}+\sum_{E_2 \in \partial K_2\cap \Gamma} \norm{ \dfrac{(\beta_3- \Pi_0 \beta_3) \partial u_{2h}}{\partial \textbf{n}_2} }_{0,E_2}\norm{ v_{2}-\mathit{I}_{2h} v_2}_{0,E_2}\\
&\leq h_{K_2} \norm{\nabla \cdot((\beta_3- \Pi_0 \beta_3) \nabla u_{2h}}_{0,K_2} \norm{\nabla v_2}_{0,\overline{\omega}_{K_2}} + \frac{1}{2} \sum_{E_2 \in \partial K_2\setminus \Gamma} h_{E_2}^{\frac{1}{2}} \norm{\jump{ \dfrac{(\beta_3- \Pi_0 \beta_3) \partial u_{2h}}{\partial \textbf{n}_2}}_{E_2}  }_{0,E_2} \norm{\nabla v_2}_{0,\overline{\omega}_{E_2}}\\
&\hspace{2in}+\sum_{E_2 \in \partial K_2 \cap \Gamma} h_{E_2}^{\frac{1}{2}} \norm{ \dfrac{(\beta_3- \Pi_0 \beta_3) \partial u_{2h}}{\partial \textbf{n}_2}  }_{0,E_2} \norm{\nabla v_2}_{0,\overline{\omega}_{E_2}}.
\end{align*}
These bounds are incorporated into the oscillation terms $osc_{K_1}$ and $osc_{K_2}$. Therefore, for the global upper bound, the new situation is carried over in Proposition~\ref{prop:GUB}. In essence, adding all those results over all elements and edges yields the following global upper bound
\begin{equation}
\label{eq:upper_a_posterori_continuous_B}
\abs{e}_{1,\Omega}+\norm{e_2}_{1,\Omega_2}+\norm{\varepsilon}_{\Lambda} \leq \overline{C} \left\lbrace 
 \sum_{K_1 \in\mathcal{T}_1} \left( \eta_{K_1}^2 + osc_{K_1}^2 \right) +
\sum_{K_2 \in \mathcal{T}_2} \left( \eta_{K_2}^2 + osc_{K_2}^2 \right) 
\right\rbrace ^{\frac{1}{2}}.
\end{equation}
 
On the other hand, with regard to the local lower bound, a certain level of caution is necessary. We locally replace $\beta$ and $\beta_3$ in each element with their local maxima, defined as follows:
\begin{align*}
\overline{\beta}_{K_1} &= \sup_{x \in K_1} \beta(x),& \overline{\beta}_{K_2} = \sup_{x \in K_2} \beta_{2}(x).
\end{align*} 
Furthermore, as required, we substitute $\beta$ and $\beta_3$ on the sets $\omega_{E_1}$ and $\omega_{E_2}$, respectively, with their local maxima, which are defined as follows
\begin{align*}
\overline{\beta}_{\omega_{E_1}} &= \sup_{x \in \omega_{E_1}} \beta(x),& \overline{\beta}_{\omega_{E_2}} = \sup_{x \in \omega_{E_2}} \beta_{3}(x),\\
\overline{\beta}_{\omega_{K_1}} &= \sup_{x \in \omega_{K_1}} \beta(x),& \overline{\beta}_{\omega_{K_2}} = \sup_{x \in \omega_{K_2}} \beta_{3}(x).
\end{align*}
Fix an arbitrary element  $K_1 \in \mathcal{T}_1$ and set $ v_{K_1}= R_{K_1}(u_h) \Psi_{K_1}$, where $\Psi_{K_1}$ is the bubble function described in Lemma~\ref{lma:psiT}. Inserting $ v_{K_1}$ we are led to
\begin{align*}
		(R_{K_1}(u_h),v_{K_1})_{K_1}&=(\Pi_0 f
											+\nabla \cdot (\Pi_0 \beta \nabla u_h),v_{K_1})_{K_1}
											-(\tilde{\lambda}_h,v_{K_1})_{K_1}\\
										&=(f+ \nabla \cdot(\Pi_0 \beta \nabla u_h)
											-(f- \Pi_0 f),v_{K_1})_{K_1}
											+(\nabla \cdot (\beta \nabla u_h)- \nabla \cdot (\beta \nabla u_h), v_{K_1})_{K_1}
											-(\tilde{\lambda}_h,v_{K_1})_{K_1} \\
										&\leq (\overline{\beta}_{K_1} \nabla e ,\nabla v_{K_1})_{K_1} 
											+ {}_{H^{-1}(K_1)}\left\langle \tilde{ \varepsilon},v_{K_1}\right\rangle_{H^1_0(K_1)} 
											-(f- \Pi_0 f ,v_{K_1})_{K_1}-(\nabla \cdot ((\beta - \Pi_0 \beta) \nabla u_h), v_{K_1})_{K_1}.
	\end{align*}
All terms have been proven before except for the extra term $(\nabla \cdot ((\beta - \Pi_0 \beta) \nabla u_h), v_{K_1})_{K_1}$. Clearly, for this term, we can do the following
\begin{align*}
(\nabla \cdot ((\beta - \Pi_0 \beta) \nabla u_h), v_{K_1})_{K_1} &\leq \norm{\nabla \cdot ((\beta - \Pi_0 \beta) \nabla u_h)}_{0,K_1} \norm{v_{K_1}}_{0,K_1}\\
&\leq \norm{\nabla \cdot ((\beta - \Pi_0 \beta) \nabla u_h)}_{0,K_1} \norm{R_{K_1}(u_h) \Psi_{K_1}}_{0,K_1}
\leq  \norm{\nabla \cdot ((\beta - \Pi_0 \beta) \nabla u_h)}_{0,K_1} \norm{R_{K_1}(u_h)}_{0,K_1}.
\end{align*}
Multiplying by $\dfrac{h_{k_1}}{\norm{R_{K_1}(u_h)}_{0,K_1}}$ we get a term that is a higher order perturbation of $\beta$ and can be fed to $osc_{K_1}$.

Similarly, fix an arbitrary interior edge $E_1 \in \mathcal{E}_1 $ and set $v_{E_1}= \Psi_{E_1} \mathbf{P}R_{E_1}(u_h)$. Inserting $v_{E_1}$ we get\\
\begin{align*}
(R_{E_1}(u_h),v_{E_1})_{E_1}&= -  ( \jump{\Pi_0 \beta \frac{\partial u_h}{\partial \textbf{n}_1 }
									}_{E_1}, v_{E_1})_{E_1}\\
									&\leq - \; {}_{H^{-1}(\omega_{E_1})}\left\langle \tilde{\varepsilon} ,
									 v_{E_1}|_{\Omega_2}\right\rangle_{H^1_0(\omega_{E_1})}
									 +( \jump{(\beta-\Pi_0 \beta)  \frac{\partial u_h}{\partial \textbf{n}_1 }
									}_{E_1}, v_{E_1})_{\omega_{E_1}}\\
									&\hspace{0.5in}+\sum_{K_1 \in \omega_{E_1}}  
									\Bigg[ ( R_{K_1} (u_h),v_{E_1})_{K_1}
									+(f- \Pi_0 f,v_{E_1})_{K_1} 
									-\overline{\beta}_{\omega_{E_1}} (\nabla e,\nabla  v_{E_1} )_{K_1}
									+(\nabla \cdot ((\beta - \Pi_0 \beta) \nabla u_h,v_{E_1})_{K_1} \Bigg].
\end{align*}

Compared to Proposition~\ref{prop:LLB}, we have two extra terms that show up and they can be bounded as follows\\
\begin{align*}
( \jump{(\beta-\Pi_0 \beta)  \frac{\partial u_h}{\partial \textbf{n}_1 }}_{E_1}, v_{E_1})_{\omega_{E_1}}
&\leq \norm{\jump{(\beta-\Pi_0 \beta)  \frac{\partial u_h}{\partial \textbf{n}_1} }}_{0,E_1} \norm{ v_{E_1}}_{ 0,\omega_{E_1}}\\
&=\norm{\jump{(\beta-\Pi_0 \beta)  \frac{\partial u_h}{\partial \textbf{n}_1} }}_{0,E_1} \norm{ \Psi_{E_1} \mathbf{P}R_{E_1}(u_h)}_{0,\omega_{E_1}}
 \leq  h_{E_1}^{\frac{1}{2}} \norm{\jump{(\beta-\Pi_0 \beta)  \frac{\partial u_h}{\partial \textbf{n}_1} }}_{0,E_1} \norm{ R_{E_1}(u_h)}_{ 0,E_1},\\
(\nabla \cdot ((\beta - \Pi_0 \beta) \nabla u_h),v_{E_1})_{K_1}&\leq \norm{\nabla \cdot ((\beta - \Pi_0 \beta) \nabla u_h)}_{0,K_1} \norm{v_{E_1}}_{0,K_1}\\
&\leq \norm{\nabla \cdot ((\beta - \Pi_0 \beta) \nabla u_h)}_{0,K_1} \norm{\Psi_{E_1} \mathbf{P}R_{E_1}(u_h)}_{0,K_1}
\leq  h_{E_1}^{\frac{1}{2}}  \norm{\nabla \cdot ((\beta - \Pi_0 \beta) \nabla u_h)}_{0,K_1} \norm{R_{E_1}(u_h)}_{0,E_1}.
\end{align*}

Multiplying by $\frac{1}{\sqrt{2}}\dfrac{h_{E_1}^{\frac{1}{2}}}{\norm{R_{E_1}(u_h)}_{0,E_1}}$ we get two terms that are higher order perturbation of $\beta$ and can be feed to $osc_{K_1}$.

In a similar manner, we proceed for elements and edges of $\mathcal{T}_2$. Putting everything together we get the following local lower bounds
\begin{equation}
\label{eq:lower_a_posterori_T1T2_continuous_B
}
	\begin{aligned}
\eta_{K_1} &\leq \underline{C}_1 \left\lbrace  
						\overline{\beta}_{\omega_{K_1}}  \norm{\nabla e}_{0,\omega_{K_1}}
						+\norm{\tilde{\varepsilon}}_{[H^1(\omega_{K_1})]^*}
						+\sum_{K'_1 \in \omega_{E_1}} osc_{K'_1}
						\right\rbrace,\\
\eta_{K_2} &\leq \underline{C}_2  \left\lbrace
						\overline{\beta} \norm{\nabla e_2}_{0
						,\omega_{K_2}}
						+\norm{\nabla e|_{\Omega_2}}_{0,\omega_{K_2}}  
						+\norm{\varepsilon}_{[H^{1}(\omega_{K_2})]^*}					
						+\sum_{K_2' \in \omega_{E_2}} osc_{K_2'}
						 \right\rbrace.
	\end{aligned}
\end{equation}
where $\overline{\beta}=\max{\left\lbrace \overline{\beta}_{\omega_{K_2}} ,1\right\rbrace }$.

The collection of all those local indicators on $\mathcal{T}_1$ and $\mathcal{T}_2$ gives the following overall estimators
	\begin{align*}
		\eta^2_1=\sum_{K_1 \in \mathcal{T}_1 } \eta_{K_1}^2,&&
		\eta^2_2=\sum_{ K_2 \in \mathcal{T}_2 } \eta_{K_2}^2.
	\end{align*}
In summary, we have successfully extended the definition of the error estimators to accommodate the general case involving smooth diffusion coefficients. Through rigorous analysis, we have confirmed the reliability and efficiency of this extended case. This comprehensive investigation reinforces the utility and applicability of our approach in accurately estimating errors in problems with smooth diffusion coefficients.
\section{Adaptivity}
\label{se:Adaptivity}
Adaptive methods aim to improve the efficiency and accuracy of numerical simulations by concentrating computational local refinement in regions where they are most needed, such as areas with steep gradients or high solution variability, while reducing computational effort in regions where the solution is smoother or less critical. A posteriori error estimators plays a crucial role in adaptive numerical simulations by providing insights into error distribution across the computational domain which direct the adaptive process.

\begin{figure}[htp]
\begin{minipage}[c]{.51\linewidth}
\begin{center}
		\includegraphics[width=0.4\linewidth]{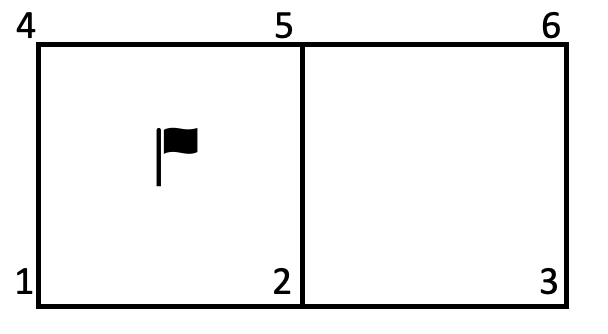} 
		\caption{Left cell marked for refinement, ($Q_1$ FE in 2D).}
		\label{fig:flagged}
\end{center}
\end{minipage}
\begin{minipage}[c]{.51\linewidth}
\begin{center}
		\includegraphics[width=0.4\linewidth]{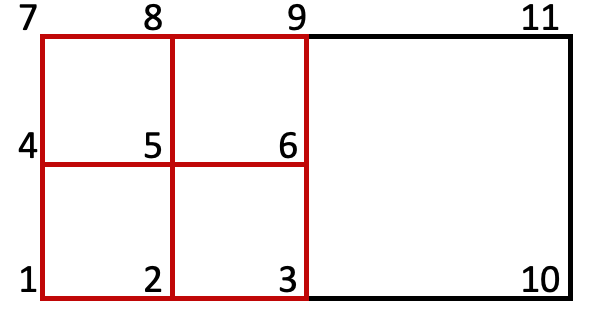}
		\caption{Refinement of left cell, ($Q_1$ FE in 2D).}
		\label{fig:hanging_nodes} 
\end{center}
\end{minipage}
\end{figure}

Adaptive methods are typically guided by the SOLVE-ESTIMATE-MARK-REFINE strategy. Initially, the problem of interest is solved on a coarser mesh. Subsequently, given that the exact error is often unavailable, error indicators are computed locally, and cells with the highest and lowest error indicators are flagged, refer to Figure~\ref{fig:flagged}. Finally, refining these flagged cells with the highest error indicators is executed, and the other case is coarsened. This iterative process continues until the overall error estimators are controlled within a specified tolerance, denoted as $tol$; i.e. $\eta_1 +\eta_2 \leq tol$.

In what follows, we explain what the flagging strategy, the refinement rule, and the coarsening rule are considered in this study.

Let $\alpha_1$ and $\alpha_2$ denote the refinement and coarsening fractions, respectively. Cells are marked for refinement based on different methods; this study adopts the bulk criterion or D\"{o}rfler marking \cite{dorfler1996convergent}, refining a fixed fraction of total error estimators $\eta_1$ and $\eta_2$. These fractions range between $0$ and $1$, where $\alpha_1=0$ signifies uniform refinement, $\alpha_1=1$ implies no refinement, and $\alpha_2=0$ means no coarsening. Here, $\alpha_1=0.6$ and $\alpha_2=0$ are chosen.

We adopt an adaptive approach wherein a parent cell is refined by subdividing it into four children in 2D or eight children in 3D. This subdivision is accomplished by connecting the midpoints of opposite edges and faces. The process for the 2D case is depicted in Figure~\ref{fig:hanging_nodes}. For coarsening, we coarsen the four neighbor cells belonging to the same parent cell. This refinement and coarsening approach results in the appearance of what is known as hanging nodes; which are nodes that appear in the faces of adjacent elements with different refinement levels, such as node 6 in Figure~\ref{fig:hanging_nodes}. 

Handling hanging nodes is essential in adaptively refined meshes. With conforming finite element methods, particularly with quadrilateral meshes, interpolation plays a crucial role in solving the issue raised by those nodes. This involves approximating the values at those nodes based on surrounding nodal values.
The goal is to ensure continuity of the discrete solution, $u_h = \sum_i c_i \phi_i$. Since $\phi_i$ are
 discontinuous at those hanging nodes, we add some constraints to choose $c_i$ to be continuous in order to have a continuous solution $u_h$. 
 In other words, we choose to eliminate the nodal unknowns associated with those nodes by restricting their value to an interpolation of the nodal values of the finite element unknown in its parent nodes. For instance, if we are using a $Q_1$ finite element then we ask that the solution at node 6, in Figure~\ref{fig:hanging_nodes}, be a linear combination of the solutions at node 9 and 3. 
 
Our marking and refinement algorithm is presented in Algorithm~\ref{alg:1}.

 \begin{algorithm}[htp]
 \caption{Flagging strategy}
 \label{alg:1}
 \begin{algorithmic}
\REQUIRE The meshes $\mathcal{T}_1$ and $\mathcal{T}_2$ with the local indicators $\eta_{K_1}$ and $\eta_{K_2}$, and the global estimators $\eta_1$ and $\eta_2$, respectively.\\ The refinement fraction $\alpha_1 \in [0,1]$, and the coarsening fraction $\alpha_2 \in [0,1]$.\\ The set of flagged cells for refinement $\mathcal{\overline{T}}$ and $\mathcal{\overline{T}}_2$, and the set for coarsening cells $\mathcal{\underline{T}}$ and $\mathcal{\underline{T}}_2$.\\ $NC$ is the number of cells in  $\mathcal{T}_1$, and $NC_2$ is the number of cells in  $\mathcal{T}_2$. \\Thresholds for refinement and coarsening $\overline{\theta}_1$, $\overline{\theta}_2$, $\underline{\theta}_1$, and $\underline{\theta}_2$.
\WHILE{$\eta_1+\eta_2 \geq tol$}
\STATE Sort cells in $\mathcal{T}_1$ in a descending local indicator order, $\{\eta_{K_1^1} \geq \eta_{K_1^2} \geq \cdots \geq \eta_{K_1^{NC}}  \}$.
\STATE Sort  cells in $\mathcal{T}_2$ in a descending local indicator order, $\{\eta_{K_2^1} \geq \eta_{K_2^2}\geq \cdots \geq \eta_{K_2^{NC_2}} \}$.
	\FOR{i=1, $\sum_{j=1}^i \eta_{K_1^j}< \alpha_1 \eta_1$, ++i}
		\STATE $\overline{\theta}_1$= i+1.
	\ENDFOR
	\FOR{i=$NC$, $\sum_{j=1}^i \eta_{K_1^j}< \alpha_2 \eta_1$, - -i}
		\STATE $\underline{\theta}_1$=i.
	\ENDFOR
	\FOR{i=1, $\sum_{j=1}^i \eta_{K_2^j}< \alpha_1 \eta_2$, ++i}
		\STATE $\overline{\theta}_2$ =i+1.
	\ENDFOR
	\FOR{i=$NC_2$ , $\sum_{j=1}^i \eta_{K_2^j}< \alpha_2 \eta_2$, - -i}
			\STATE $\underline{\theta}_2$=i.
	\ENDFOR
	\STATE Add $K_1^j, j=1,\cdots, \overline{\theta}_1$ to $\mathcal{\overline{T}}$.
	\STATE Add $K_1^j, j=\underline{\theta}_1, \cdots, NC$ to $\mathcal{\underline{T}}$.
	\STATE Add $K_1^j, j=1,\cdots, \overline{\theta}_2$ to $\mathcal{\overline{T}}_2$.
	\STATE Add $K_1^j, j=\underline{\theta}_2, \cdots, NC_2$ to $\mathcal{\underline{T}}_2$.
	\STATE Refine $\mathcal{\overline{T}}$ and $\mathcal{\overline{T}}_2$. Coarsen $\mathcal{\underline{T}}$ and $\mathcal{\underline{T}}_2$.
	\STATE Update the meshes $\mathcal{T}_1$ and $\mathcal{T}_2$ and empty the sets $\mathcal{\overline{T}}$, $\mathcal{\overline{T}}_2$,$\mathcal{\underline{T}}$ and $\mathcal{\underline{T}}_2$.
\ENDWHILE
 \end{algorithmic}
 \end{algorithm}
\section{Numerical results}
\label{se:numerical_results}
In this section, we employ the finite element library, \texttt{deal.II}~\cite{dealII95}, to solve the examples under study. The system is tackled utilizing a GMRES solver, with a triangular-blocked preconditioner as detailed in \cite{boffi2023parallel}.

\begin{figure}[htp]
\begin{minipage}[c]{.51\linewidth}
\begin{center}
		\includegraphics[width=0.9\linewidth]{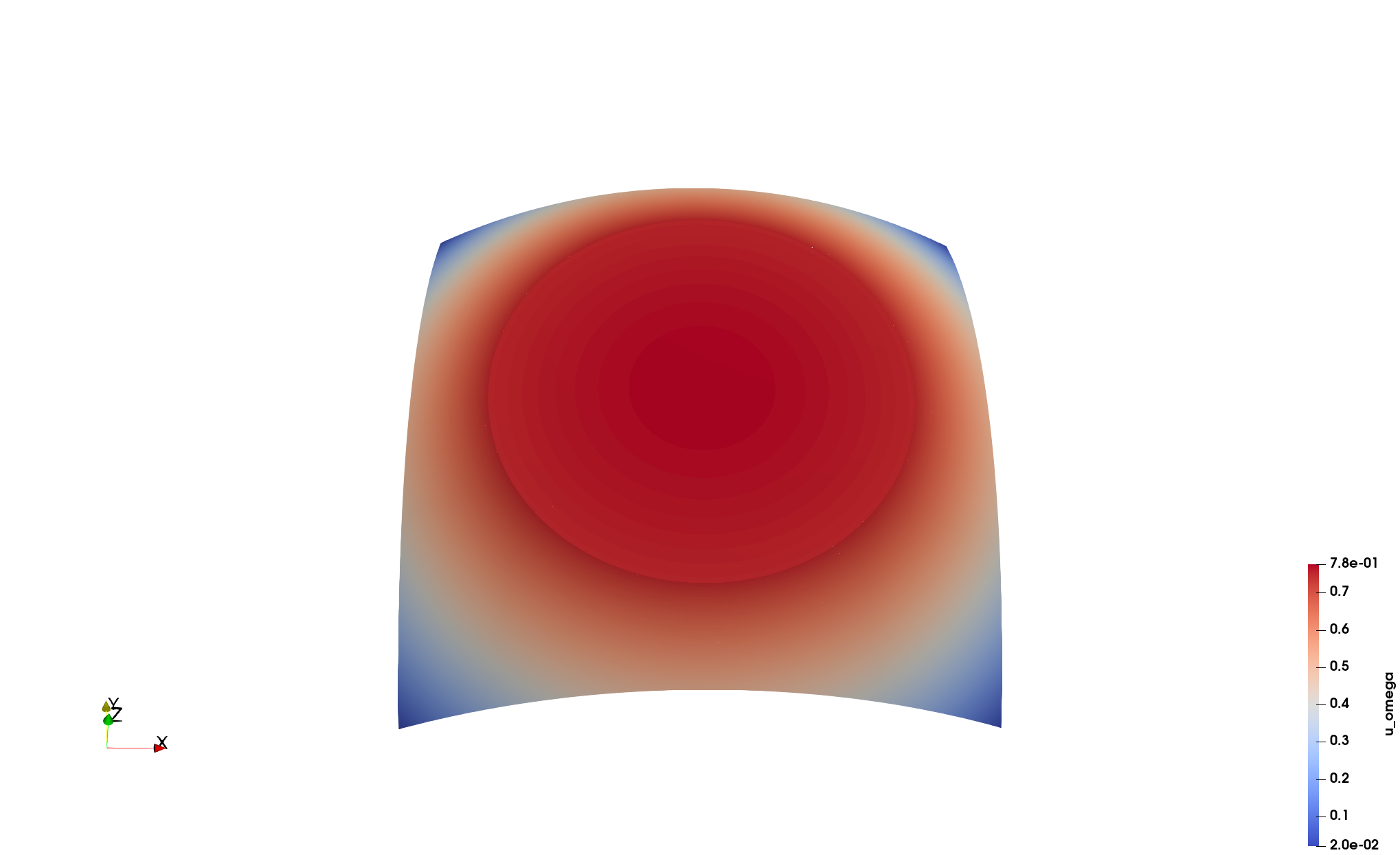} 
		\caption{Solution: immersed circle shape.}
		\label{fig:circle}
\end{center}
\end{minipage}
\begin{minipage}[c]{.51\linewidth}
\begin{center}
		\includegraphics[width=0.9\linewidth]{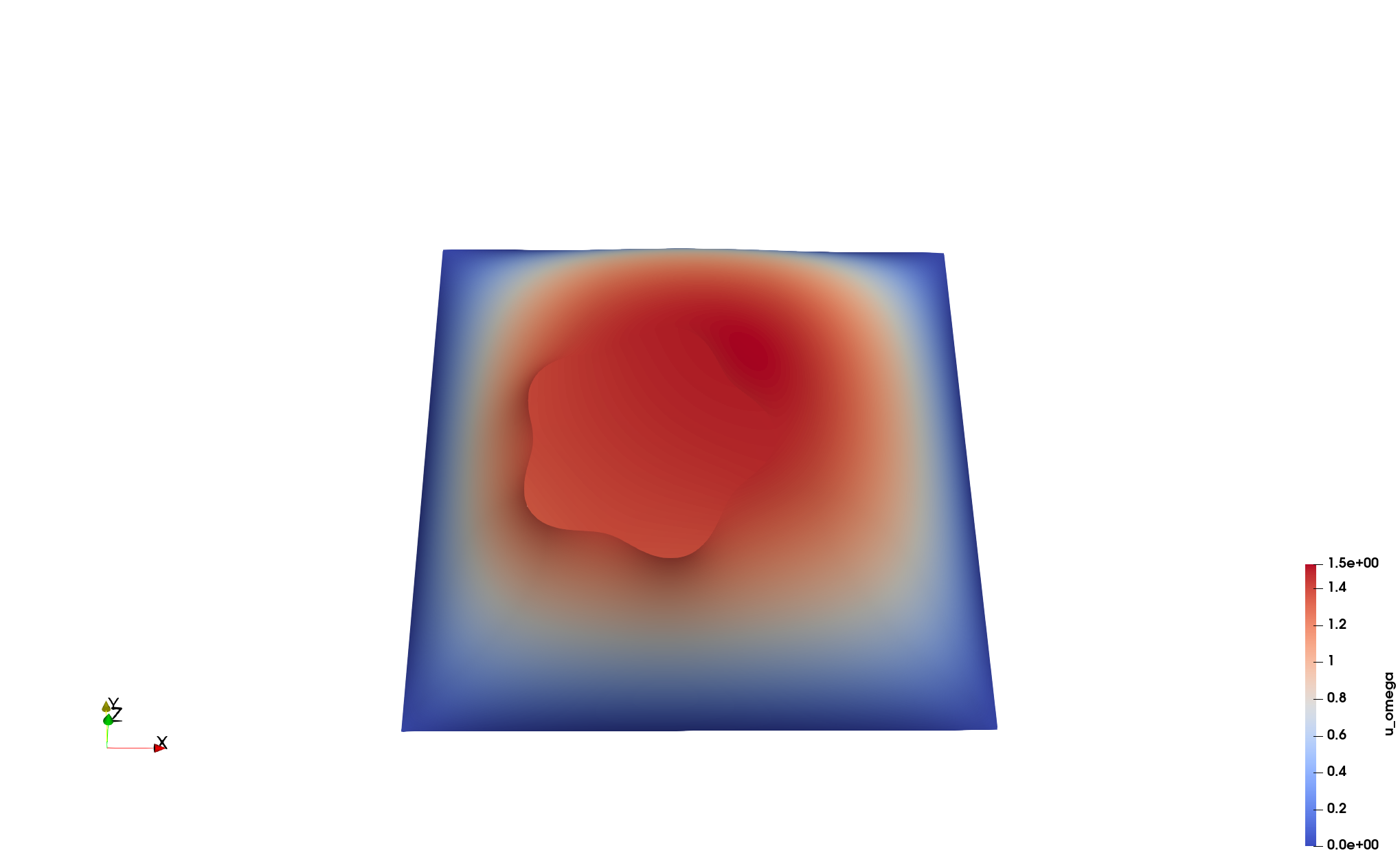}
		\caption{Solution: immersed flower shape.}
		\label{fig:flower} 
\end{center}
\end{minipage}
\begin{minipage}[c]{.51\linewidth}
\begin{center}
		\includegraphics[width=0.9\linewidth]{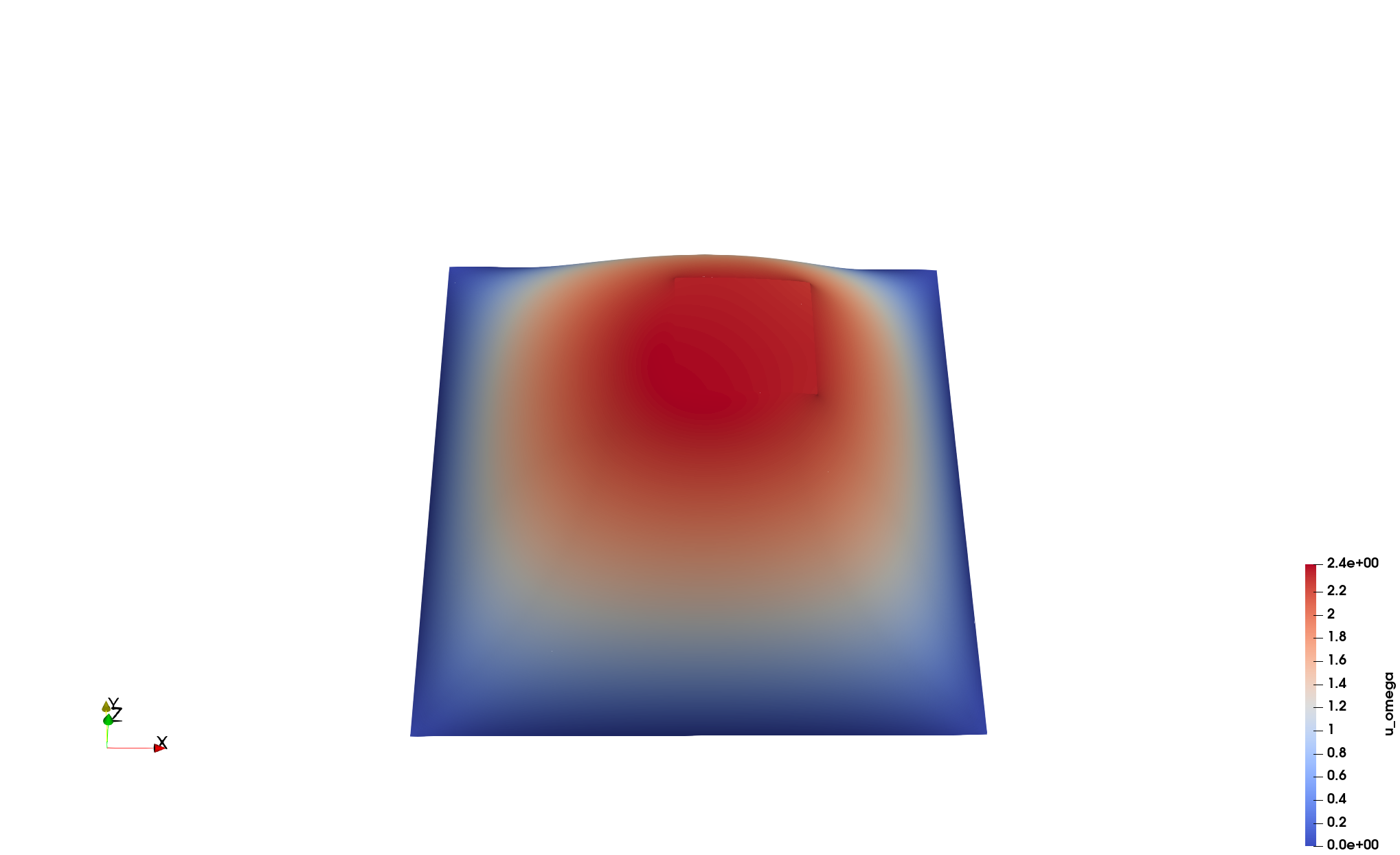} 
		\caption{Solution: immersed square shape.}
		\label{fig:square}
\end{center}
\end{minipage}
\begin{minipage}[c]{.51\linewidth}
\begin{center}
		\includegraphics[width=0.9\linewidth]{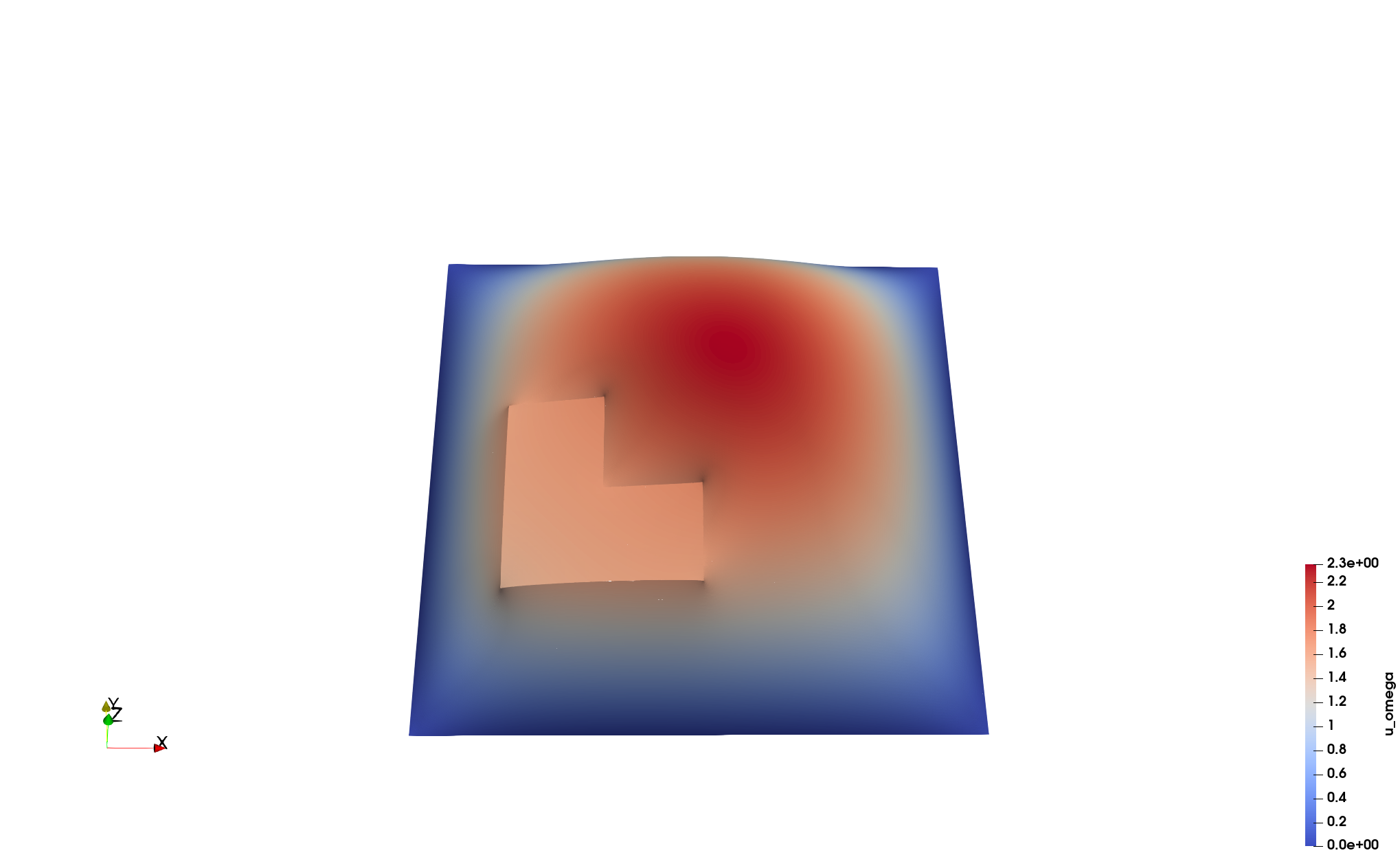}
		\caption{Solution: immersed L-shape.}
		\label{fig:L} 
\end{center}
\end{minipage}
\end{figure}

We showcase four examples that not only illustrate the reliability and efficiency of the estimator, supporting theoretical proofs, but also demonstrate the capability to govern a self-adaptive algorithm and show the robustness of our approach concerning different immersed geometries. These examples include:
\begin{itemize}
\item A rectangular mesh incorporating a square immersed shape mesh, where the domains are defined as $\Omega=[0,6]^2$ and $\Omega_2=[e, 1+\pi]^2$.
\item Another rectangular mesh featuring an immersed L-shape mesh, defined within $\Omega=[0,6]^2$ and $\Omega_2=[1,3]^2\setminus[2,3]^2$.
\item A rectangular mesh combined with a circular immersed shape mesh, spanning $\Omega=[-1.4,1.4]^2$ and $\Omega_2 =\mathcal{B}((0,0),1)$.
\item Lastly, a rectangular mesh involving a flower-shaped immersed mesh, where $\Omega=[-2,3]^2$ and $\Omega_2$ represents a flower centered at (0,0) with a radius $1+0.1\cos(5 \, \theta)$, for $ 0\leq \theta <2 \, \pi$.
\end{itemize}

Among these scenarios, the example with an immersed circle shape presents exact analytic solutions. Specifically, for $\beta=1$ and $\beta_2=10$, the solution is given by:
\begin{align*}
\begin{cases}
\dfrac{4-x^2-y^2}{4}& \textit{in } \Omega_1,\\
\dfrac{31-x^2-y^2}{40}& \textit{in } \Omega_2.
\end{cases}
\end{align*}
Similarly, for $\beta=1$ and $\beta_2=1000$, the solution becomes:
\begin{align*}
\begin{cases}
\dfrac{4-x^2-y^2}{4}& \textit{in } \Omega_1,\\
\dfrac{3001-x^2-y^2}{4000}& \textit{in } \Omega_2.
\end{cases}
\end{align*}
When the coefficients violate the discrete Elker necessity condition, i.e. when $\beta_2<\beta$, the solution of the case $\beta=10$ and $\beta_2=1$ is altered to:
\begin{align*}
\begin{cases}
\dfrac{4-x^2-y^2}{40}& \textit{in } \Omega_1,\\
\dfrac{13-10x^2-10y^2}{40}& \textit{in } \Omega_2.
\end{cases}
\end{align*}

For the remaining examples lacking known solutions, we employed reference solutions computed using over 16 million degrees of freedom to compute errors. All solutions were plotted in Figures~\ref{fig:circle}-\ref{fig:L} for the case $\beta=1$ and $\beta_2=10$. 

Let us denote the refinements as levels. The coarsest mesh would then be at level 0. We tested different ratios of mesh sizes in the initial level, level 0, including $\dfrac{h}{h_2}\approx 1,\dfrac{h}{h_2}< 1$, and $\dfrac{h}{h_2}> 1$.

\begin{figure}[htp]
\begin{minipage}[c]{.33\linewidth}
\begin{center}
		\includegraphics[width=1\linewidth]{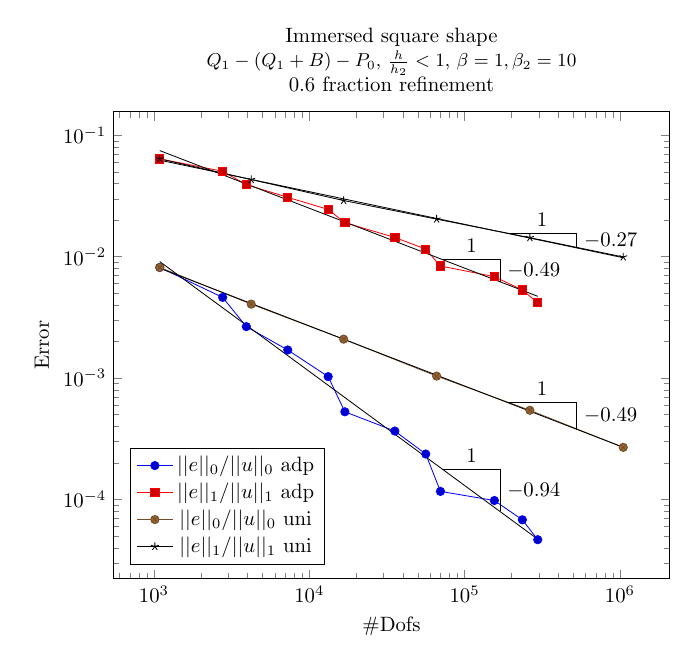} 
\end{center}
\end{minipage}
\begin{minipage}[c]{.33\linewidth}
\begin{center}
		\includegraphics[width=1\linewidth]{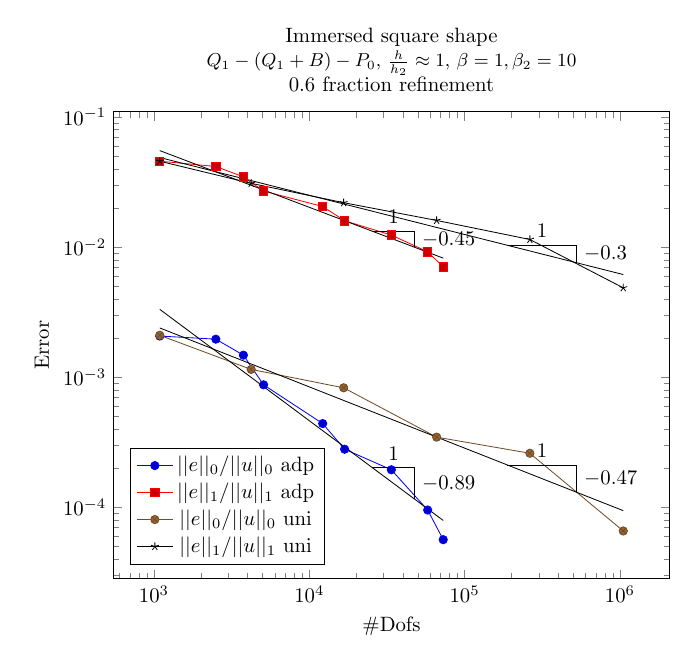}
\end{center}
\end{minipage}
\begin{minipage}[c]{.33\linewidth}
\begin{center}
		\includegraphics[width=1\linewidth]{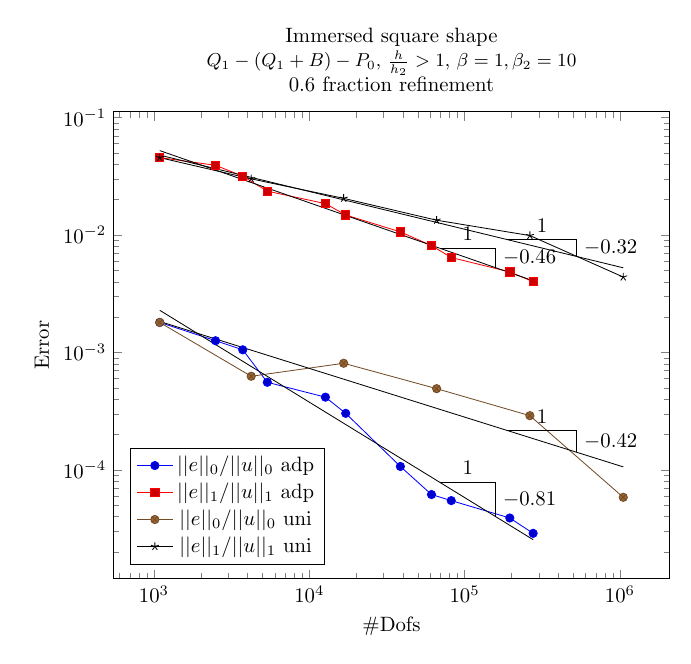} 
\end{center}
\end{minipage}
		\caption{Error of u: immersed square shape.}
		\label{fig:rate_square}
\end{figure}

\begin{figure}[htp]
\begin{minipage}[c]{.33\linewidth}
\begin{center}
		\includegraphics[width=1\linewidth]{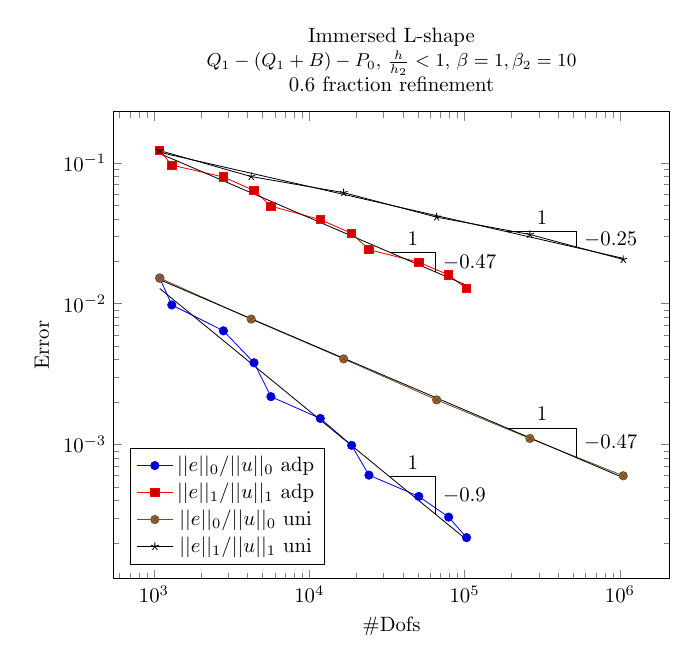} 
\end{center}
\end{minipage}
\begin{minipage}[c]{.33\linewidth}
\begin{center}
		\includegraphics[width=1\linewidth]{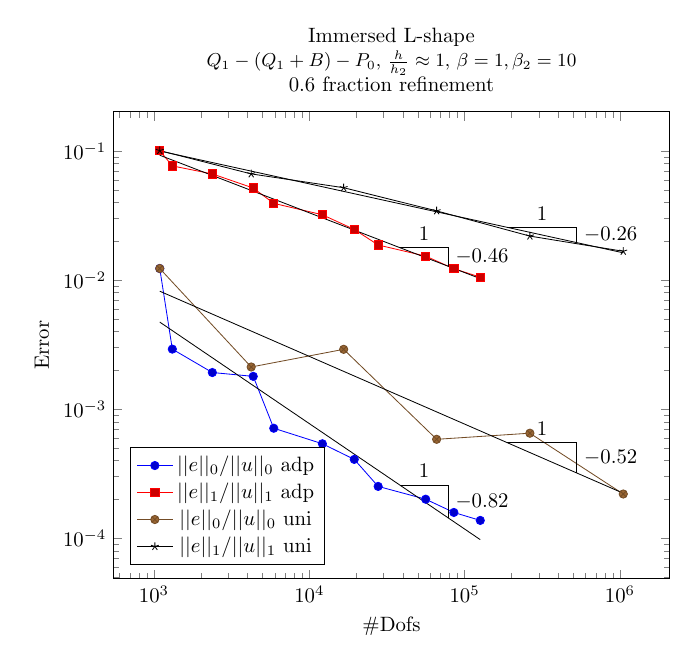}
\end{center}
\end{minipage}
\begin{minipage}[c]{.33\linewidth}
\begin{center}
		\includegraphics[width=1\linewidth]{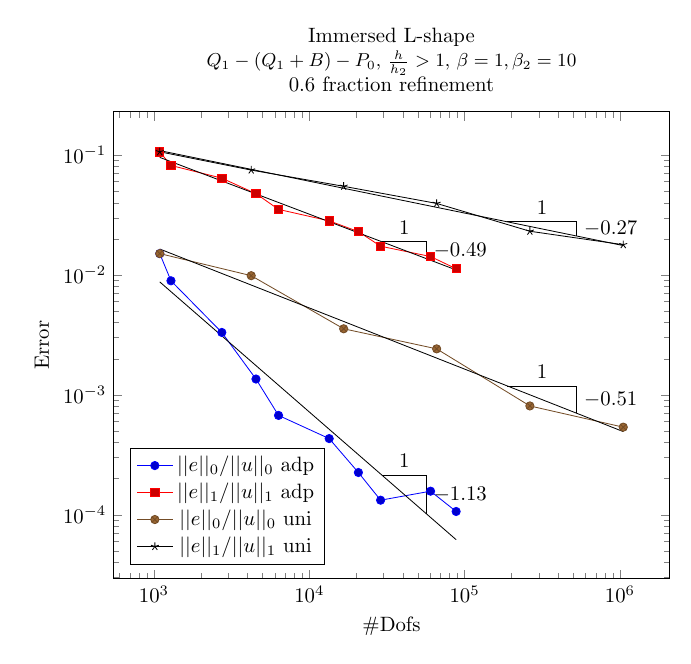} 
\end{center}
\end{minipage}
		\caption{Error of u: immersed L-shape.}
		\label{fig:rate_L}
\end{figure}

\begin{figure}[htp]
\begin{minipage}[c]{.33\linewidth}
\begin{center}
		\includegraphics[width=1\linewidth]{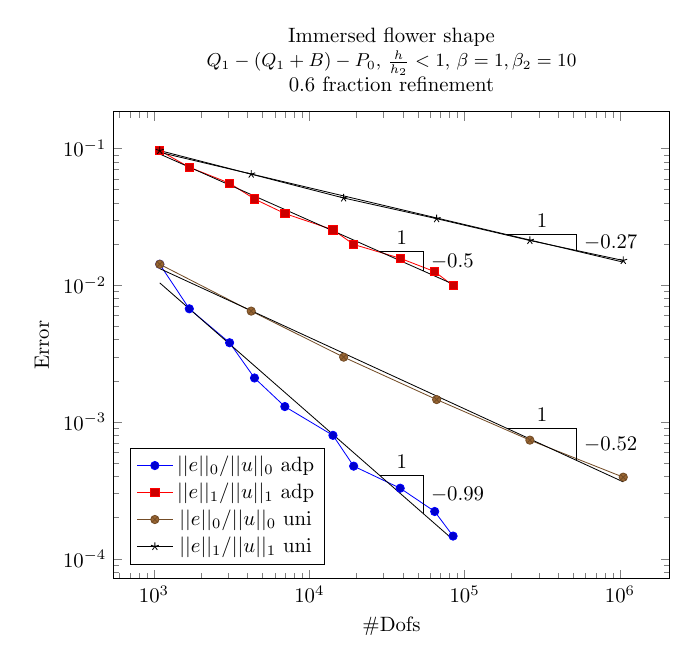} 
\end{center}
\end{minipage}
\begin{minipage}[c]{.33\linewidth}
\begin{center}
		\includegraphics[width=1\linewidth]{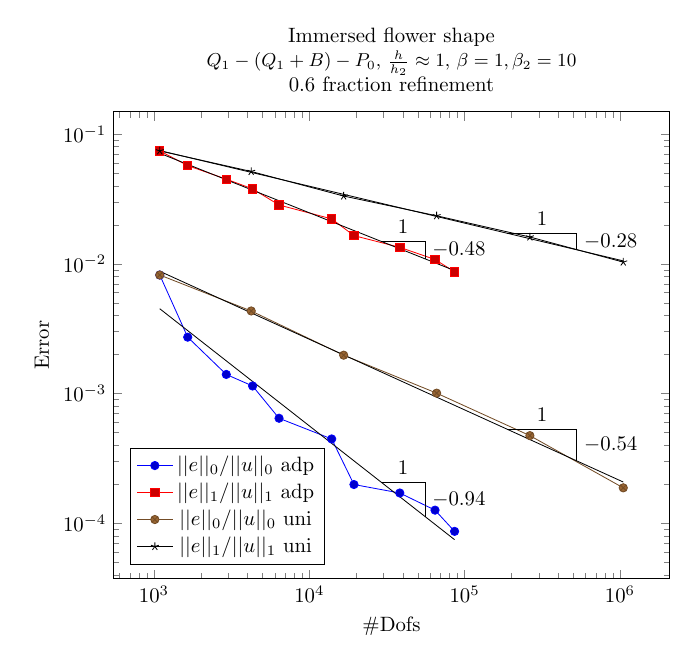}
\end{center}
\end{minipage}
\begin{minipage}[c]{.33\linewidth}
\begin{center}
		\includegraphics[width=1\linewidth]{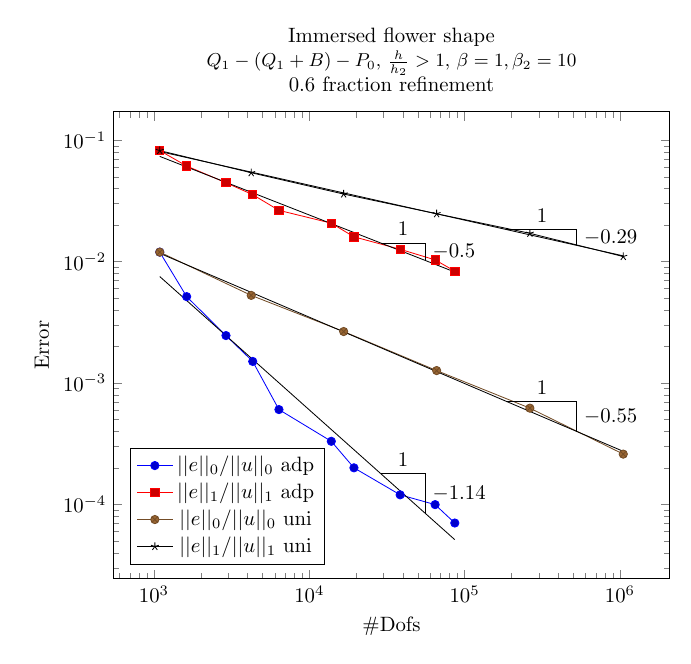} 
\end{center}
\end{minipage}
		\caption{Error of u: immersed flower shape.}
		\label{fig:rate_flower}
\end{figure}

\begin{figure}[htp]
\begin{minipage}[c]{.33\linewidth}
\begin{center}
		\includegraphics[width=1\linewidth]{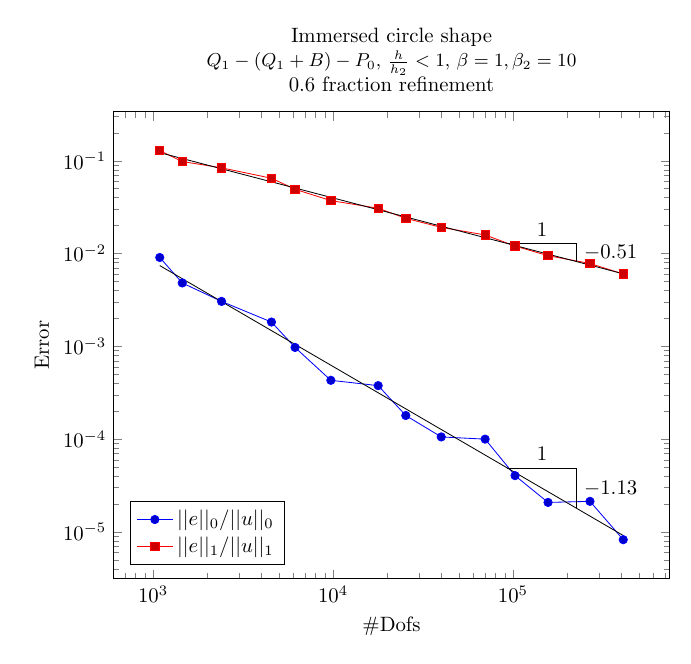} 
\end{center}
\end{minipage}
\begin{minipage}[c]{.33\linewidth}
\begin{center}
		\includegraphics[width=1\linewidth]{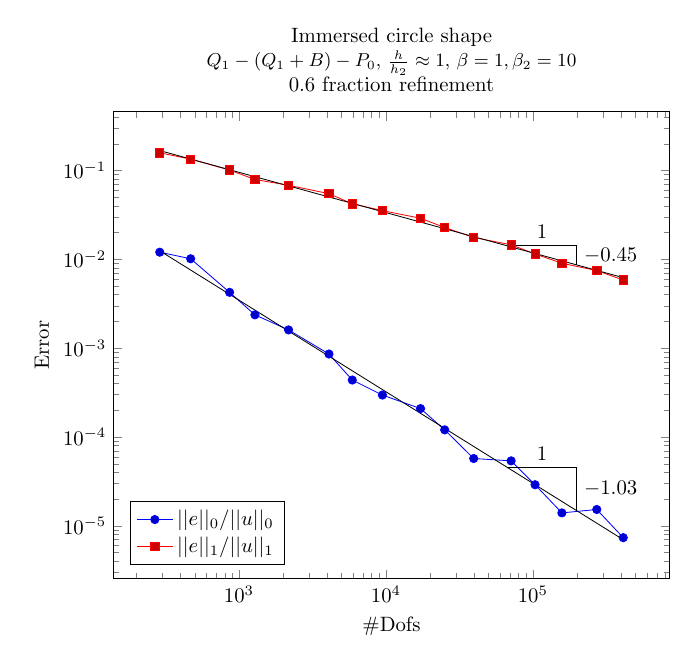}
\end{center}
\end{minipage}
\begin{minipage}[c]{.33\linewidth}
\begin{center}
		\includegraphics[width=1\linewidth]{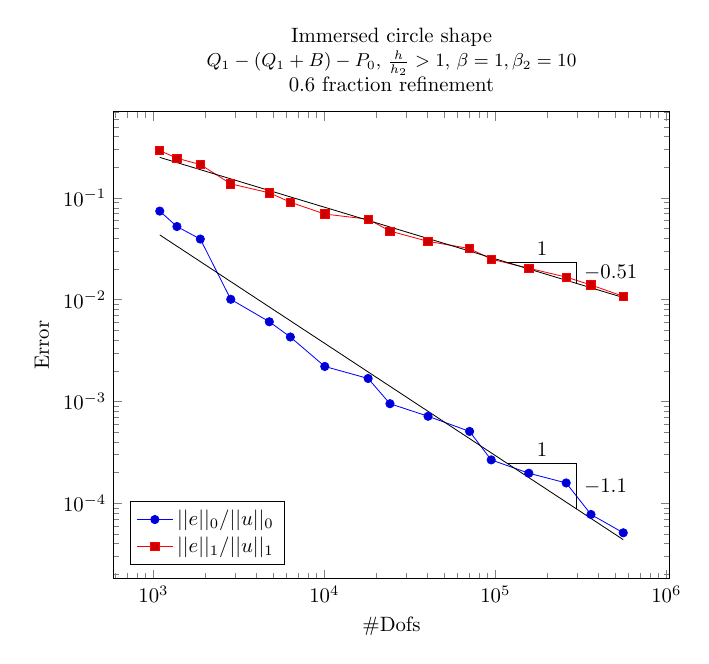} 
\end{center}
\end{minipage}
		\caption{Error of u: immersed circle shape.}
		\label{fig:rate_circle_u}
\end{figure}

\begin{figure}[htp]
\begin{minipage}[c]{.33\linewidth}
\begin{center}
		\includegraphics[width=1\linewidth]{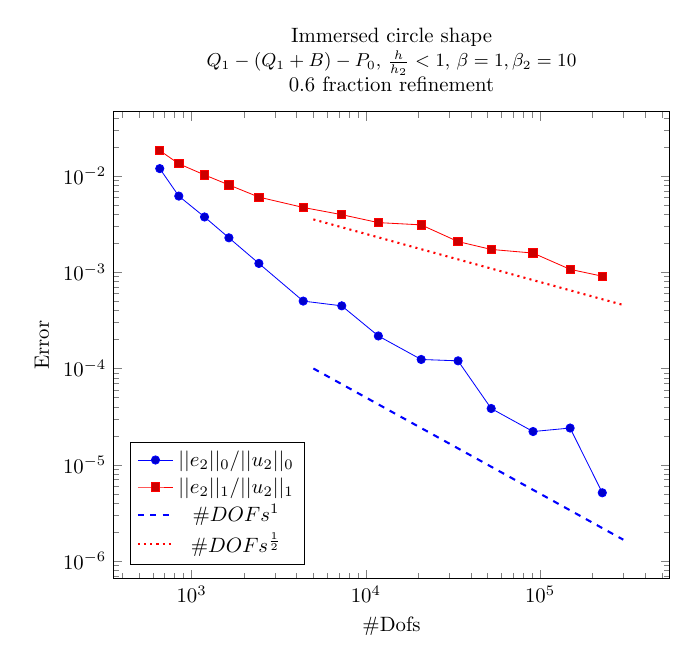} 
\end{center}
\end{minipage}
\begin{minipage}[c]{.33\linewidth}
\begin{center}
		\includegraphics[width=1\linewidth]{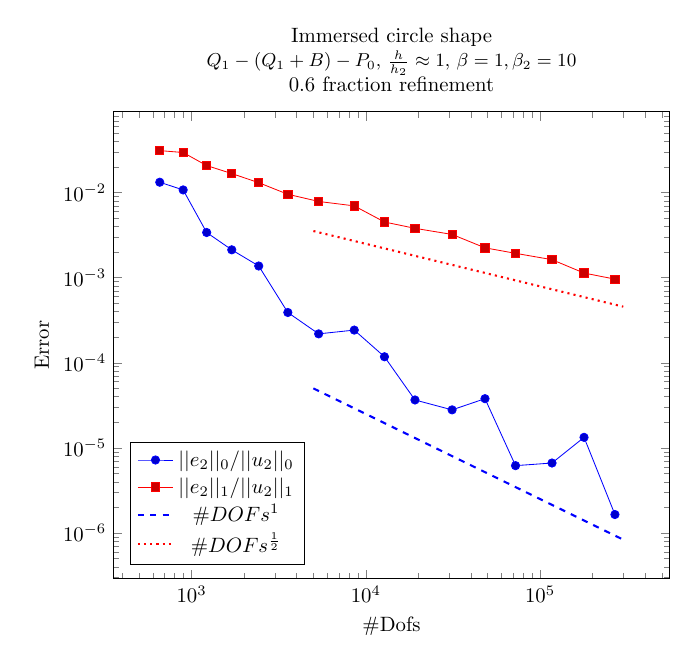}
\end{center}
\end{minipage}
\begin{minipage}[c]{.33\linewidth}
\begin{center}
		\includegraphics[width=1\linewidth]{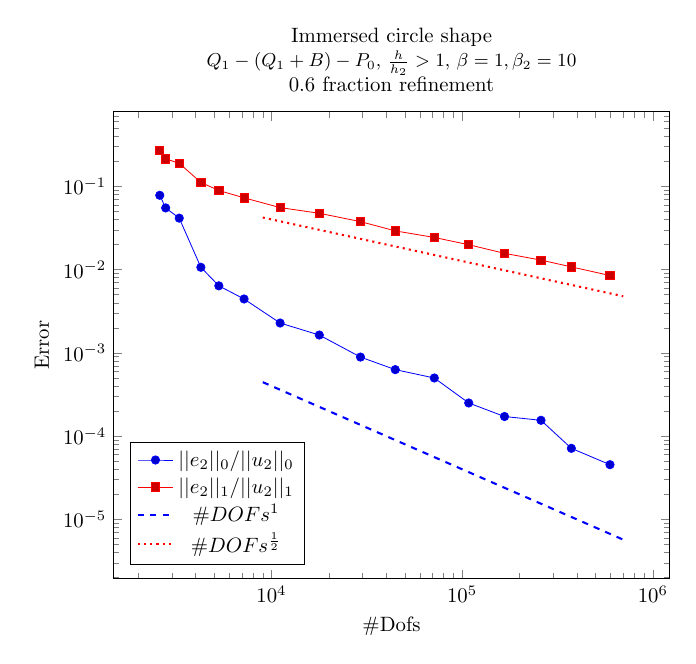} 
\end{center}
\end{minipage}
		\caption{Error of $u_2$: immersed circle shape.}
		\label{fig:rate_circle_u2}
\end{figure}

Our examples align with the geometries in \cite{Najwa2022elliptic}, where uniform refinement is applied. In this work, we apply an adaptive refinement algorithm, introduced in Section~\ref{se:Adaptivity}, to study convergence rates, reliability, and efficiency of the estimators. Remarkably, with our proposed adaptive algorithm, we were able to outperform the results in \cite{Najwa2022elliptic}, achieving optimal convergence rates even with less regular solutions and with the presence of singularities in the solutions. Specifically, $\norm{e}_0$ converges with order $O(h^2)$ and $\norm{e}_1$ converges with order $O(h)$, as depicted in Figures~\ref{fig:rate_square}-\ref{fig:rate_flower}. 
In particular, for the case of the immersed circle shape example, illustrated by Figure~\ref{fig:rate_circle_u}, $\norm{e}_0$ and $\norm{e}_1$ also exhibit the optimal convergence rate when utilizing the proposed adaptive algorithm. Notably, $\norm{e_2}_0$ and $\norm{e_2}_1$ initially display some pre-asymptotic behavior but converge to the correct rate with further refinement (Figure~\ref{fig:rate_circle_u2}). Super-convergence at the first few levels of adaptive refinements is a common behavior.

\begin{figure}[htp]
\begin{minipage}[c]{.33\linewidth}
\begin{center}
		\includegraphics[width=1\linewidth]{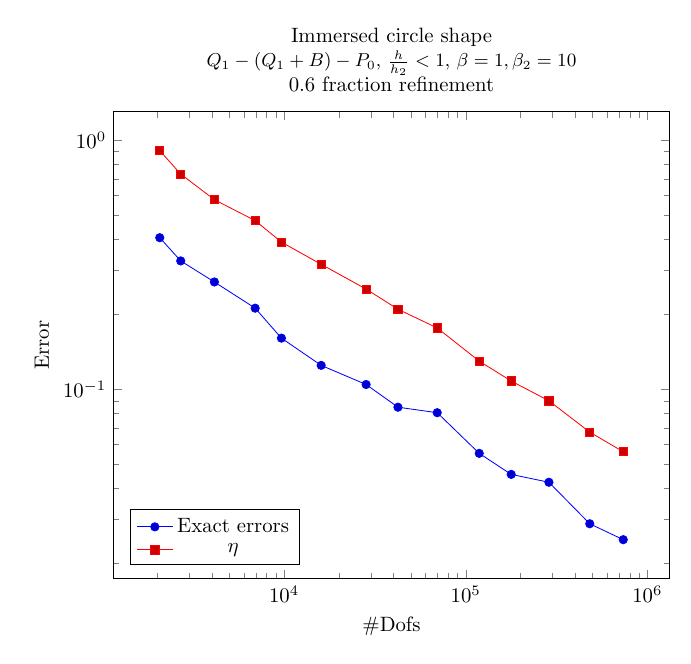} 
\end{center}
\end{minipage}
\begin{minipage}[c]{.33\linewidth}
\begin{center}
		\includegraphics[width=1\linewidth]{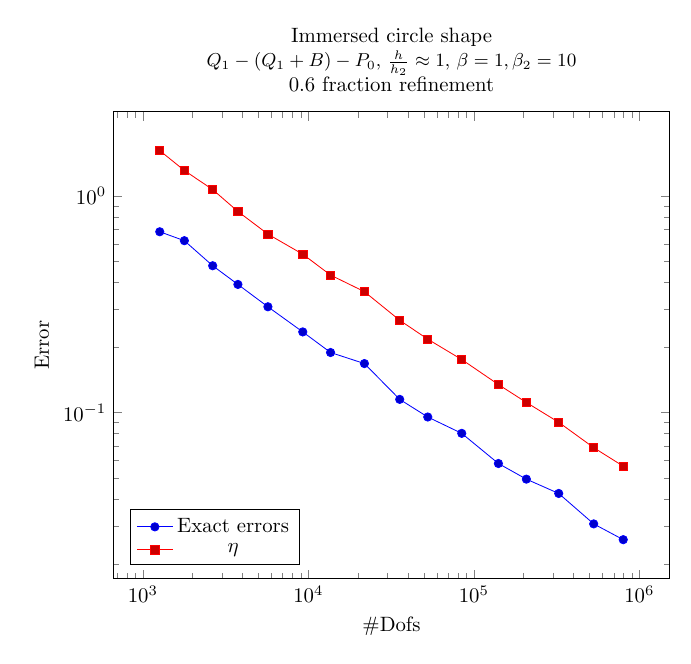}
\end{center}
\end{minipage}
\begin{minipage}[c]{.33\linewidth}
\begin{center}
		\includegraphics[width=1\linewidth]{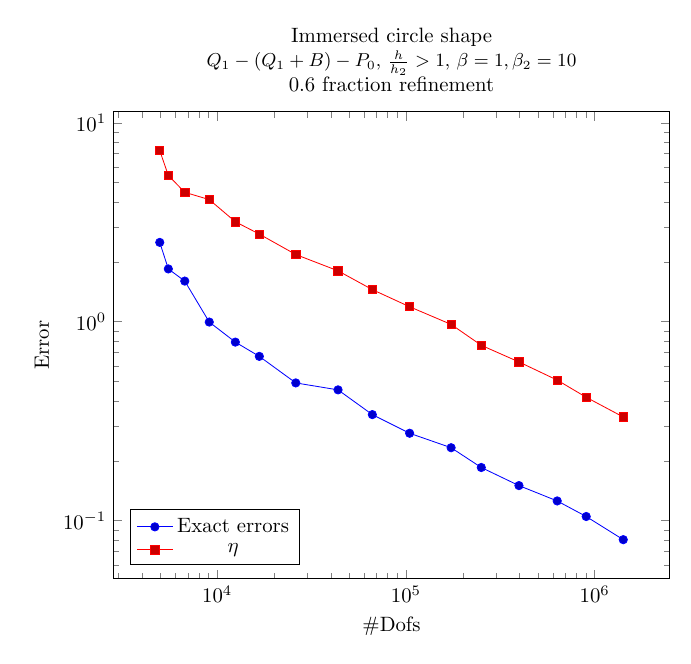} 
\end{center}
\end{minipage}
		\caption{Errors vs estimators: immersed circle shape.}
		\label{fig:rate_circle_all}
\end{figure}

\begin{figure}[htp]
\begin{minipage}[c]{.33\linewidth}
\begin{center}
		\includegraphics[width=1\linewidth]{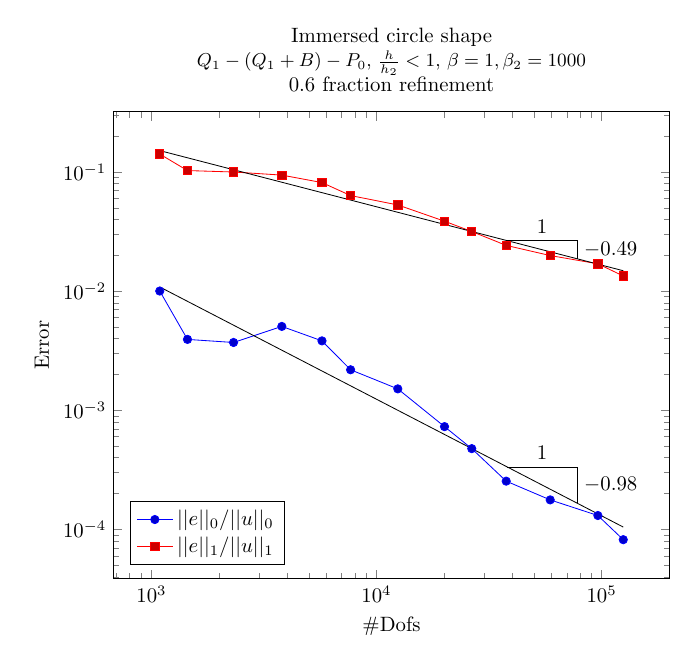} 
\end{center}
\end{minipage}
\begin{minipage}[c]{.33\linewidth}
\begin{center}
		\includegraphics[width=1\linewidth]{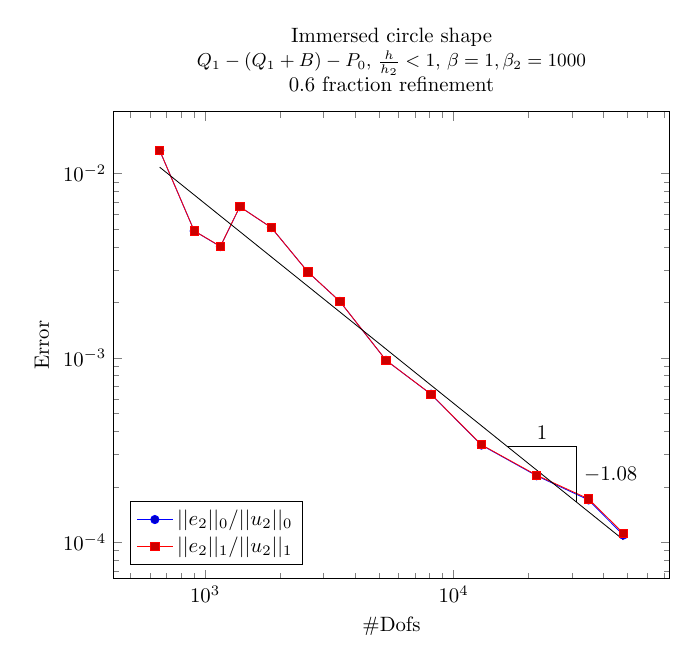}
\end{center}
\end{minipage}
\begin{minipage}[c]{.33\linewidth}
\begin{center}
		\includegraphics[width=1\linewidth]{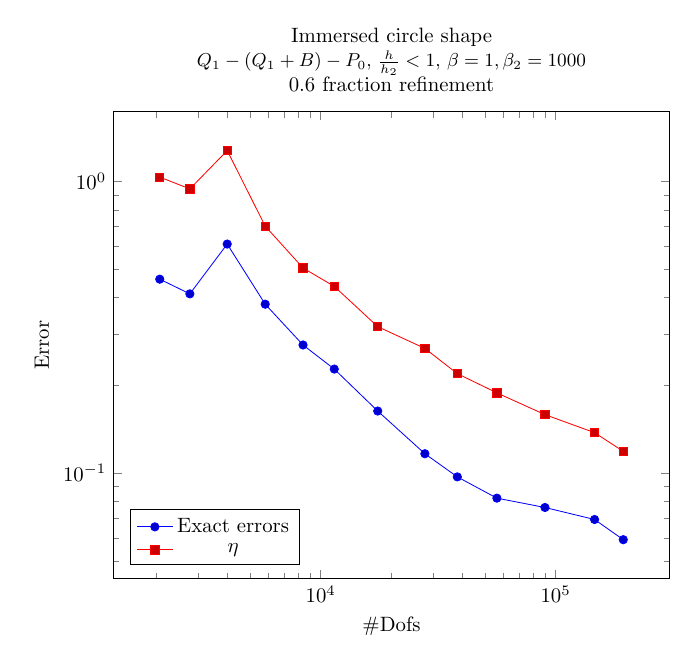} 
\end{center}
\end{minipage}
		\caption{Error of $u$, $u_2$, and errors vs estimators: immersed circle shape.}
		\label{fig:circle1}
\end{figure}

\begin{figure}[htp]
\begin{minipage}[c]{.33\linewidth}
\begin{center}
		\includegraphics[width=1\linewidth]{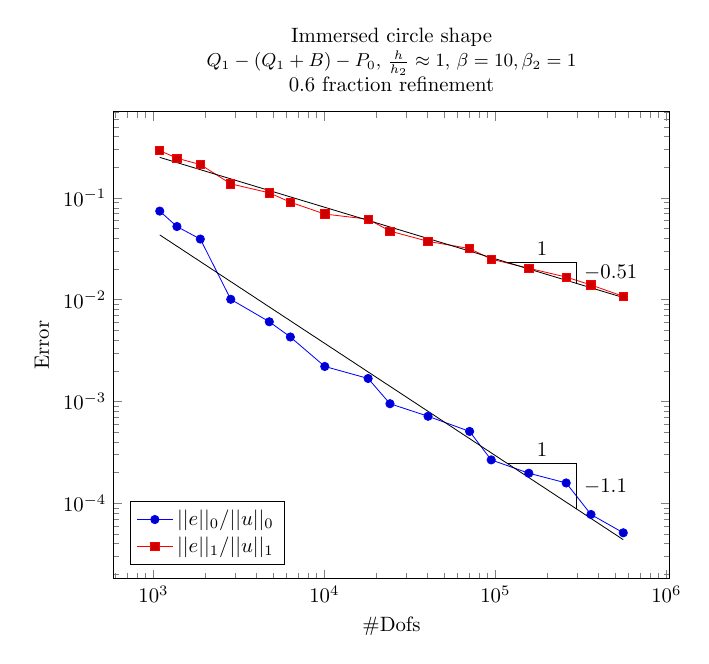} 
\end{center}
\end{minipage}
\begin{minipage}[c]{.33\linewidth}
\begin{center}
		\includegraphics[width=1\linewidth]{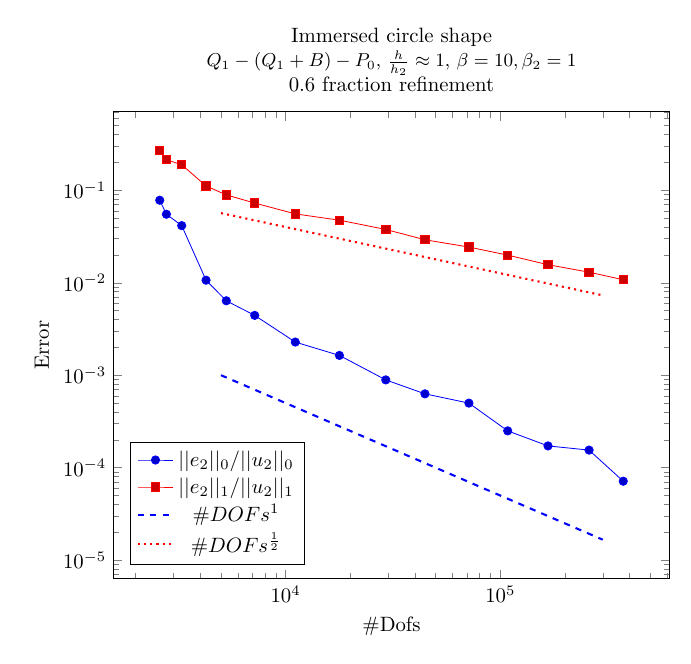}
\end{center}
\end{minipage}
\begin{minipage}[c]{.33\linewidth}
\begin{center}
		\includegraphics[width=1\linewidth]{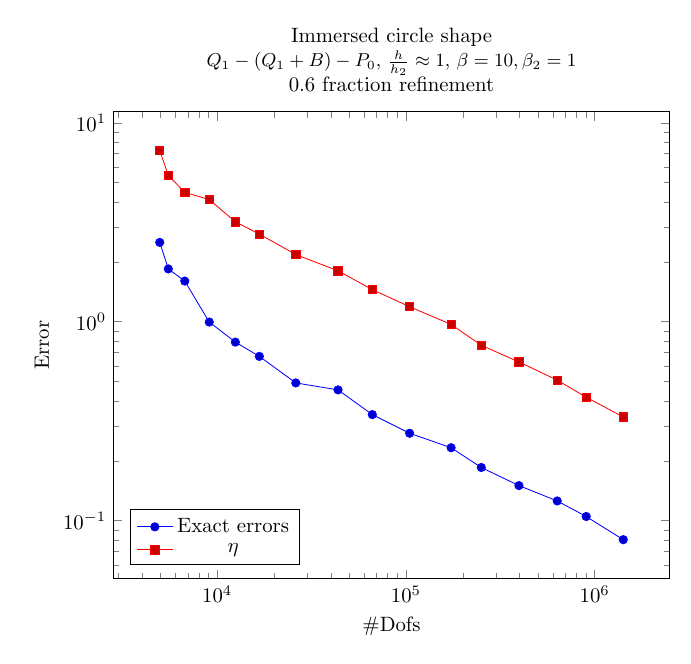} 
\end{center}
\end{minipage}
		\caption{Error of $u$, $u_2$, and errors vs estimators in the example with immersed circle shape.}
		\label{fig:circle2}
\end{figure}

Since exact solutions are known for the example with an immersed circle shape, we tested the efficiency of our estimator by comparing the convergence of the sum of all error norms and the sum of all local indicators in both meshes. As illustrated in Figure~\ref{fig:rate_circle_all}, the estimator demonstrates equivalence to the exact errors with a stable efficiency index. We also  tested our approach for a higher jump of the coefficients ( $\beta=1$ and $\beta_2=1000$) and for the case that violates the Elker necessary condition ( $\beta=10$ and $\beta_2=1$). As depicted in Figures~\ref{fig:circle1} and \ref{fig:circle2}, the error converges optimally and the indicators are both reliable and efficient.

\begin{figure}[htp]
\begin{minipage}[c]{.5\linewidth}
\begin{center}
		\includegraphics[width=1\linewidth]{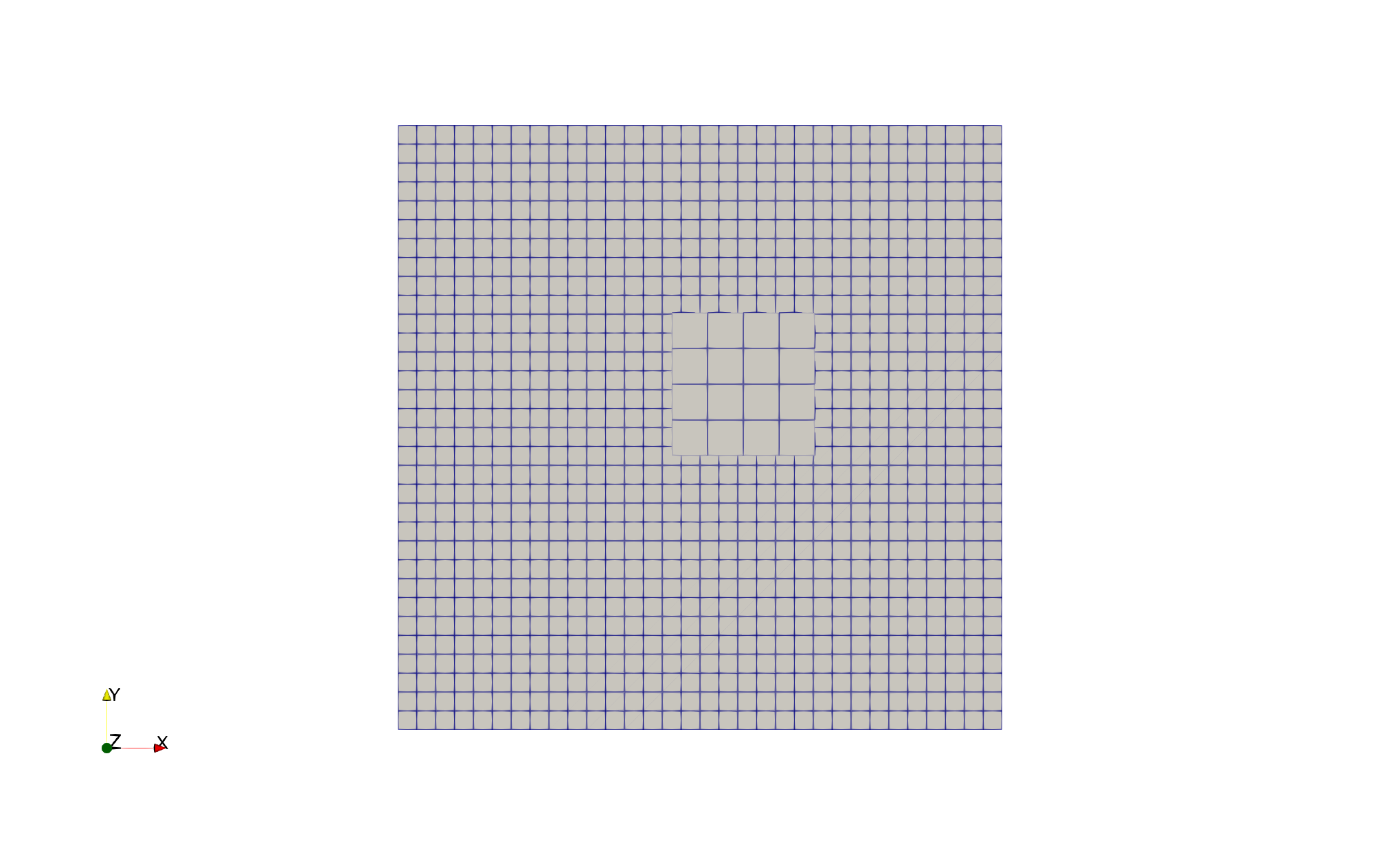} 
\end{center}
\end{minipage}
\begin{minipage}[c]{.5\linewidth}
\begin{center}
		\includegraphics[width=1\linewidth]{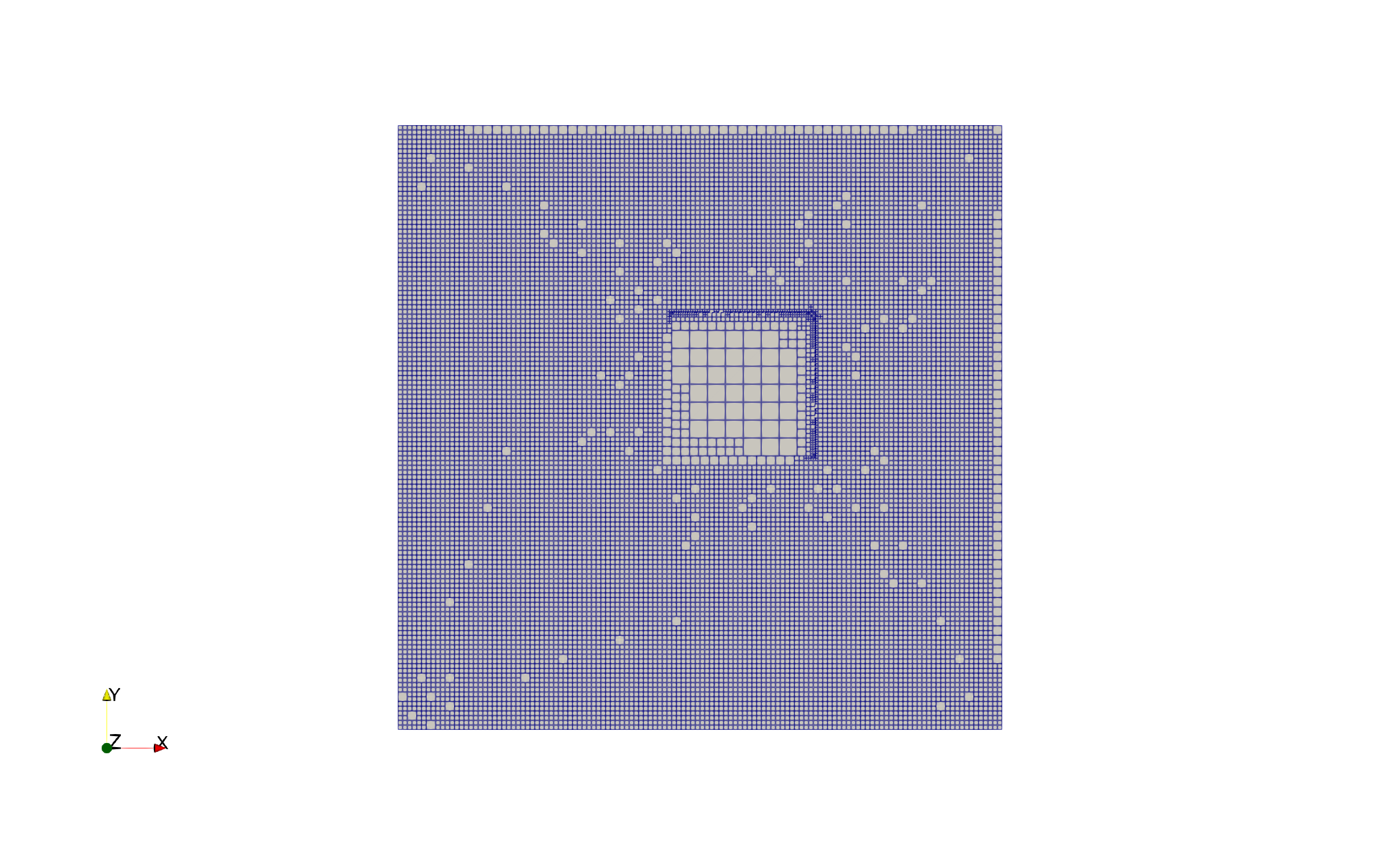}
\end{center}
\end{minipage}
		\caption{Mesh refinements level 0 and Level 5: immersed square shape.}
		\label{fig:mesh_square}
\end{figure}

\begin{figure}[htp]
\begin{minipage}[c]{.5\linewidth}
\begin{center}
		\includegraphics[width=1\linewidth]{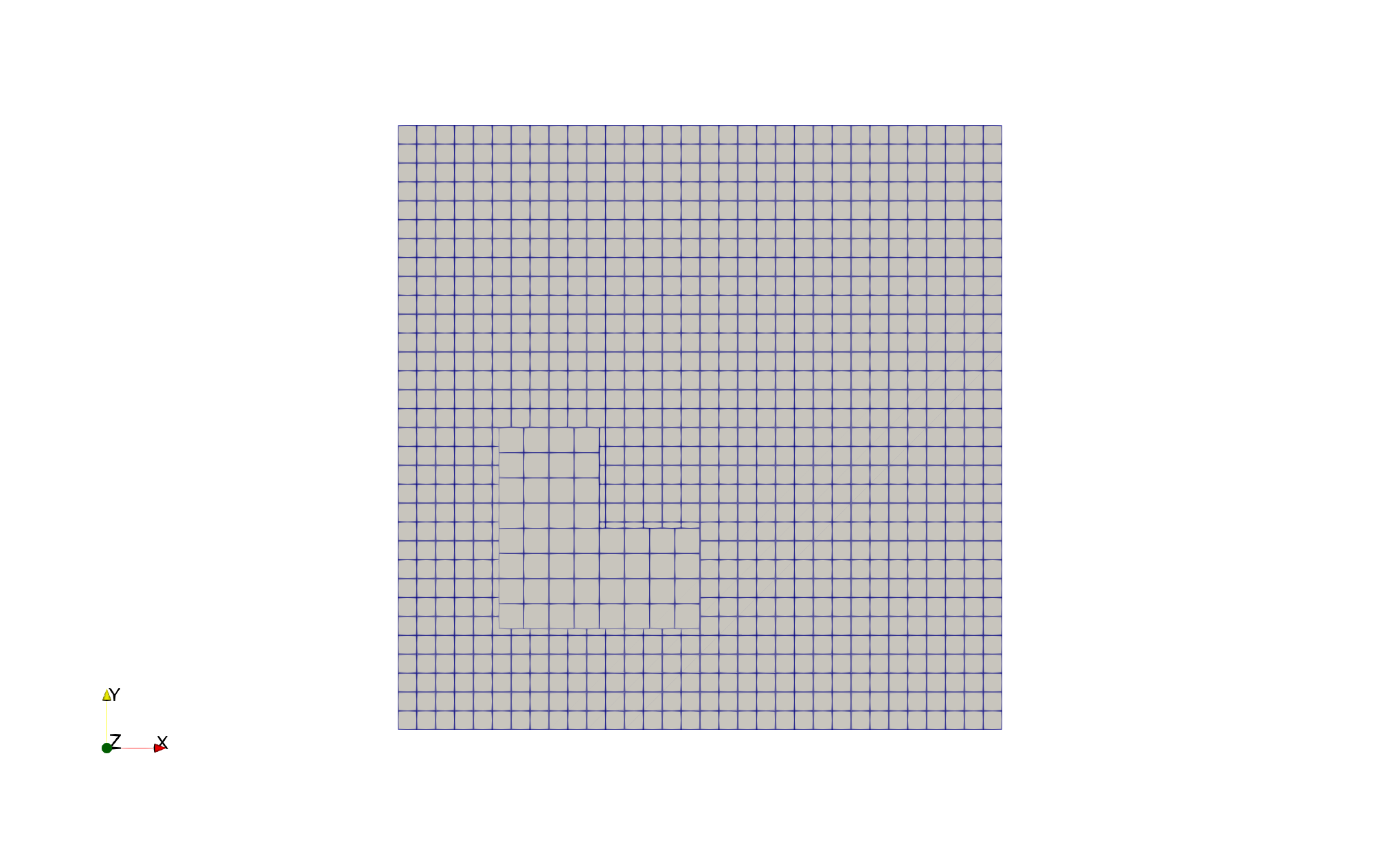} 
\end{center}
\end{minipage}
\begin{minipage}[c]{.5\linewidth}
\begin{center}
		\includegraphics[width=1\linewidth]{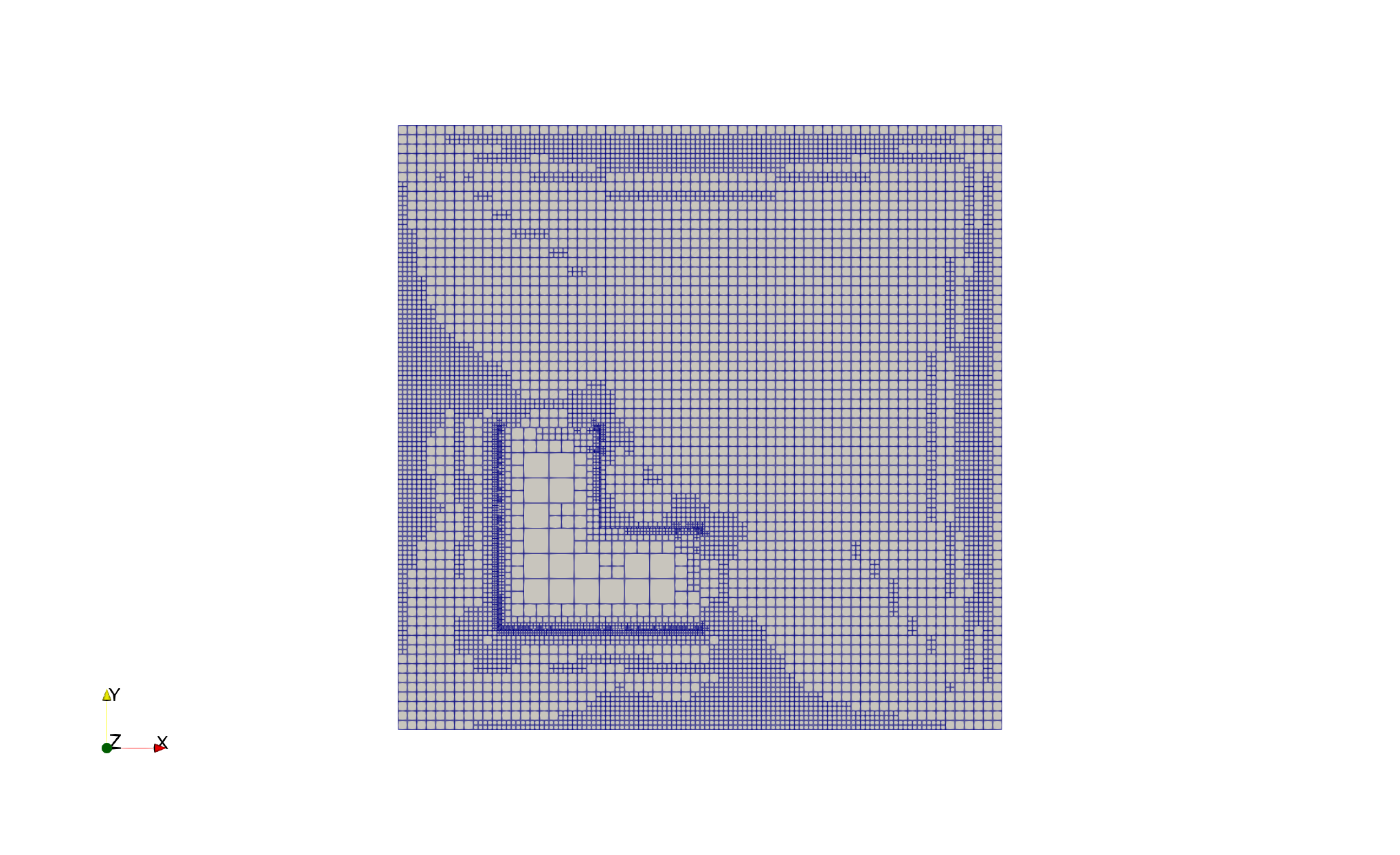}
\end{center}
\end{minipage}
		\caption{Mesh refinements level 0 and Level 5: immersed L-shape.}
		\label{fig:mesh_L}
\end{figure}

\begin{figure}[htp]
\begin{minipage}[c]{.5\linewidth}
\begin{center}
		\includegraphics[width=1\linewidth]{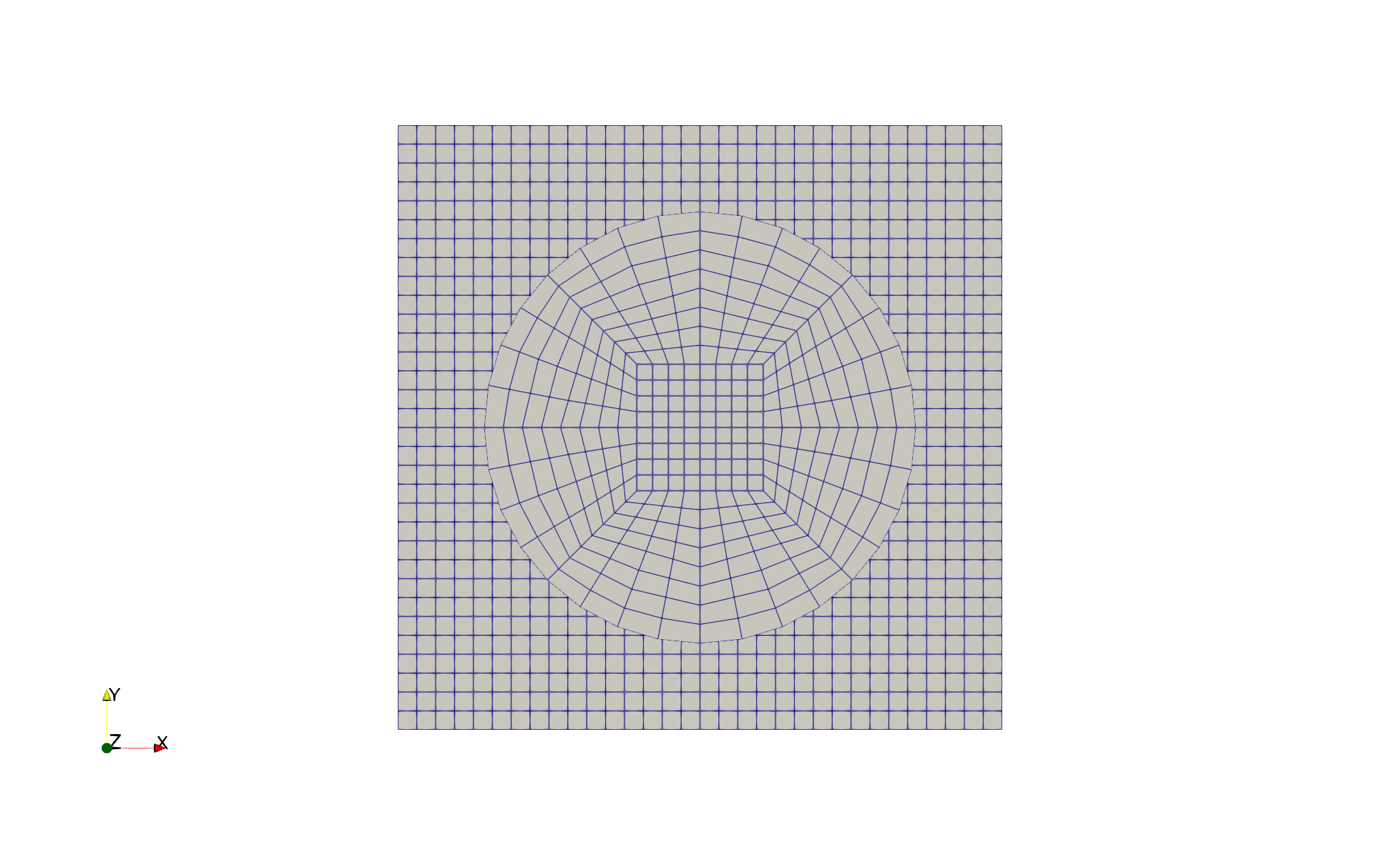}
\end{center}
\end{minipage}
\begin{minipage}[c]{.5\linewidth}
\begin{center}
		\includegraphics[width=1\linewidth]{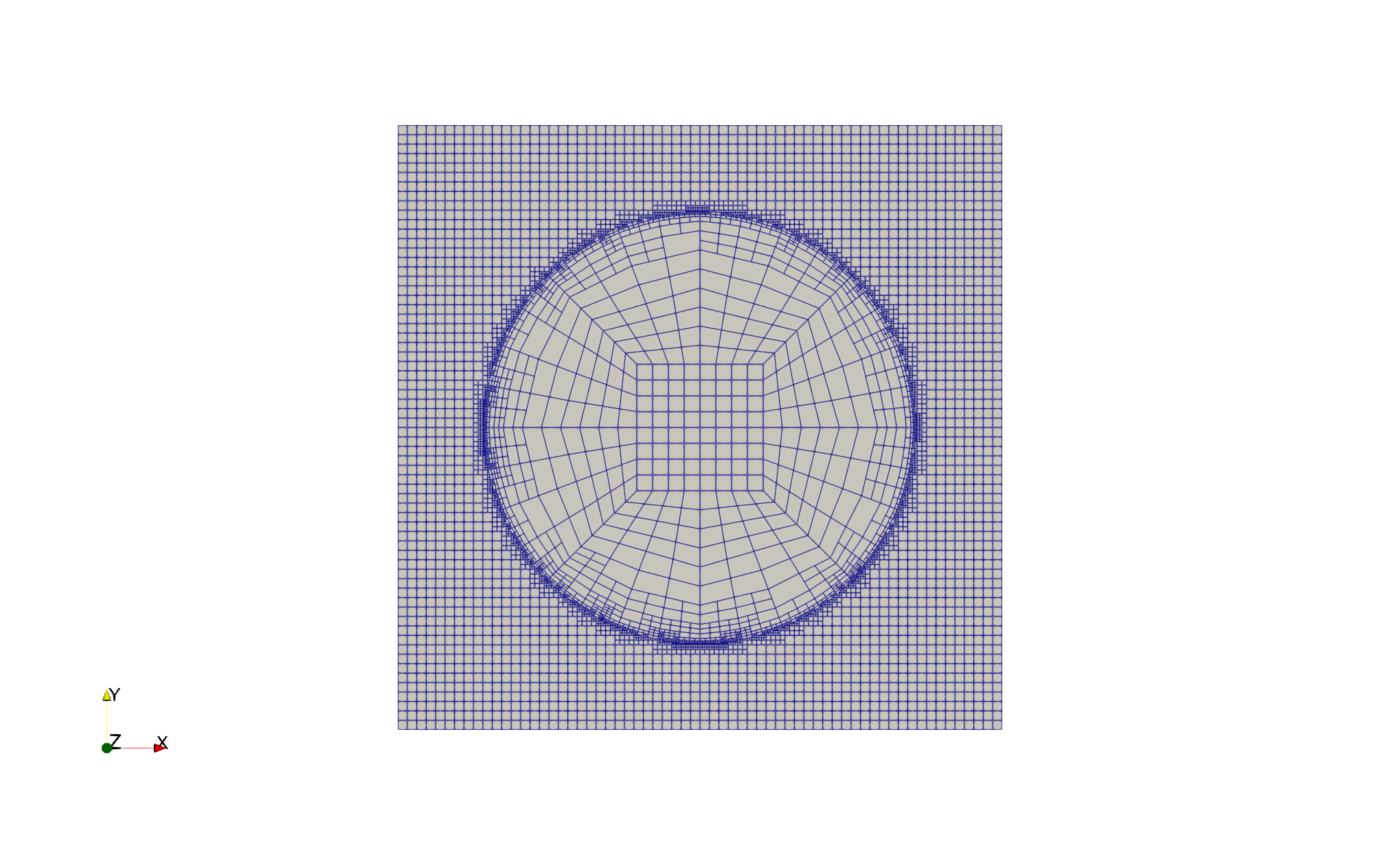}
\end{center}
\end{minipage}
		\caption{Mesh refinements level 0 and Level 5: immersed circle shape.}
		\label{fig:mesh_circle}
\end{figure}

\begin{figure}[htp]
\begin{minipage}[c]{.5\linewidth}
\begin{center}
		\includegraphics[width=1\linewidth]{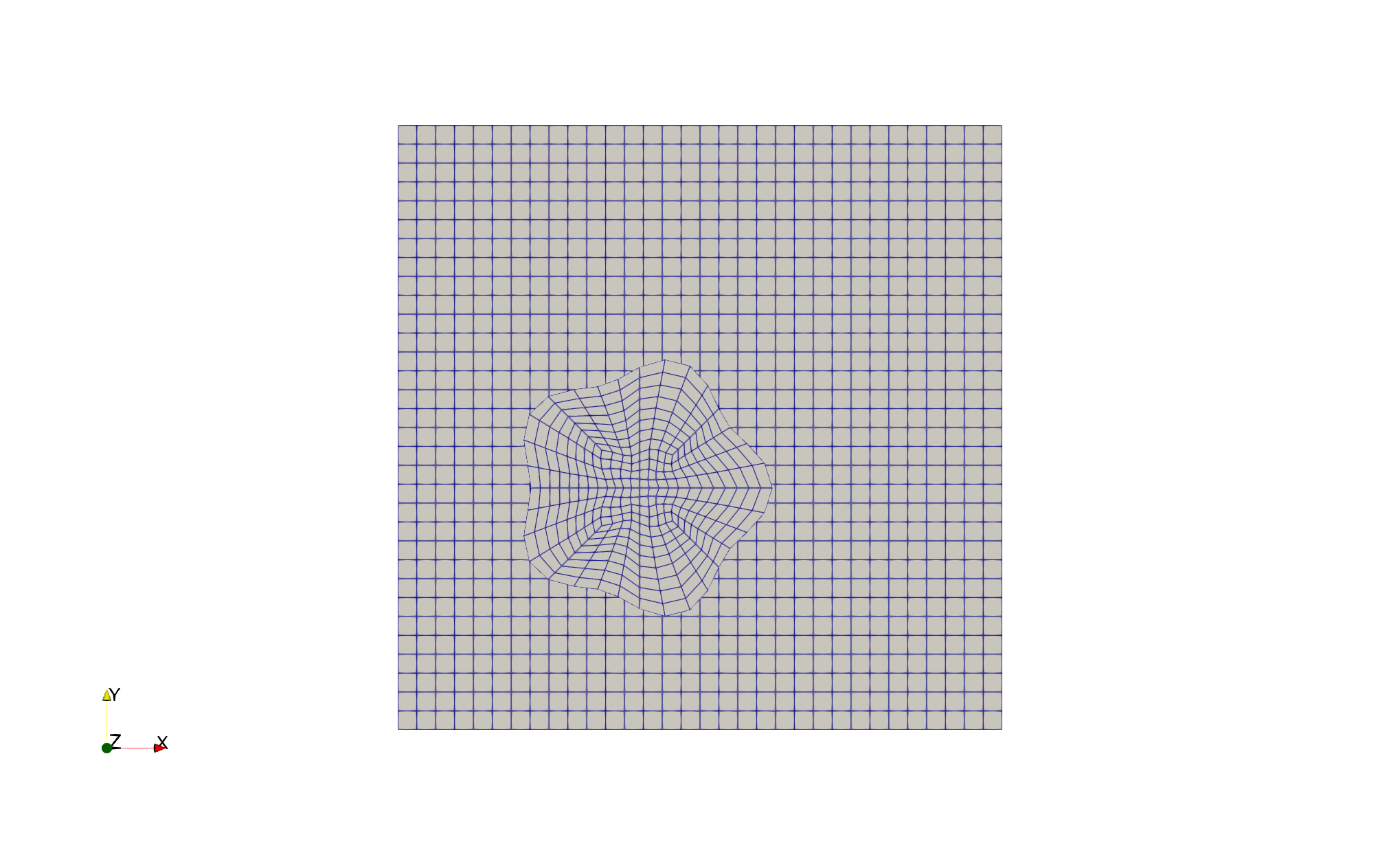} 
\end{center}
\end{minipage}
\begin{minipage}[c]{.5\linewidth}
\begin{center}
		\includegraphics[width=1\linewidth]{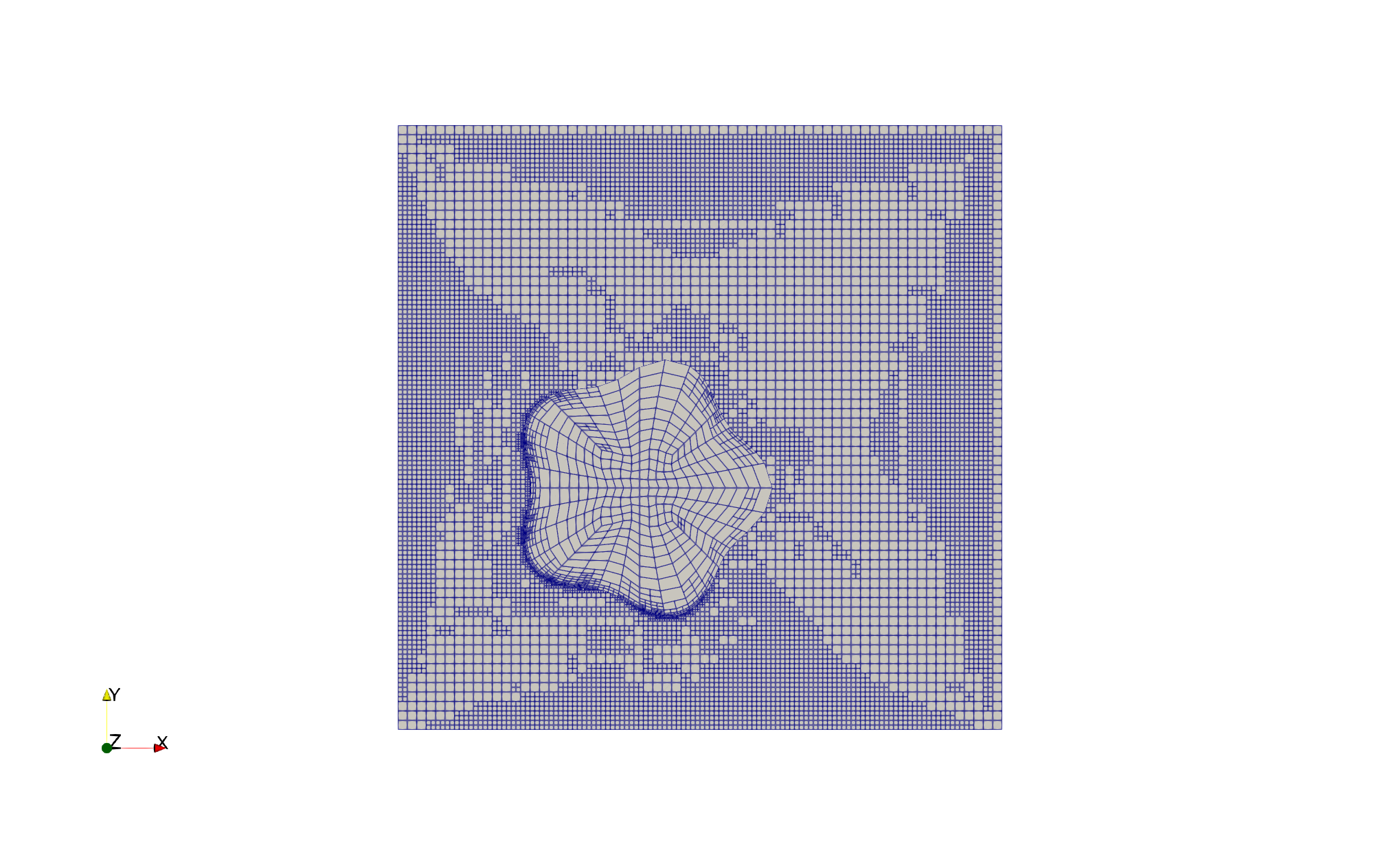}
\end{center}
\end{minipage}
		\caption{Mesh refinements level 0 and Level 5: immersed flower shape.}
		\label{fig:mesh_flower}
\end{figure}

\begin{figure}[htp]
\begin{minipage}[c]{.5\linewidth}
\begin{center}
		\includegraphics[width=1\linewidth]{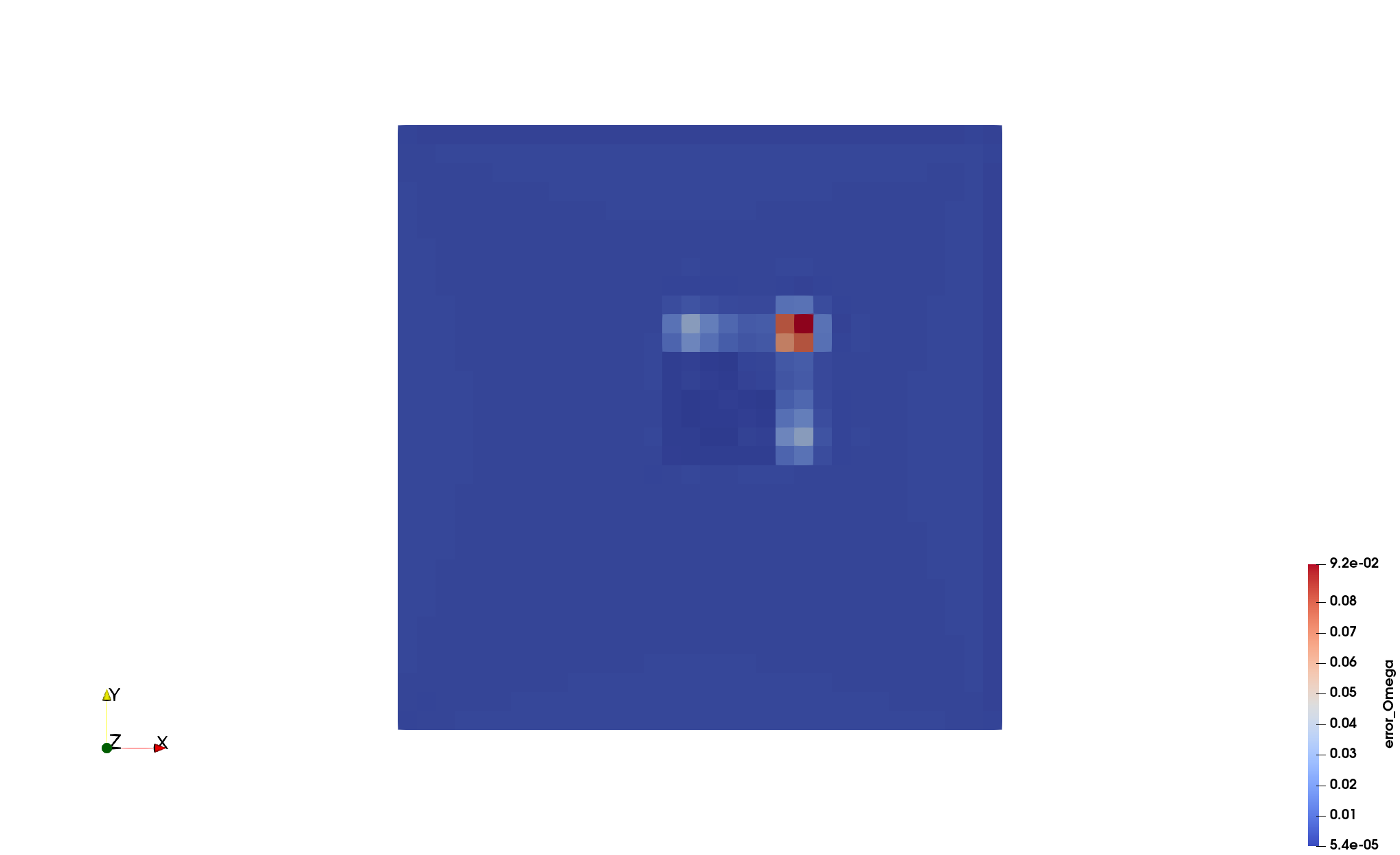} 
\end{center}
\end{minipage}
\begin{minipage}[c]{.5\linewidth}
\begin{center}
		\includegraphics[width=1\linewidth]{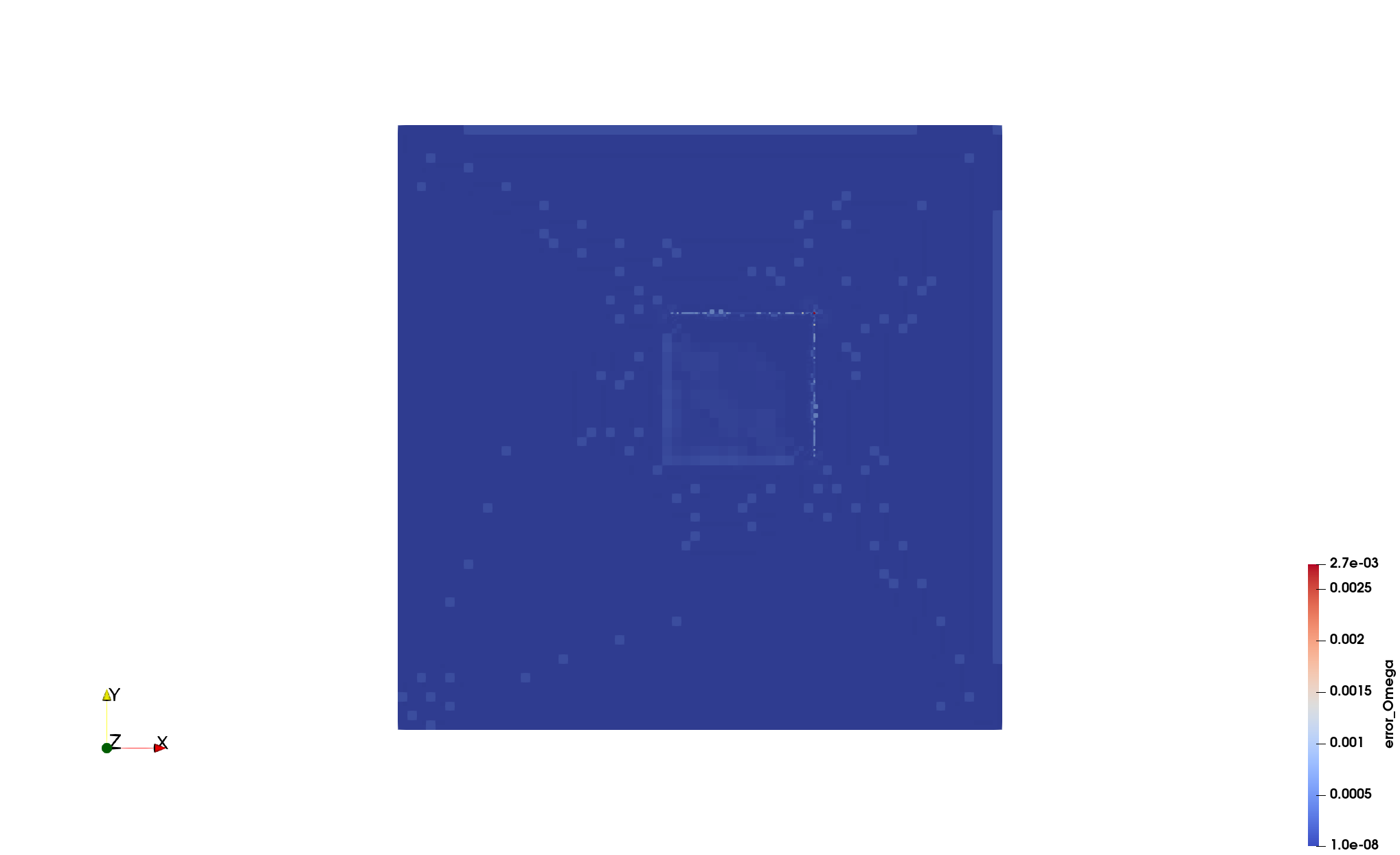}
\end{center}
\end{minipage}
		\caption{Indicators level 0 and level 5: immersed square shape.}
		\label{fig:error_square}
\end{figure}

\begin{figure}[htp]
\begin{minipage}[c]{.5\linewidth}
\begin{center}
		\includegraphics[width=1\linewidth]{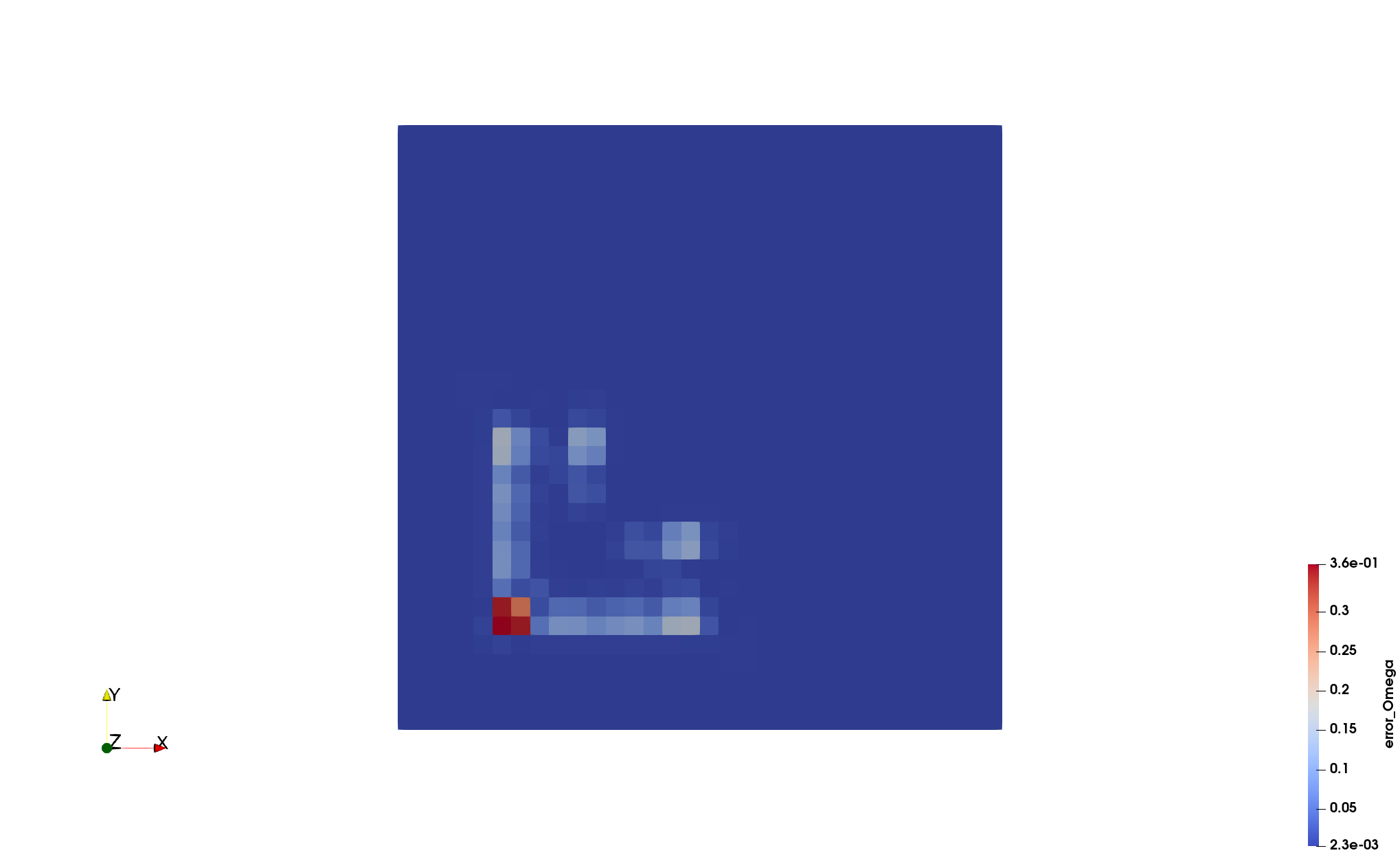} 
\end{center}
\end{minipage}
\begin{minipage}[c]{.5\linewidth}
\begin{center}
		\includegraphics[width=1\linewidth]{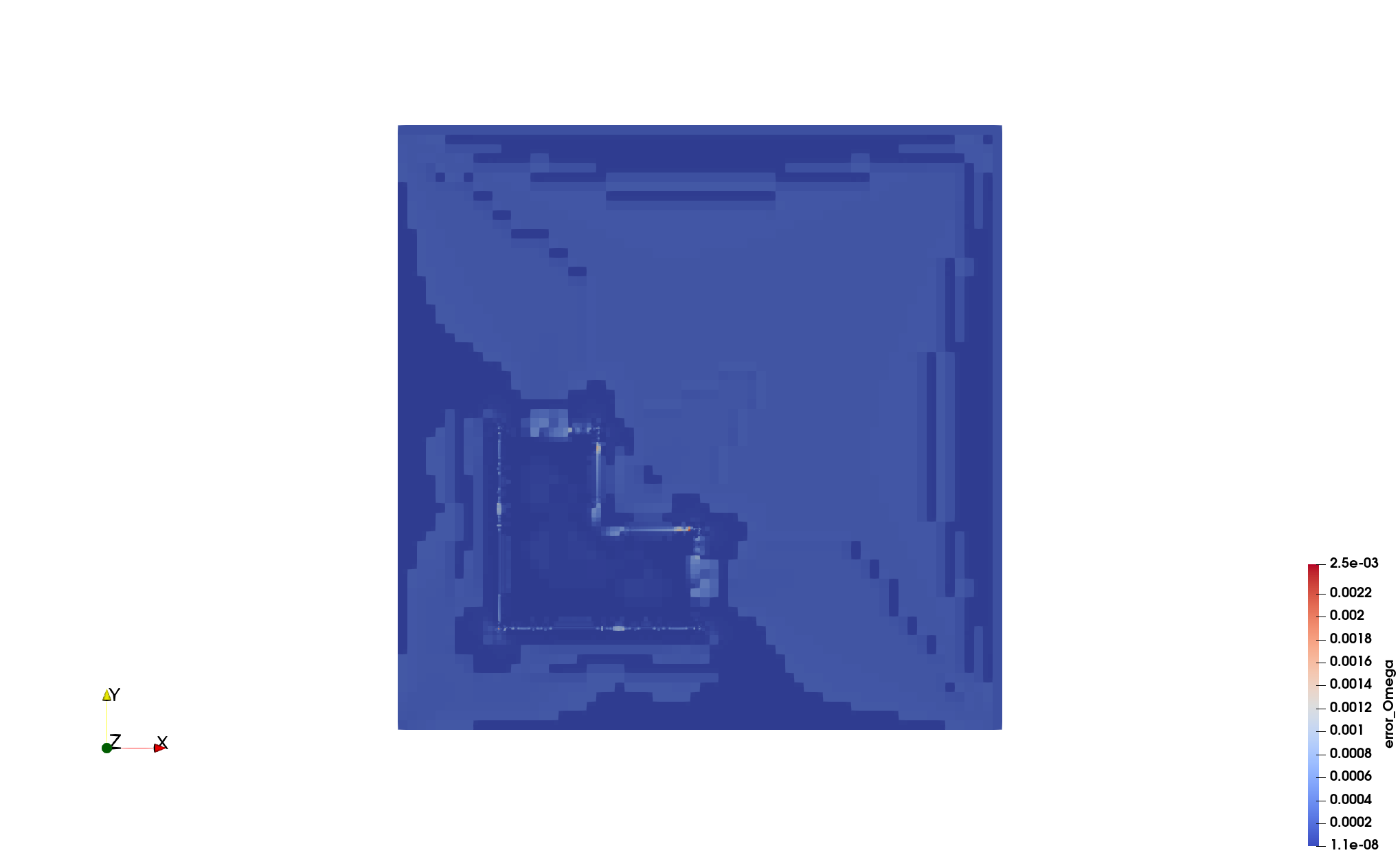}
\end{center}
\end{minipage}
		\caption{Indicators level 0 and level 5: immersed L-shape.}
		\label{fig:error_L}
\end{figure}

\begin{figure}[htp]
\begin{minipage}[c]{.5\linewidth}
\begin{center}
		\includegraphics[width=1\linewidth]{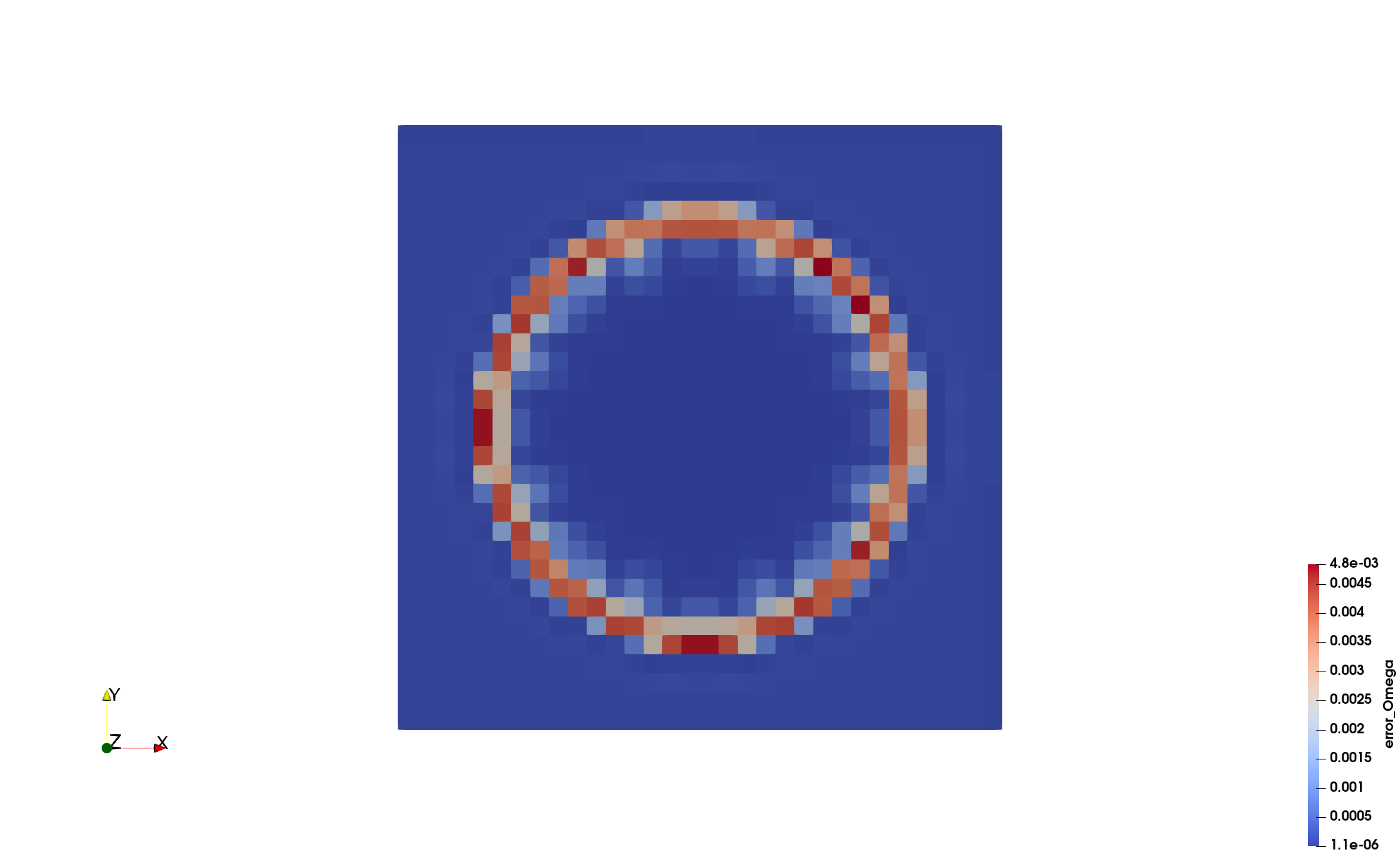} 
\end{center}
\end{minipage}
\begin{minipage}[c]{.5\linewidth}
\begin{center}
		\includegraphics[width=1\linewidth]{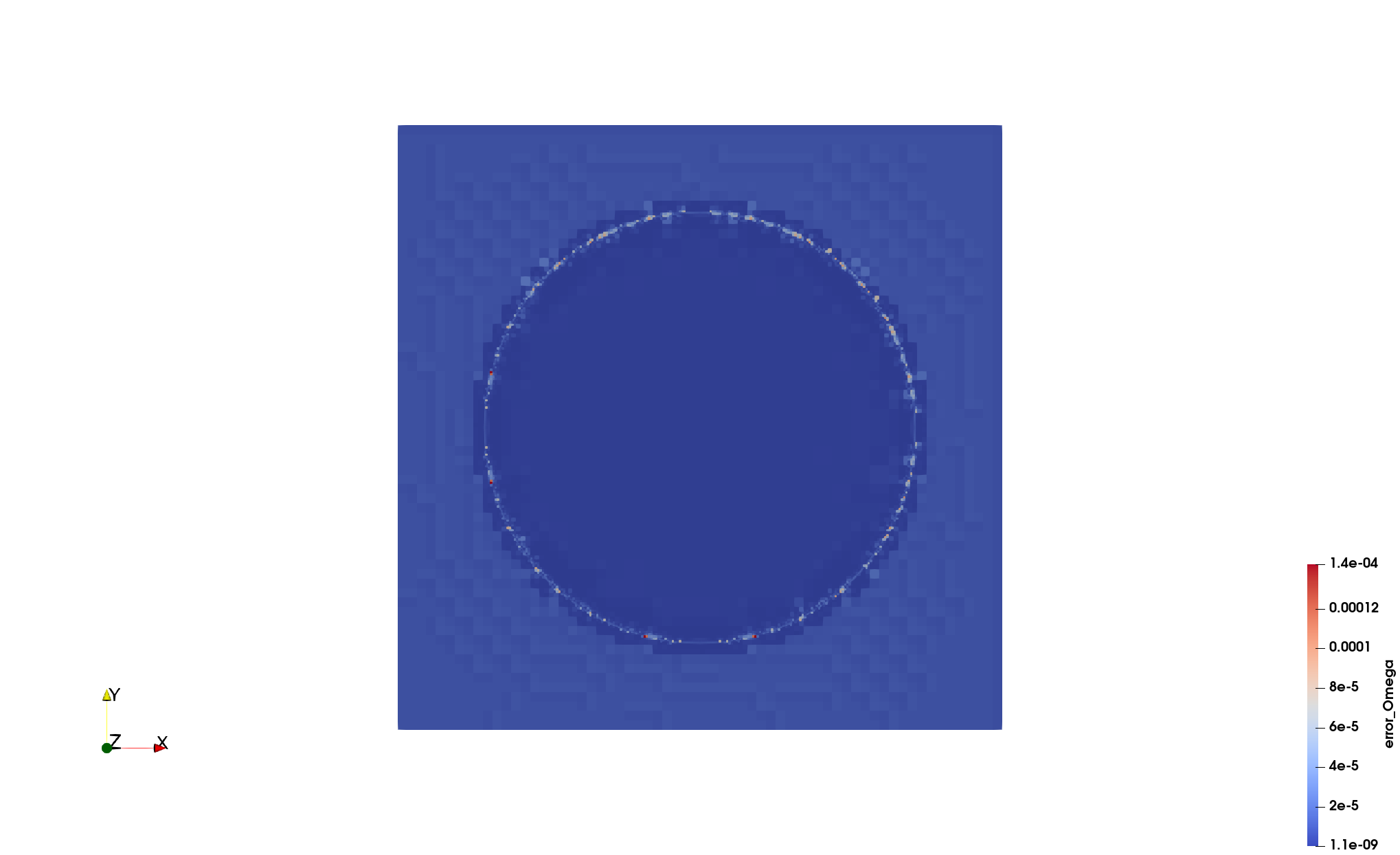}
\end{center}
\end{minipage}
		\caption{Indicators level 0 and level 5: immersed circle shape.}
		\label{fig:error_circle}
\end{figure}

\begin{figure}[htp]
\begin{minipage}[c]{.5\linewidth}
\begin{center}
		\includegraphics[width=1\linewidth]{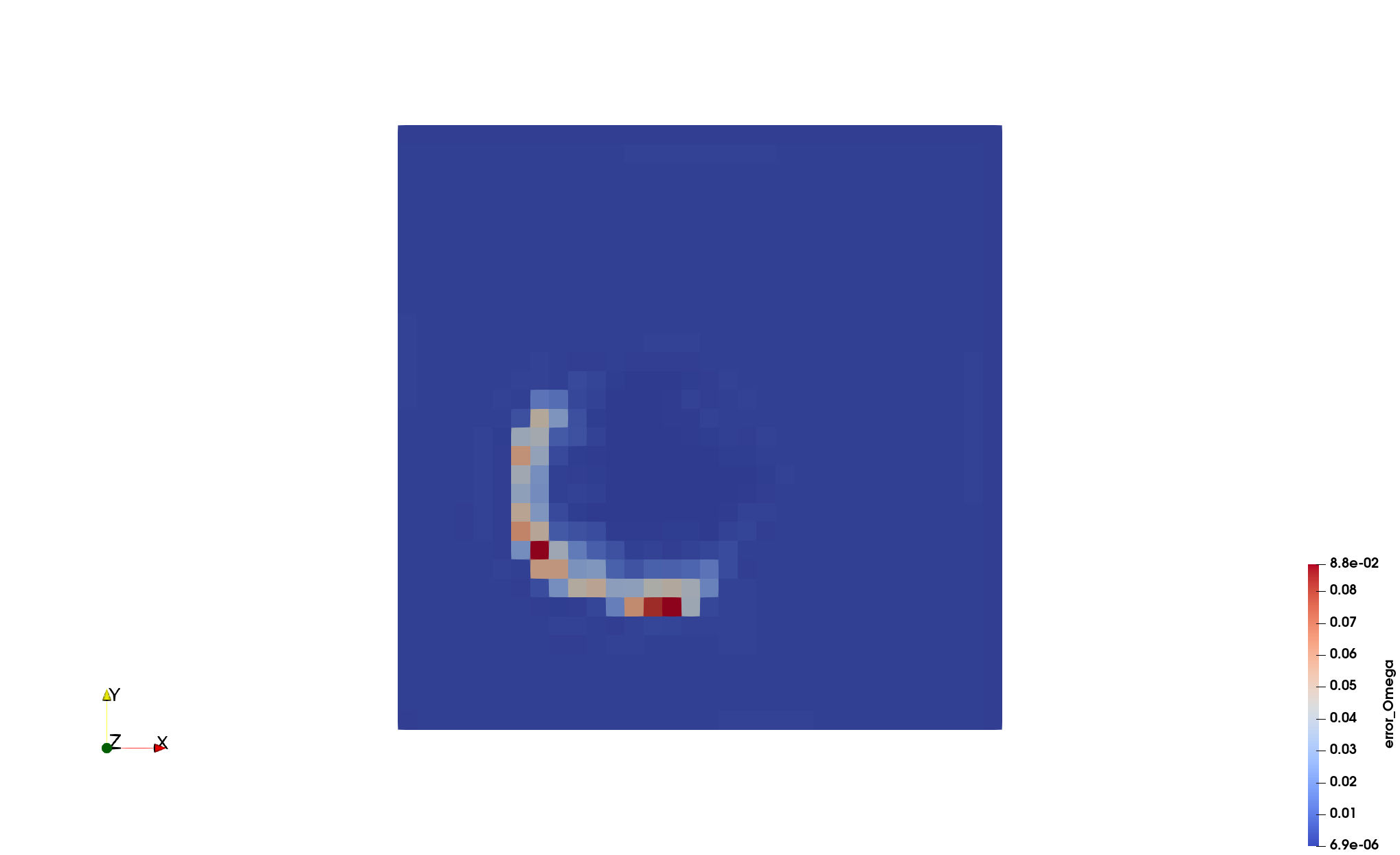} 
\end{center}
\end{minipage}
\begin{minipage}[c]{.5\linewidth}
\begin{center}
		\includegraphics[width=1\linewidth]{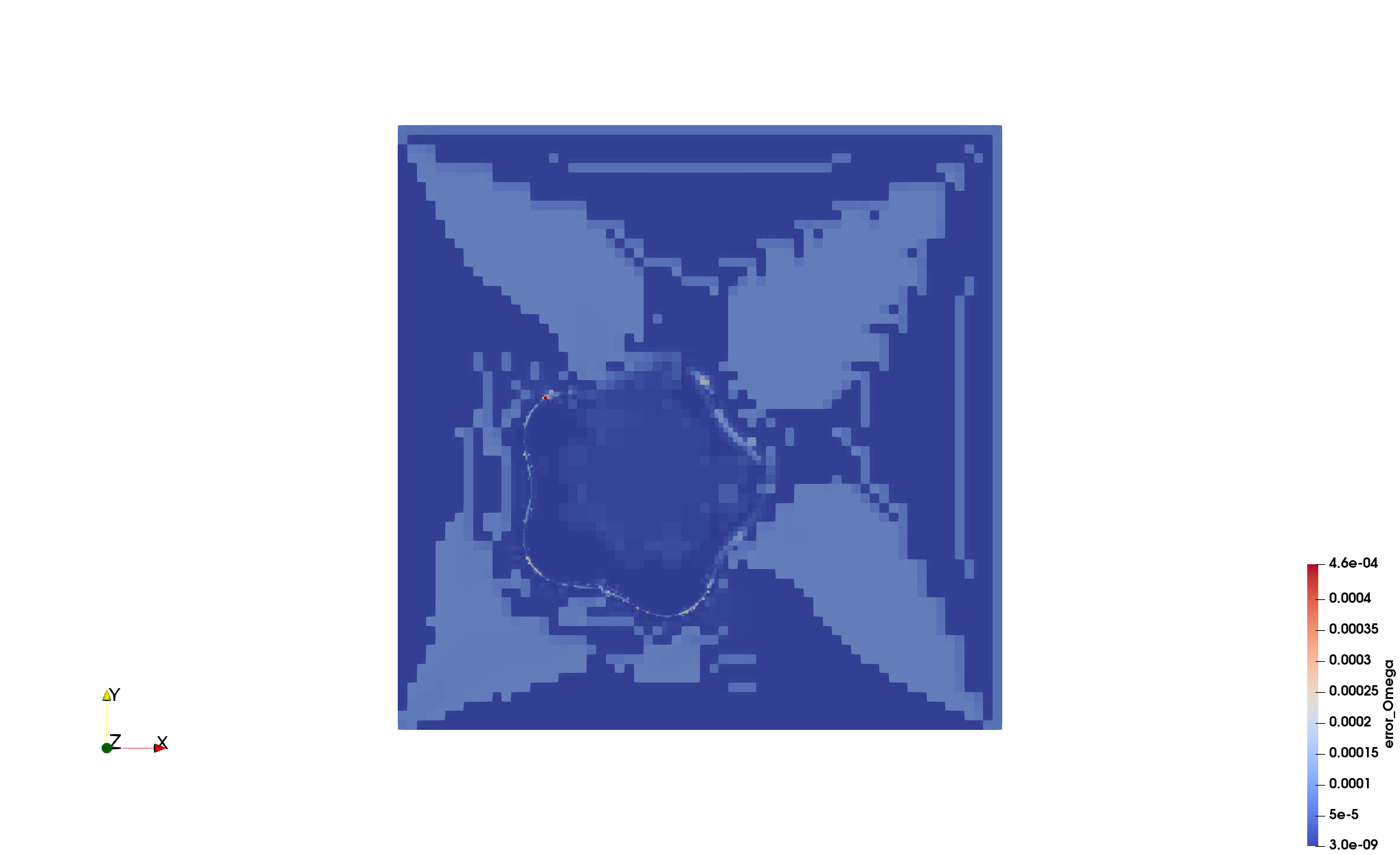}
\end{center}
\end{minipage}
		\caption{Indicators level 0 and level 5: immersed flower shape.}
		\label{fig:error_flower}
\end{figure}

In Figures~\ref{fig:mesh_square}-\ref{fig:mesh_flower}, levels 0 and 5 of refinement of the mesh are presented. On the other hand, the local error indicators are plotted in Figures~\ref{fig:error_square}-\ref{fig:error_flower}. It is evident that the estimator is spotting the regions with higher error, as expected, closer to the interface. Refining those regions results in lowering the local error indicator. 
It is worth highlighting that the edges of the interface nearer to the boundary, where the solution exhibits a larger jump in the normal derivative at the interface, tend to incur higher errors, necessitating additional refinement in those areas. Conversely, in the part of the interface closer to the center, where the transmission of the solution is smoother there, errors are smaller, resulting in a reduced need for refinement.

Furthermore, when there is a substantial discrepancy in coefficient jumps, the areas with higher coefficients often yield solutions that approach constancy. This is because we are applying a much smaller force that inversely correlates with the coefficients values. As the value of the jump of the coefficients increase, we approach a scenario where the problem is effectively solved with negligible applied force on the region with the higher coefficient.

In conclusion, the numerical findings in this section corroborate the theoretical results, demonstrating that our a posteriori error estimator is reliable and efficient. Additionally, they showcase the robustness of our proposed approach across various immersed geometries and in the presence of strong singularities in the solution.





\bibliographystyle{apalike}
\bibliography{References}

\end{document}